\documentclass{article}

\usepackage[width=6in]{geometry}
\usepackage{authblk}
\usepackage{amssymb,amsthm,amsmath}
\usepackage{threeparttable}
\usepackage{tabularx} 
\usepackage{booktabs}
\usepackage{multirow}
\usepackage[ruled]{algorithm2e} 
\usepackage{stmaryrd}
\usepackage{bbm}
\usepackage{subfig}
\usepackage{graphicx}
\usepackage{color}
\usepackage{hyperref}
\title{Solving Multi-Group Neutron Diffusion Eigenvalue Problem with Decoupling Residual Loss Function}
\author[1]{Shupei Yu}
\author[1]{Qiaolin He}
\author[1]{Shiquan Zhang}
\author[1]{Qihong Yang}
\author[1]{Yu Yang}
\author[2]{Helin Gong}
\affil[1]{School of Mathematics, Sichuan University, Chengdu 610065, China.}
\affil[2]{Paris Elite Institute of Technology, Shanghai Jiao Tong University, Shanghai 200240, China.}

\begin{document}
	\maketitle
	
	\begin{abstract}
		In the midst of the neural network's success in solving partial differential equations, tackling eigenvalue problems using neural networks remains a challenging task. However, the Physics Constrained-General Inverse Power Method Neural Network (PC-GIPMNN) approach was proposed and successfully applied to solve the single-group critical problems in reactor physics. This paper aims to solve critical problems in multi-group scenarios and in more complex geometries. Hence, inspired by the merits of traditional source iterative method, which can overcome the ill-condition of the right side of the equations effectively and solve the multi-group problem effectively, we propose two residual loss function called Decoupling Residual loss function and Direct Iterative loss function. Our loss function can deal with multi-group eigenvalue problem, and also single-group eigenvalue problem. Using the new residual loss functions, our study solves one-dimensional, two-dimensional, and three-dimensional multi-group problems in nuclear reactor physics without prior data. In numerical experiments, our approach demonstrates superior generalization capabilities compared to previous work.
	\end{abstract}
	\begin{center}
		\textbf{keywords:} deep learning, eigenvalue problem, nuclear reactors, multi-group problem
	\end{center}
	
	\section{Introduction}\label{Introduction}
	The neutron diffusion equations are fundamental equations in nuclear reactor physics, which are used to describe the transport behavior of neutrons in a nuclear reactor and derived by neutron diffusion theory \cite{osti_4074688}. Depending on different physical situations and assumptions, the neutron diffusion equations can be further developed into more complex equations, such as the multi-group diffusion equations, neutron noise equations, and so on.  
	{The neutron diffusion equation can be simplified from the Boltzmann transport equation, which accurately describes the neutron transport process, and the neutron diffusion equation includes single-group, multi-group, transient and steady-state problems.}
	The neutron diffusion equations can be used for fuel management to determine the optimal arrangement of fuel rods, thereby maximizing fuel utilization. The neutron diffusion equations can also be used for reactor operation control to maintain the stability and safety performance of the reactor, such as controlling the power of the reactor by controlling the neutron flux in the reactor.
	
	Over the past few decades, researchers have developed a range of solution methods for problems such as finite difference\cite{hamada2022higher}, finite element\cite{azekura1980new}, finite volume\cite{demaziereCORESIMMultipurpose2011} \cite{bernal2014resolution} ,nodal expansion\cite{osti_10191160} and methods of characteristic\cite{tao2022neutron}, among others. However, with the potential of neural networks being explored in various fields, there is an urgent need for research in using neural networks to solve physically relevant multi-group neutron diffusion eigenvalue problems (NDEPs) in nuclear reactor scenarios.
	
	As early as 1994, Dissanayake et al.\cite{dissanayake1994neural} attempted to use neural network methods to solve simple cases of linear and nonlinear problems. In 1998, Lagris et al.\cite{lagaris1998artificial} presented a more comprehensive algorithmic framework for using neural networks to solve partial differential equations. Numerous researchers have made continuous efforts in the development of related works. 
	Currently, the most popular framework is the Physics-Informed Neural Networks (PINNs) proposed by Raissi et al.\cite{raissi2019physics} This framework directly links neural networks with the information of physical equations through the form of a loss function.
	It is evident that neural network based methods for solving partial differential equations offer several advantages:
	\begin{itemize}
		\item Independence from mesh generation: Neural networks do not rely on mesh files for solving PDEs. Instead, they utilize collected training samples or data points, which eliminates the need for grid generation and adapts well to irregular or complex geometries.
		\item Integration of observational data: Neural networks have the capability to incorporate observational or experimental data into the learning process. This allows them to effectively combine the physical information from the equations with the available data, resulting in enhanced accuracy and predictive capabilities compared to traditional solver algorithms that solely rely on the physics of the equations.
		\item Handling high-dimensional problems: Neural networks demonstrate their advantage in dealing with high-dimensional problems, overcoming the curse of dimensionality. Neural networks can effectively learn and represent complex relationships in high-dimensional spaces, making them suitable for a wide range of problems with multiple variables or parameters.
	\end{itemize}
	
	{The PINNs method has also been applied in inverse problems \cite{raissi2020hidden} \cite{chen2020physics}, as well as in numerical solutions for solving stochastic differential equations \cite{doi:10.1137/19M1260141} \cite{yang2020physics}. For interface problems that are more commonly encountered in practical applications, there are also related studies \cite{wuINNInterfacedNeural2022} \cite{hu2022discontinuity} \cite{hu2022shallow}.}
	
	{Recently, neural network algorithms have been gradually maturing, and they have been instrumental in solving engineering problems with increased robustness and efficiency.} {Cheng et al.\cite{cheng2023generalised} \cite{cheng2024efficient} simplified the process of data assimilation (DA) using neural networks.} {Gong et al.\cite{gong2022data} utilized the data-enabled physics-informed approach to solve the neutron field and coefficients in nuclear reactor physics.}
	{Phillips et al.\cite{phillips2021} replaced the discrete forms in traditional numerical methods with convolutional kernels from CNNs for numerical computation, yielding satisfactory results.} 
	
	{Based on the aforementioned advantages, there also have been numerous works applying neural network algorithms to solve PDEs in various problem domains. Researchers utilize Physics-informed PointNet \cite{kashefi2023prediction} to predict the velocity and pressure fields of two-dimensional steady incompressible flow in porous media. Zhou et al.\cite{zhou2024multi} proposed the Multi-Scale Physics Constrained Neural Network (MSPCNN) to solve fluid dynamics problems by integrating fidelity terms at multiple scales. Mao et al.\cite{mao2020physics} utilized the PINN framework with the integration of prior data point losses to solve one-dimensional and two-dimensional forward and inverse problems in high-speed flow.}
	{These examples highlight the strong capability of neural network algorithms to integrate physical information with data points, which is crucial for solving real-world physical problems.}
	
	Specifically, there is still considerable research being conducted on PDEs in the field of nuclear reactor physics, which is also the focus of this paper. We have listed the relevant articles in Table \ref{tab-intro}.
	
	The eigenvalue problem in one dimension was solved for a single material by Wang et al\cite{wang2022surrogate}. Yang et al.\cite{yangDataenabledPhysicsinformedNeural2022} utilized the Data Enabled Physics Informed Neural Network (DEPINN) approach, where observed data points are incorporated into the neural network for solving the problem. This approach eliminates the need for regularization techniques. On the other hand, Yang et al.\cite{yang2023neural} \cite{yang2023physics} focused on solving eigenvalue problems solely based on the information provided by the physical equations, without considering any prior data points:
	\begin{itemize}
		\item The main emphasis of Generalized Inverse Power Method is to combine neural networks with the traditional power iteration method. It represents a practical implementation of neural networks within the computational framework of traditional eigenvalue problems.
		\item Physics Constrained refers to the application of neural networks to enforce the continuity of energy groups and flux, which effectively enhances the fidelity of the physical modeling itself.
		\item Based on the definition of eigenvalues and in combination with control equations, this paper adopts Rayleigh-Quotient expression to accelerate the convergence of eigenvalues. The objective is to improve the convergence rate of the eigenvalue solution.
	\end{itemize}
	
	However, there have been no reported studies utilizing neural networks to solve the multi-group eigenvalue problem without the introduction of prior data points. The multi-group problem is particularly relevant in nuclear reactor physics, where there is a strong engineering demand for its solution. Therefore, the main contributions of this paper are as follows:
	\begin{itemize}
		\item This paper introduces a novel loss function called decoupling loss function that effectively handles the ill-conditioned structure of the governing equations when solving multi-group neutron diffusion problems. It is designed to be applicable even in the case of degenerating into a single-group scenario.
		\item In this paper, the decoupling loss function is applied to multidimensional (1D, 2D, and 3D) multi-group complex interface problems. This marks the first application of neural network-based solutions to multi-group neutron diffusion problems in real 3D physical scenarios of nuclear reactors.
		\item In this paper, detailed numerical experiments are conducted for different reactor geometries, including sampling from training points and applying interface conditions. These experiments provide a solid foundation for future research in this field.
	\end{itemize}
	
	These contributions highlight the advancements made in applying neural networks to solve the multi-group eigenvalue problem, shedding light on the challenges and potential solutions in the field of nuclear reactor physics.
	
	The subsequent content of this paper is as follows.	In Section \ref{Problem Statement}, we provide an overview of problems being solved. In Section \ref{Neural Networks}, we delve into the details of the neural network algorithm employed in this study. We presents the numerical results obtained from applying the neural network algorithm to specific problems in Section \ref{Numerical Experiments}. Finally, we draw conclusions in Section \ref{Conclusion}.
	
	\begin{table}[!htb]
		\caption{Main researches in neural network algorithms solving nuclear reactor physics PDE problem. \label{tab-intro}}
		\centering
		\begin{threeparttable}
			\begin{tabularx}{.815\textwidth}{|c|c|c|c|c|c|c|}
				\cline{1-7}
				Method & NoEG & NoM & Dim & Eig & TD & Model \\
				\cline{1-7}
				\multirow{2}{*}{PIDL\cite{xie2021neural}} & 1 & 1 & 2 & - & $\circ$ & Diffusion\\
				\cline{2-7}
				& 1 & 2 & 2 & - & - & Diffusion\\
				\cline{1-7}
				\multirow{2}{*}{SVD-AE\cite{phillips2021autoencoder}} & 1 & 2 & 1 & $\circ$ & - & Diffusion\\
				\cline{2-7}
				& 1 & 4 & 2 & $\circ$ & - & Diffusion\\
				\cline{1-7}
				\multirow{3}{*}{BDPINN\cite{xie2024boundary}} & 1 & 2 & 1 & - & - & Transport\\
				\cline{2-7}
				& 1 & 2 & 2 & - & - & Transport\\
				\cline{2-7}
				& 2 & 3 & 2 & - & - & Transport\\
				\cline{1-7}
				\multirow{2}{*}{BC-cPINN\cite{wang2022surrogate}} & 1 & 1/2/5 & 1 & - & - & Diffusion\\
				\cline{2-7}
				& 1 & 1/2 & 1 & $\circ$ & - & Diffusion\\
				\cline{1-7}
				\multirow{2}{*}{DEPINN\cite{yangDataenabledPhysicsinformedNeural2022}}& 1 & 1 & 1/2 & $\circ$ & - & Diffusion\\
				\cline{2-7}
				& 2 & 4 & 2 & $\circ$ & - & Diffusion\\
				\cline{1-7}
				\multirow{2}{*}{PC-GIPMNN\cite{yang2023physics}} & 1 & 3 & 1 & $\circ$ & - & Diffusion \\
				\cline{2-7}
				& 1 & 4/6 & 2 & $\circ$ & - & Diffusion \\
				\cline{1-7}
				\multirow{5}{*}{PINNs\cite{doi:10.1080/00295639.2022.2123211}} & 1 & 2 & 2 & - & - & Diffusion \\
				\cline{2-7}
				& 1 & 2/3 & 2 & $\circ$ & - & Diffusion \\
				\cline{2-7}
				& 2 & 1/3 & 1 & - & - & Diffusion \\
				\cline{2-7}
				& 2 & 2 & 2 & - & - & Diffusion \\
				\cline{2-7}
				& 2 & 3 & 1/2 & $\circ$ & - & Diffusion \\
				\cline{1-7}
				{\textbf{Present work}} & 2 & 2/3/4/5 & 1/2/3 & $\circ$ & - & Diffusion\\
				\cline{1-7}
			\end{tabularx}
			\begin{tablenotes}
				\item[*] `NoEG', `NoM', `Dim', `Eig' and `TD' mean number of energy group, number of material, dimension, eigenvalue and time-dependent, respectively. 
				\item[*] In the table, we use `$\circ$' and `-' to indicate whether the model covers a particular scenario. `$\circ$' represents `yes' and indicates that the model covers the scenario, while `-' represents `no' and indicates that the model does not cover the scenario. 
				\item[*] Among the aforementioned work, the problems being solved in this paper are based on realistic models of nuclear reactor cores \cite{theler2011solution} and are solved without introducing any observation points.
			\end{tablenotes}
		\end{threeparttable}
	\end{table}
	
	\section{Problem Statement} \label{Problem Statement}
	First, let us consider a non time-dependent two-group diffusion problem where the physical equations include the fast group $\phi_1(\boldsymbol{x})$ and thermal group $\phi_2(\boldsymbol{x})$ in \eqref{eq-diff-eigen}, which is defined in domain $\Omega$. {The diffusion of neutrons in two groups can be described by coupled diffusion equations that account for processes such as neutron sources, absorption, scattering, and moderation. The left-hand side of the equations considers the diffusion and reaction processes for both fast and thermal neutron groups, with neutrons from the fast group transitioning into the thermal group. This system of equations assumes that fission neutrons only fall within the fast group range; hence, the source term in the thermal group equation originates solely from the moderation process of the fast group.}
	\begin{equation}
		\left\{
		\begin{aligned}
			-\nabla \cdot D_1 \nabla \phi_1(\boldsymbol{x}) + (\Sigma_{a,1}+\Sigma_{1\rightarrow 2})\phi_1(\boldsymbol{x}) & = \frac{1}{k}(\nu\Sigma_{f,1}\phi_1(\boldsymbol{x})+\nu\Sigma_{f,2}\phi_2(\boldsymbol{x})), & &~ \text{in} ~ \Omega,\\
			-\nabla \cdot D_2 \nabla \phi_2(\boldsymbol{x}) + \Sigma_{a,2}\phi_2(\boldsymbol{x}) - \Sigma_{1\rightarrow 2}\phi_1(\boldsymbol{x}) & = 0, & &~ \text{in} ~\Omega,
		\end{aligned}
		\right.
		\label{eq-diff-eigen}
	\end{equation}
	
	\begin{equation}
		\left\{
		\begin{aligned}
			\frac{\partial\phi_e}{\partial \boldsymbol{n}} & = 0,  &~ \text{on} ~\partial\Omega_{sym},\\
			\phi_e & = f,  &~ \text{on} ~\partial\Omega_{ext,1},\\
			\frac{\partial\phi_e}{\partial \boldsymbol{n}} & = -\frac{c_{bou}}{D_e}\phi_e,  &~ \text{on} ~\partial\Omega_{ext,2},
		\end{aligned}
		\right.
		\label{eq-diff-bou}
	\end{equation}
	where $D_i, \Sigma_{a,i}, \Sigma_{f,i}$ represent diffusion coefficient, macroscopic absorption cross section and fission cross section, respectively. $\Sigma_{1\rightarrow 2}$ stands for scattering cross section from fast group to thermal group. {$\frac{1}{k}$ and $\phi_1, \phi_2$ denote the eigenvalue and eigenfunctions, respectively.}
	The boundary conditions for the problem depend on the specific physical problem at hand. However, we can summarize them into \eqref{eq-diff-bou}, where $e = 1,2$ represent different energy groups, $\boldsymbol{n}$ represents the normal direction, $c_{bou}$ is a constant coefficient for the Robin condition, $\partial\Omega_{sym}$ indicates the symmetric boundary where the homogeneous Neumann condition is applied, while $\partial\Omega_{ext,1}$ and $\partial\Omega_{ext,2}$ represent different types of external boundaries, where $\partial\Omega_{ext,1}$ applies the Dirichlet boundary condition and $\partial\Omega_{ext,2}$ applies the Robin boundary condition. For different test cases, the applied boundary conditions vary. The specific details will be discussed in Section \ref{Numerical Experiments}.
	
	From mathematical perspectives, \eqref{eq-diff-eigen} is a system of eigenvalue partial differential equations, defined on the domain $\Omega$. 
	{The number of solutions to such problems are infinite, but for such a critical system of neutron diffusion problems, the primary concern in nuclear reactors is to determine the maximum value of $k$ and the corresponding values of $\phi_1$ and $\phi_2$ when $k$ reaches its maximum. This value of $k$ is referred to as $k_{\text{eff}}$, the effective multiplication factor. It can be observed that obtaining the maximum value of $k$ and the corresponding $\phi_1$ and $\phi_2$ is equivalent to solving for the smallest eigenvalue and eigenvector of equation \eqref{eq-diff-eigen}. Since this paper only considers the solution of the effective multiplication factor, the term $k_{\text{eff}}$ will be used to denote $k$ hereinafter.}
	
	\begin{figure*}[h]
		\centering
		\includegraphics[width=15cm]{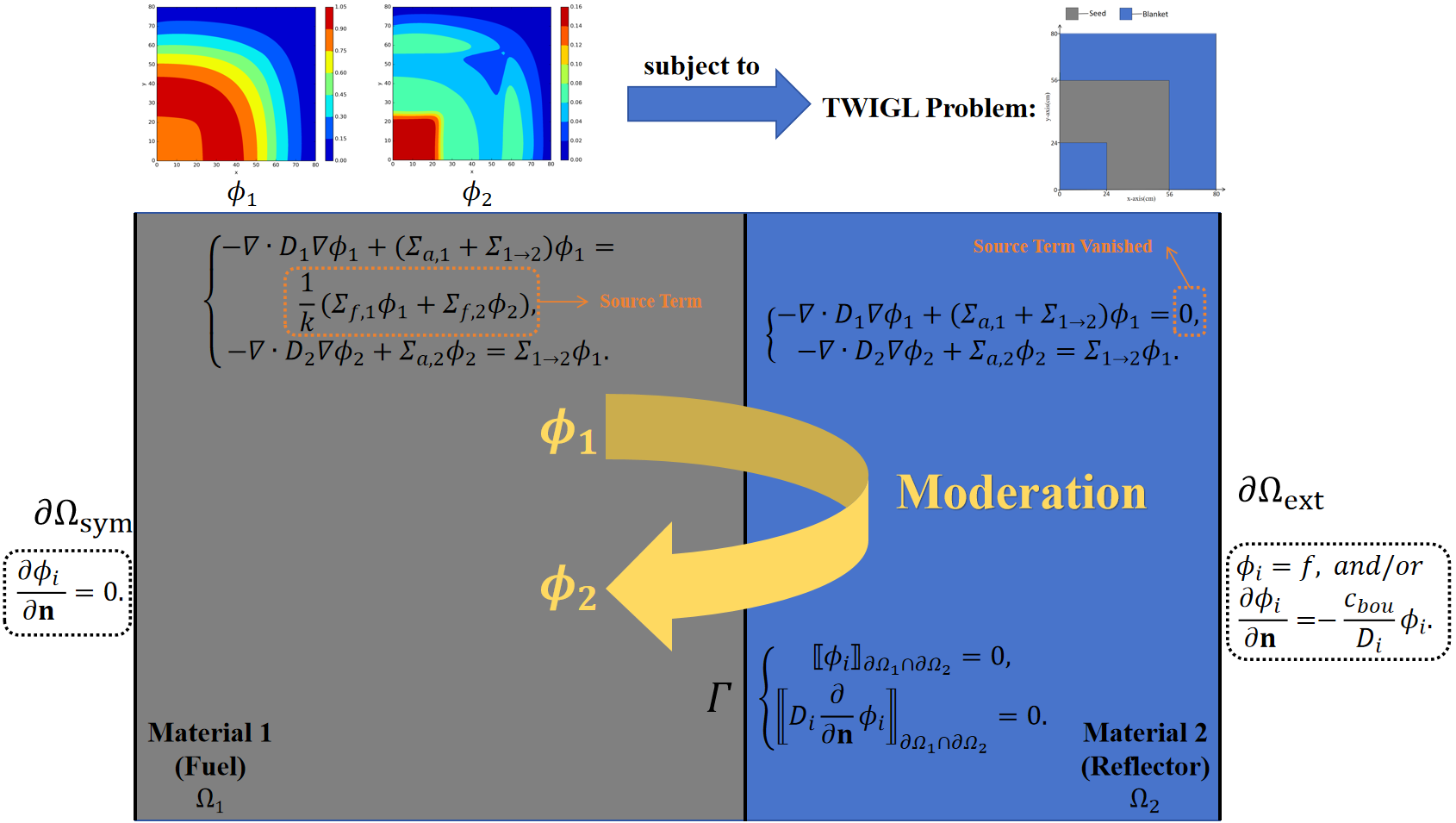}
		\caption{{Schematic diagram of the problem physical meaning.}}
		\label{fig-interface}
	\end{figure*}
	
	{Neutron diffusion eigenvalue problems originate from nuclear reactor cores. The nuclear reactor core is a complex geometrical region composed of multiple materials, which leads to the coefficients in the equations appearing as piecewise constants, and the smoothness of the solutions is also affected at the interfaces. Therefore, we will introduce interface conditions suitable for neutron diffusion eigenvalue problems, which are crucial for solving such problems.}
	{For a problem consisting of two regions (Figure.\ref{fig-interface})}, we impose the following constraints on the solution at the interface:
	\begin{equation}
		\left\{
		\begin{aligned}
			\llbracket \phi_i \rrbracket\big|_{\partial \Omega_1 \cap \partial \Omega_2} & = 0, &~ \text{in} ~ \Omega,\\
			\bigg\llbracket D_i\frac{\partial}{\partial \mathbf{n}}\phi_i \bigg\rrbracket\bigg|_{\partial \Omega_1 \cap \partial \Omega_2} & = 0, &~ \text{in} ~\Omega,\\
		\end{aligned}
		\right.
		\label{eq-int}
	\end{equation}
	for $i = 1, 2$, where 
	\begin{equation}
		\big\llbracket u \big\rrbracket\big|_{\Gamma} := u_1|_{{\Gamma}} - {u_2}|_{{\Gamma}},
	\end{equation}
	\begin{equation}
		\bigg\llbracket a\frac{\partial}{\partial \mathbf{n}}{u} \bigg\rrbracket\bigg|_{\Gamma} := (a_1\frac{\partial}{\partial \mathbf{n}}{u_1})|_{\Gamma} - (a_2\frac{\partial}{\partial \mathbf{n}}{u_2})|_{\Gamma},
	\end{equation}
	\begin{equation}
		u(\boldsymbol{x}) = 
		\left\{ 
		\begin{aligned}
			u_1(\boldsymbol{x}),~~~ & \boldsymbol{x} ~\in~\Omega_1,\\
			u_2(\boldsymbol{x}),~~~ & \boldsymbol{x} ~\in~\Omega_2,
		\end{aligned}\right.
	\end{equation}
	and $\Gamma$ means interface between any pair of materials. If it is a multi-interface problem, we can apply the same conditions \eqref{eq-int} at the interfaces between each pair of materials.
	
	\section{Neural Networks} \label{Neural Networks}
	\begin{figure*}[h]
		\centering
		\includegraphics[width=\linewidth]{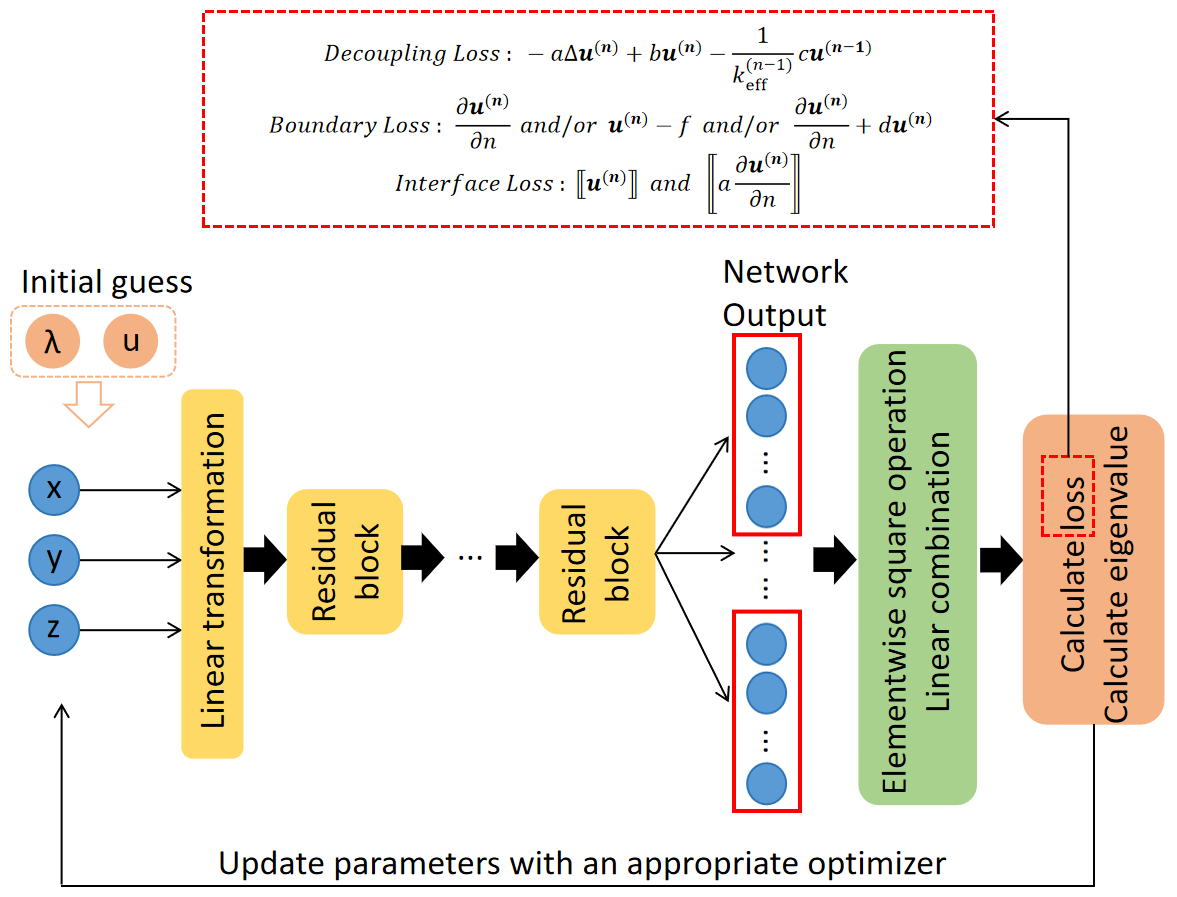}
		\caption{Neural network to solve PDE flow chart.}
		\label{fig-diagram}
	\end{figure*}
	
	In this paper, we choose a residual neural network as an approximator, which typically consists of an input layer, residual block(s), identity connections and output layer. Here, $\boldsymbol{W}_i, \boldsymbol{b}_i$ denote weight matrix and bias, respectively, {and $\rho$ is activation function.} Meanwhile we define $\boldsymbol{h}_{i}$ as residual block. For residual neural networks, the first step is to use linear transformations to change the dimensionality of the input data to match the number of neurons within the residual block. We have the flexibility to choose the number of residual blocks and the number of neurons in every residual block, while the sizes of the input and output need to be determined based on the specific problem at hand. 
	
	If the input data $\boldsymbol{x}_{Input} \in \mathbbm{R}^{d\times 1}$ and the number of neurons in every residual block is $t$, the $\boldsymbol{W}_0 \in \mathbbm{R}^{t\times d}, \boldsymbol{b}_0 \in \mathbbm{R}^{t\times 1}$ and $\boldsymbol{W}_k, \boldsymbol{\tilde {W}}_k \in \mathbbm{R}^{t\times t}, \boldsymbol{b}_k, \boldsymbol{\tilde{b}}_k \in \mathbbm{R}^{t\times 1}, k = 1,\cdots,m-1$. In the residual block, there is also an identity mapping, which we denote as $\boldsymbol{I}$. In the case where we need to solve a system of partial differential equations involving $n$ functions, the output layer should produce $n$ results. Therefore, the dimensions of the weight matrix $\boldsymbol{W}_m \in \mathbbm{R}^{n\times t}$ and bias $\boldsymbol{b}_m \in \mathbbm{R}^{n}$ in the output layer will be determined by the number $n$. However, there is no {activation function $\rho$} before the output layer, as shown in \eqref{eq-outputlayer}. For clarity, we have put a schematic diagram of the algorithm in Figure \ref{fig-diagram}.
	
	\begin{equation}
		\boldsymbol{y}_{0} = \boldsymbol{W}_0 \cdot \boldsymbol{x}_{Input} + \boldsymbol{b}_0,
	\end{equation}
	
	\begin{equation}
		\left\{
		\begin{aligned}
			\boldsymbol{h}_{k-1} &= \boldsymbol{W}_k \cdot \boldsymbol{y}_{k-1} + \boldsymbol{b}_k,&\\
			\boldsymbol{y}_{k} &= \rho \circ (\boldsymbol{\tilde{W}}_k \cdot \boldsymbol{h}_{k-1} + \boldsymbol{\tilde{b}}_k + \boldsymbol{I} \cdot \boldsymbol{y}_{k-1}),& k = 1,\cdots,m-1,\\
		\end{aligned}
		\right.
	\end{equation}
	
	\begin{equation}
		\mathcal{NN}(\boldsymbol{x}) = \boldsymbol{y}_{m} = \boldsymbol{W}_{m}\cdot\boldsymbol{y}_{m-1}+ \boldsymbol{b}_m.
		\label{eq-outputlayer}
	\end{equation}
	
	\subsection{Network Structure}
	For a PDE system with $n$ unknown functions, a straightforward approach is to define the neural network as a $d$-dimensional input and $n$-dimensional output network. However, this approach clearly cannot achieve the crucial interface conditions. For a multi-region problem with $p$ regions, in practice, we can adjust the setup to one neural network by using $n\cdot p$ neurons in the output layer to represent $n$ solutions in each of the $p$ regions. Therefore, we want to build a neural network with $d$-dimension input and $n\cdot p$-dimension output, where $p$ is the number of different regions. $\boldsymbol{u}^{NN}_i = \boldsymbol{u}^{NN}|_{\Omega_i}$ is the restriction of the function $\boldsymbol{u}^{NN}$ to the region $\Omega_i$. Here, $\boldsymbol{u}^{NN}_i$ is a vector of $n$ components. Thus output of neural network is presented as \eqref{eq-NN-output} and interface condition can be implemented as \eqref{eq-NN-int-0} \eqref{eq-NN-int-1}.
	
	\begin{equation}
		\mathcal{NN}(\boldsymbol{x}) = [\boldsymbol{u}^{NN}_1(\boldsymbol{x}); \boldsymbol{u}^{NN}_2(\boldsymbol{x});\cdots;\boldsymbol{u}^{NN}_p(\boldsymbol{x})].
		\label{eq-NN-output} 
	\end{equation}
	
	\begin{equation}
		\big\llbracket \boldsymbol{u}^{NN}\big\rrbracket_{\partial \Omega_i \cap \partial \Omega_j} =\boldsymbol{u}^{NN}_i\big|_{\partial \Omega_i \cap \partial \Omega_j}-\boldsymbol{u}^{NN}_j\big|_{\partial \Omega_i \cap \partial \Omega_j} = 0,~~\forall i\neq j.
		\label{eq-NN-int-0} 
	\end{equation}
	
	\begin{equation}
		\bigg\llbracket \boldsymbol{D}(\boldsymbol{x})\frac{\partial}{\partial \boldsymbol{n}}\boldsymbol{u}^{NN}\bigg\rrbracket_{\partial \Omega_i \cap \partial \Omega_j} = (\boldsymbol{D}_i(\boldsymbol{x})\frac{\partial}{\partial \boldsymbol{n}}\boldsymbol{{u}}^{NN}_i)\big|_{\partial \Omega_i \cap \partial \Omega_j}- (\boldsymbol{D}_j(\boldsymbol{x})\frac{\partial}{\partial \boldsymbol{n}}\boldsymbol{{u}}^{NN}_j)\big|_{\partial \Omega_i \cap \partial \Omega_j} = 0,~~\forall i\neq j.
		\label{eq-NN-int-1} 
	\end{equation}
	The numerical solution by neural network is a weighted linear combination of $\boldsymbol{{u}}^{NN}_i$ respect to $\mathbbm{1}_{\Omega_i}$:
	\begin{equation}
		\boldsymbol{{u}}^{NN} =\sum\limits_{i=1}^p\mathbbm{1}_{\Omega_i}\boldsymbol{{u}}^{NN}_i.
	\end{equation}
	However, in this paper, we only consider the case where $n=2$, as the case where $n=1$ has already been mentioned in \cite{yang2023physics}.
	
	\subsection{Loss function}
	In this section, our loss functions will be presented in accordance with the partitioning of our training points. $\mathbf{X}_{Res}, \mathbf{X}_{Bou}$ and $\mathbf{X}_{Int}$ are used to represent the residual training points, boundary training points, and interface condition training points of the equation, respectively.
	\subsubsection{Residual loss}
	\begin{figure*}[h]
		\centering
		\includegraphics[width=\linewidth]{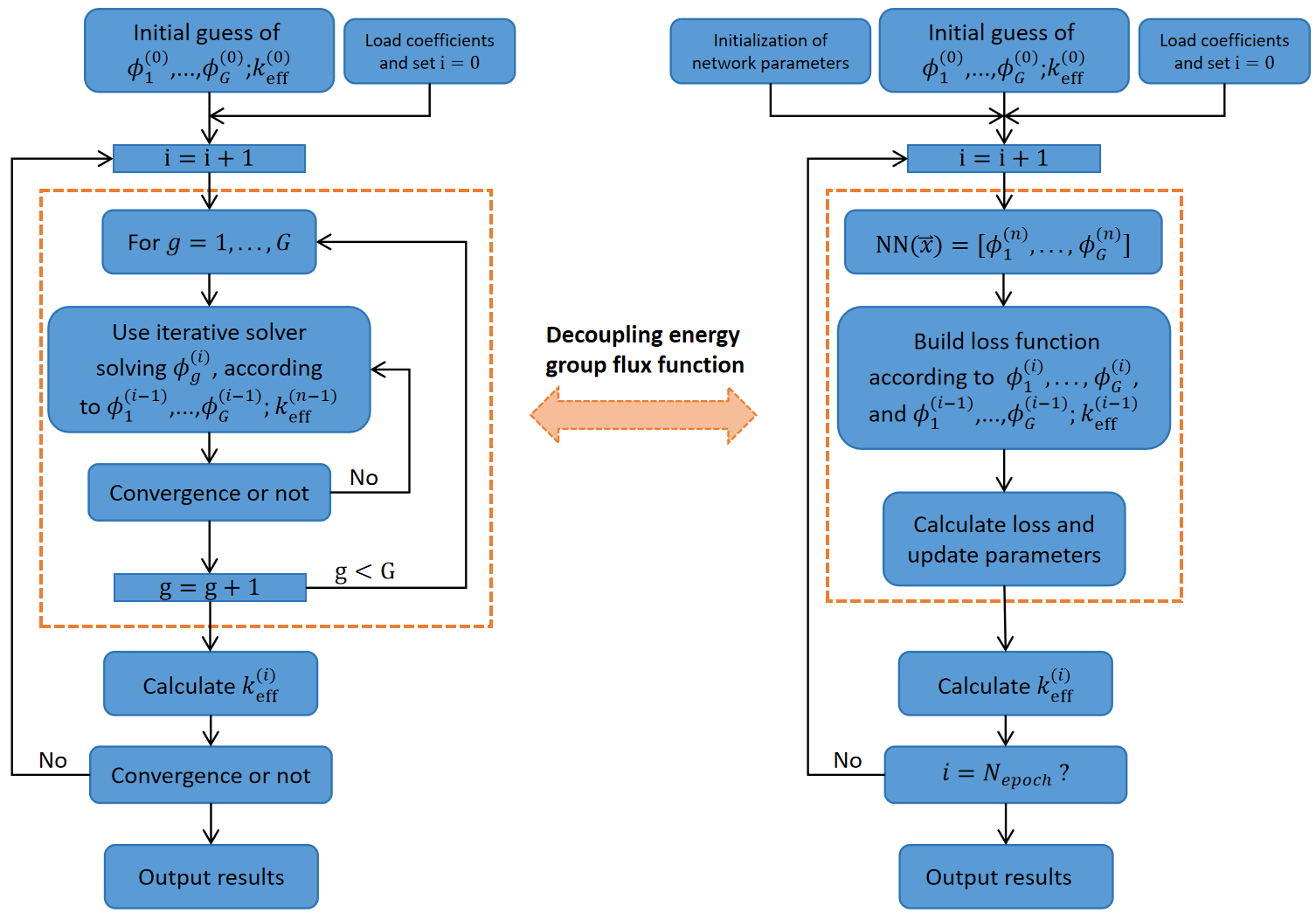}
		\caption{Left: the algorithm flowchart for traditional methods of solving multi-group neutron diffusion problems. Right: the flowchart for solving the same problem using neural networks. The orange dashed boxes represent the decoupling process in both algorithms.}
		\label{fig-tra-NN}
	\end{figure*}
	
	The critical state of the steady-state neutron diffusion equation is characterized by a system of eigenvalue equations involving multiple energy group variables \eqref{eq-multi}. Each equation contains variables from different energy groups and here we consider a general form with a shift of eigenvalue:
	\begin{equation}
		\begin{aligned}
			-\nabla \cdot D_g \nabla \phi_g(\boldsymbol{x}) + \Sigma_{a,g}\phi_g(\boldsymbol{x})-\sum\limits_{g'=1}^{G}\Sigma_{g'\rightarrow g}\phi_{g'}(\boldsymbol{x})+&\sigma\chi_g\sum\limits_{g'=1}^{G}\nu\Sigma_{f,g'}\phi_{g'}(\boldsymbol{x}) \\
			=&(\frac{1}{k_{\text{eff}}}+\sigma)\chi_g\sum\limits_{g'=1}^{G}\nu\Sigma_{f,g'}\phi_{g'}(\boldsymbol{x}), ~ g=1,\cdots,G.
		\end{aligned}
		\label{eq-multi}
	\end{equation}
	
	In conventional computational methods , researchers commonly use source iteration, also known as power iteration \cite{staceyNuclearReactorPhysics}, to solve the aforementioned problem. Firstly, conventional methods need to discretize the equation on the grid points with some specific finite difference scheme. We derive the discretized form into following: 
	\begin{equation}
		(\boldsymbol{D}_g + \boldsymbol{A}_g)\boldsymbol{\phi}_g +\sum\limits_{g'=1}^G\boldsymbol{S}_{g'\rightarrow g}\boldsymbol{\phi}_{g'}+\sigma\chi_g\sum\limits_{g'=1}^G\boldsymbol{F}_{g'}\boldsymbol{\phi}_{g'} = (\frac{1}{k_{\text{eff}}}+\sigma)\chi_g\sum\limits_{g'=1}^G\boldsymbol{F}_{g'}\boldsymbol{\phi}_{g'}, ~ g=1,\cdots,G.
	\end{equation}
	where matrix $\boldsymbol{D}_g, \boldsymbol{A}_g$ and $\boldsymbol{S}_g$ contain diffusion absorption and scattering terms from \eqref{eq-multi} left hand side, matrix $\boldsymbol{F}_g$ represents fission terms from right hand side of \eqref{eq-multi} and vector $\boldsymbol{\phi}_g$ contains the value on the grid points.
	
	The source iteration approach involves providing an initial guess, substituting it into the right-hand side of the equation system, and sequentially solving $G$ equations for different energy groups. This process yields a new set of $\boldsymbol{\phi}_g$ values. Then, based on the physical meaning of $k_{\text{eff}}$, the $k_{\text{eff}}$ value for this iteration step is computed, resulting in a new right-hand side term. During the $i^{th}$ iteration, if we disregard the scattering of neutrons from the lower energy group to the higher energy group \cite{marguet2018physics}, the process can be simplified as follows:
	\begin{equation}
		(\boldsymbol{D}_g + \boldsymbol{A}_g+\sigma\chi_g\boldsymbol{F}_g)\boldsymbol{\phi}_g^{(i)} = \mathcal{RHS}_g^{(i-1)}, ~ g=1,\cdots,G.
		\label{eq-multi-conventional-1}
	\end{equation}
	where 
	\begin{equation}
		\mathcal{RHS}_g^{(i-1)} =  \frac{1}{k^{(i-1)}_{\text{eff}}}\chi_g\sum\limits_{g'=1}^{G}\boldsymbol{F}_{g'}\boldsymbol{\phi}_{g'}^{(i-1)}+\sigma\chi_g\boldsymbol{F}_{g}\boldsymbol{\phi}_{g}^{(i-1)}-\sum\limits_{g'=1}^{g-1}\boldsymbol{S}_{g'\rightarrow g}\boldsymbol{\phi}_{g'}^{(i)}.
		\label{eq-multi-conventional-2}
	\end{equation}
	We rewrite the formation \eqref{eq-multi-conventional-1} \eqref{eq-multi-conventional-2} into block-matrix form for further discussion:
	\begin{equation}
		(\mathcal{D}+\mathcal{A}+\mathcal{S}+\sigma\mathcal{F}_{Diag}){\phi}^{(i)} = \frac{1}{k^{(i-1)}_{\text{eff}}}\mathcal{F}{\phi}^{(i-1)}+\sigma\mathcal{F}_{Diag}{\phi}^{(i-1)},
		\label{multi-eq-mat}
	\end{equation}
	where 
	\begin{equation}
		\begin{aligned}
			\mathcal{D} = 
			\left[ 
			\begin{array}{cccc}
				\boldsymbol{D}_1 & 0 & \cdots & 0\\
				0 & \boldsymbol{D}_2 & \cdots & 0\\
				\vdots & \vdots & \ddots & \vdots\\
				0 & 0 & \cdots & \boldsymbol{D}_G\\
			\end{array} 
			\right ],
			\mathcal{A} = 
			\left[ 
			\begin{array}{cccc}
				\boldsymbol{A}_1 & 0 & \cdots & 0\\
				0 & \boldsymbol{A}_2 & \cdots & 0\\
				\vdots & \vdots & \ddots & \vdots\\
				0 & 0 & \cdots & \boldsymbol{A}_G\\
			\end{array} 
			\right ], \\
			\mathcal{S} = 
			\left[ 
			\begin{array}{cccc}
				0 & 0 & \cdots & 0\\
				\boldsymbol{S}_{1\rightarrow 2} & 0 & \cdots & 0\\
				\vdots & \vdots & \ddots & \vdots\\
				\boldsymbol{S}_{1\rightarrow G} & \boldsymbol{S}_{2\rightarrow G} & \cdots & 0\\
			\end{array} 
			\right ],
			\mathcal{F} = 
			\left[ 
			\begin{array}{cccc}
				\chi_1\boldsymbol{F}_1 & \chi_1\boldsymbol{F}_2 & \cdots & \chi_1\boldsymbol{F}_G\\
				\chi_2\boldsymbol{F}_1 & \chi_2\boldsymbol{F}_2 & \cdots & \chi_2\boldsymbol{F}_G\\
				\vdots & \vdots & \ddots & \vdots\\
				\chi_G\boldsymbol{F}_1 & \chi_G\boldsymbol{F}_2 & \cdots & \chi_G\boldsymbol{F}_G\\
			\end{array} 
			\right ],
		\end{aligned}
	\end{equation}
	and $\phi^{(i)} = [{\boldsymbol{\phi}_1^{(i)}}^T,{\boldsymbol{\phi}_2^{(i)}}^T,\cdots,{\boldsymbol{\phi}_G^{(i)}}^T]^T$, where $\mathcal{D},\mathcal{A},\mathcal{S},\mathcal{F}$ are block matrix and ${\phi}^{(i)}$ are results in $i^{th}$ iteration. The key in the source iteration approach is to modify the left side of \eqref{multi-eq-mat} into diagonal block matrix, and the right hand side of \eqref{multi-eq-mat} is replaced by a known source term:
	
	\begin{equation}
		(\mathcal{D}+\mathcal{A}+\sigma\mathcal{F}_{Diag}){\phi}^{(i)} = (\sigma\mathcal{F}_{Diag}+\frac{1}{k^{(i-1)}_{\text{eff}}}\mathcal{F})\phi^{(i-1)}-\mathcal{S}\phi^{(i)}.
		\label{multi-eq-mat-mod}
	\end{equation}
	
	As a result, we can equivalently transform the solution of equation \eqref{eq-multi} into a sequential solution of equation \eqref{eq-multi-conventional-1} or \eqref{multi-eq-mat-mod} for $g= 1,\cdots,G$. In this case, the right-hand side becomes a fixed source term rather than an eigenvalue term. 
	In this paper, we consider a neutron diffusion equation with $G=2$, as described in \eqref{eq-multi}:
	\begin{equation}
		\left\{
		\begin{aligned}
			-\nabla \cdot D_1 \nabla \phi_1(\boldsymbol{x}) + (\Sigma_{a,1}+\Sigma_{1\rightarrow 2})\phi_1(\boldsymbol{x})+&\sigma(\nu\Sigma_{f,1}\phi_1(\boldsymbol{x})+\nu\Sigma_{f,2}\phi_2(\boldsymbol{x})) &\\ 
			=&(\frac{1}{k_{\text{eff}}}+\sigma)(\nu\Sigma_{f,1}\phi_1(\boldsymbol{x})+\nu\Sigma_{f,2}\phi_2(\boldsymbol{x})), & ~ \text{in} ~ \Omega,\\
			\centering
			-\nabla \cdot D_2 \nabla \phi_2(\boldsymbol{x}) + \Sigma_{a,2}\phi_2(\boldsymbol{x}) - \Sigma_{1\rightarrow 2}\phi_1(\boldsymbol{x}) &= 0, & ~ \text{in} ~\Omega.
		\end{aligned}
		\right.
		\label{eq-diff-general}
	\end{equation}
	
	Based on the decoupling concept of conventional solving methods for multi-group neutron diffusion equations \eqref{eq-multi-conventional-1}-\eqref{multi-eq-mat-mod}, we have also adopted a similar technique in designing the loss function for \eqref{eq-diff-general}, where the fast and thermal group fluxes are independently solved using two distinct iterative equations. Therefore, we replace the remaining irrelevant group fluxes with the functions obtained in the previous step of the neural network \eqref{De_lhs2} \eqref{De_rhs2} while the relevant group functions are represented by the current output of the neural network \eqref{De_lhs1} \eqref{De_rhs1}. This ensures that the entire residual function corresponds to a partial differential equation with a single unknown function. By utilizing \eqref{loss-res-3}, we obtain the final form of the loss function, which is named the Decoupling loss function. It is worth noting that in the loss function of neural networks, there is no need for discretization of the equations. We typically employ automatic differentiation \cite{paszke2017automatic} to handle the derivative operations that arise in the equations.
	
	\begin{subequations}
		\begin{align}
			\mathcal{LHS}_{De,1} &= -\nabla \cdot D_1 \nabla \phi_1^{NN(i)} + (\Sigma_{a,1}+\Sigma_{1\rightarrow 2})\phi_1^{NN(i)} + \sigma\nu\Sigma_{f,1}\phi_1^{NN(i)}\label{De_lhs1}\\
			\mathcal{LHS}_{De,2} &= -\nabla \cdot D_2 \nabla \phi_2^{NN(i)} + \Sigma_{a,2}\phi_2^{NN(i)}\label{De_lhs2}\\
			\mathcal{RHS}_{De,1} &= \sigma\nu\Sigma_{f,1}\phi_1^{NN(i-1)}+ \frac{1}{k^{(i-1)}_{\text{eff}}}(\nu\Sigma_{f,1}\phi_1^{NN(i-1)}+\nu\Sigma_{f,2}\phi_2^{NN(i-1)})\label{De_rhs1}\\
			\mathcal{RHS}_{De,2} &= \Sigma_{1\rightarrow 2}\phi_1^{NN(i-1)}\label{De_rhs2}
		\end{align}
	\end{subequations}
	
	\begin{equation}
		\begin{aligned}
			Loss_{De} = & \sum\limits_{\boldsymbol{x} \in \mathbf{X}_{res}} \big(\mathcal{LHS}_{De,1}-\mathcal{RHS}_{De,1}\big)^2_{\boldsymbol{x}}  + \sum\limits_{\boldsymbol{x} \in \mathbf{X}_{res}}\big(\mathcal{LHS}_{De,2}-\mathcal{RHS}_{De,2}\big)^2_{\boldsymbol{x}}, 
		\end{aligned}
		\label{loss-res-3}
	\end{equation}
	
	We present the application of the neural network solving process using the new Decoupling loss function and the traditional solving method in Figure \ref{fig-tra-NN} to illustrate their connection and differences.
	
	Alternatively, in the case where decoupling is not considered, if the inverse power method is used to solve for the minimum eigenvalue, the terms involving eigenvalue pairs usually rely on information from the previous iteration step. However, in certain problem settings, the information of $\phi_1$ can be lost due to coefficients $\nu\Sigma_{f,1}$ being zero. In \eqref{DI_lhs1}-\eqref{DI_rhs2}, we separate the equation into parts that are related to the eigenvalue and those that are unrelated. We place them on the left and right sides of the equation respectively. Then, we replace $\phi_1^{NN(i)}$ in the second equation with $\phi_1^{NN(i-1)}$ \eqref{DI_rhs1}, in order to maintain the iteration information of the previous step's $\phi_1$ when $\nu\Sigma_{f,1}=0$. Therefore, we also propose a Direct Iterative loss function \eqref{loss-res-2} that addresses this issue:
	
	\begin{subequations}
		\begin{align}
			\mathcal{LHS}_{DI,1} &= -\nabla \cdot D_1 \nabla \phi_1^{NN(i)} + (\Sigma_{a,1}+\Sigma_{1\rightarrow 2})\phi_1^{NN(i)} + \sigma(\nu\Sigma_{f,1}\phi_1^{NN(i)}+\nu\Sigma_{f,2}\phi_2^{NN(i)})\label{DI_lhs1}\\
			\mathcal{LHS}_{DI,2} &= -\nabla \cdot D_2 \nabla \phi_2^{NN(i)} + \Sigma_{a,2}\phi_2^{NN(i)}-\Sigma_{1\rightarrow 2}\phi_1^{NN(i-1)}\label{DI_lhs2}\\
			\mathcal{RHS}_{DI,1} &= (\frac{1}{k^{(i-1)}_{\text{eff}}}+\sigma)(\nu\Sigma_{f,1}\phi_1^{NN(i-1)}+\nu\Sigma_{f,2}\phi_2^{NN(i-1)})\label{DI_rhs1}\\
			\mathcal{RHS}_{DI,2} &= 0 \label{DI_rhs2}
		\end{align}
		\label{eq-DI-LRHS}
	\end{subequations}
	
	\begin{equation}
		\begin{aligned}
			Loss_{DI} = & \sum\limits_{\boldsymbol{x} \in \mathbf{X}_{res}} \big(\mathcal{LHS}_{DI,1}-\mathcal{RHS}_{DI,1}\big)^2_{\boldsymbol{x}}  + \sum\limits_{\boldsymbol{x} \in \mathbf{X}_{res}}\big(\mathcal{LHS}_{DI,2}-\mathcal{RHS}_{DI,2}\big)^2_{\boldsymbol{x}}.
		\end{aligned}
		\label{loss-res-2}
	\end{equation}
	
	In the subsequent numerical experiments, we will use these two types of residual functions (decoupling loss function\eqref{loss-res-3} and direct iterative loss function \eqref{loss-res-2}) as residual loss function, which are denoted as $Loss_{Res}$. We will then proceed to compare the performance of these two loss functions.
	
	\subsubsection{Boundary loss}
	
	However, for the boundary training points $\mathbf{X}_{Bou}$, the specific form of the loss function needs to be determined based on the applied boundary conditions. Therefore, we provide a more detailed classification here, dividing $\mathbf{X}_{Bou}$ into $\mathbf{X}_{D}$, $\mathbf{X}_{N}$, and $\mathbf{X}_{R}$, and defining separate boundary loss functions \eqref{loss-diri} \eqref{loss-neumann} \eqref{loss-robin} for each. Clearly, $\mathbf{X}_{Bou} = \mathbf{X}_D \cup \mathbf{X}_N \cup \mathbf{X}_R$. For the sake of convenience in the subsequent discussion, we collectively refer to the loss functions \eqref{loss-diri} \eqref{loss-neumann} \eqref{loss-robin} corresponding to $\mathbf{X}_{D}$, $\mathbf{X}_{N}$, and $\mathbf{X}_{R}$ as $Loss_{Bou}$. 
	
	\begin{subequations}
		\begin{align}
			Loss_{D} =& \sum\limits_{\boldsymbol{x} \in \mathbf{X}_{D}} \big(\phi_1^{NN}-f\big)^2_{\boldsymbol{x}}+\big(\phi_2^{NN}-f\big)^2_{\boldsymbol{x}}.\label{loss-diri}\\
			Loss_{N} =& \sum\limits_{\boldsymbol{x} \in \mathbf{X}_{N}} \big(\frac{\partial \phi_1^{NN}}{\partial \boldsymbol{n}}-0\big)^2_{\boldsymbol{x}}+\big(\frac{\partial \phi_2^{NN}}{\partial \boldsymbol{n}}-0\big)^2_{\boldsymbol{x}}.\label{loss-neumann}\\
			Loss_{R} =& \sum\limits_{\boldsymbol{x} \in \mathbf{X}_{R}} \big( \frac{\partial \phi_1^{NN}}{\partial \boldsymbol{n}} +\frac{c_{bou}}{D_1}\phi_1^{NN}\big)^2_{\boldsymbol{x}}+\big( \frac{\partial \phi_2^{NN}}{\partial \boldsymbol{n}} +\frac{c_{bou}}{D_2}\phi_2^{NN}\big)^2_{\boldsymbol{x}}.\label{loss-robin} 
		\end{align}
	\end{subequations}
	
	\subsubsection{Interface loss}
	
	When it comes to setting the loss function, solely relying on boundary conditions and the governing equations is often insufficient to solve the multi-region eigenvalue problem. We must introduce crucial interface conditions \eqref{eq-int} between different materials to accurately capture the behavior of the system.
	For the interface condition training points $\mathbf{X}_{Int}$, since there are constraints on both the primitive function and its first-order derivative, we have defined two types of loss functions, \eqref{loss-primitive} and \eqref{loss-derivative}, respectively.
	
	\begin{equation}
		Loss_{Int0} = \sum\limits_{\boldsymbol{x}\in\mathbf{X}_{Int}}\sum\limits_{\forall p\neq q} \big(\big\llbracket \phi_1^{NN} \big\rrbracket _{\partial \Omega_p \cap \partial \Omega_q}\big)^2_{\boldsymbol{x}}+\big(\big\llbracket \phi_2^{NN} \big\rrbracket _{\partial \Omega_p \cap \partial \Omega_q}\big)^2_{\boldsymbol{x}}.\label{loss-primitive}\\
	\end{equation}
	\begin{equation}
		Loss_{Int1} = \sum\limits_{\boldsymbol{x}\in\mathbf{X}_{Int}}\sum\limits_{\forall p\neq q} \big(\bigg\llbracket D_1\frac{\partial \phi_1^{NN}}{\partial \boldsymbol{n}} \bigg\rrbracket _{\partial \Omega_p \cap \partial \Omega_q}\big)^2_{\boldsymbol{x}} + \big(\bigg\llbracket D_2\frac{\partial \phi_2^{NN}}{\partial \boldsymbol{n}} \bigg\rrbracket _{\partial \Omega_p \cap \partial \Omega_q}\big)^2_{\boldsymbol{x}}.\label{loss-derivative}
	\end{equation}
	In the subsequent numerical experiments section, different scenarios will correspond to different boundary losses. The specific settings will be elaborated in Section \ref{Numerical Experiments} of the paper. The total loss should be written in the following form, where the $\alpha$ represent the weights for different loss terms:
	\begin{equation}
		\begin{aligned}
			Loss =& \alpha_{Res}Loss_{Res}+\alpha_{Bou}Loss_{Bou}+\alpha_{Int0}Loss_{Int0}+\alpha_{Int1}Loss_{Int1}.
		\end{aligned}
		\label{loss-total-1}
	\end{equation}
	
	\subsection{Eigenvalue}
	For the critical problem in nuclear reactors, it is not only necessary to determine the unknown functions but also to find the eigenvalues. In non-eigenvalue problems, each coefficient in the equation is known, which significantly simplifies the problem-solving process. However, in eigenvalue problems, the challenge lies in efficiently obtaining eigenvalues to an acceptable level of accuracy and subsequently utilizing neural networks to compute the corresponding eigenfunctions. This is where the introduction of the Rayleigh-Quotient \eqref{eq-RayleighQ} becomes crucial. 
	{Here, only the set of points used for the residual part $\mathbf{X}_{Res}$ has been selected to approximate the expression for the Rayleigh-Quotient.}
	\begin{equation}
		\begin{aligned}
			\tilde{\lambda}:=\frac{1}{k_{\text{eff}}}+\sigma =& \frac{\int_\Omega (-\nabla \cdot D_1 \nabla \phi_1 + (\Sigma_{a,1}+\Sigma_{1\rightarrow 2})\phi_1 + \sigma\chi_1(\nu\Sigma_{f,1}\phi_1+\nu\Sigma_{f,2}\phi_2))\phi_1 d\boldsymbol{x}}{\int_\Omega(\chi_1\phi_1+\chi_2\phi_2)(\nu\Sigma_{f,1}\phi_1+\nu\Sigma_{f,2}\phi_2)d\boldsymbol{x}}\\
			+&\frac{\int_\Omega(-\nabla \cdot D_2 \nabla \phi_2 + \Sigma_{a,2}\phi_2-\Sigma_{1\rightarrow 2}\phi_1+\sigma\chi_2(\nu\Sigma_{f,1}\phi_1+\nu\Sigma_{f,2}\phi_2))\phi_2 d\boldsymbol{x}}{\int_\Omega(\chi_1\phi_1+\chi_2\phi_2)(\nu\Sigma_{f,1}\phi_1+\nu\Sigma_{f,2}\phi_2)d\boldsymbol{x}}\\
			\approx& \frac{\sum\limits_{\boldsymbol{x}\in\mathbf{X}_{Res}} \big[(-\nabla \cdot D_1 \nabla \phi_1 + (\Sigma_{a,1}+\Sigma_{1\rightarrow 2})\phi_1 + \sigma\chi_1(\nu\Sigma_{f,1}\phi_1+\nu\Sigma_{f,2}\phi_2))\phi_1\big]_{\boldsymbol{x}}} {\sum\limits_{i=1}^M \big((\chi_1\phi_1+\chi_2\phi_2)(\nu\Sigma_{f,1}\phi_1+\nu\Sigma_{f,2}\phi_2)\big)_{\boldsymbol{x}}}\\
			+&\frac{\sum\limits_{\boldsymbol{x}\in\mathbf{X}_{Res}} \big[(-\nabla \cdot D_2 \nabla \phi_2 + \Sigma_{a,2}\phi_2-\Sigma_{1\rightarrow 2}\phi_1+\sigma\chi_2(\nu\Sigma_{f,1}\phi_1+\nu\Sigma_{f,2}\phi_2))\phi_2\big]_{\boldsymbol{x}}} {\sum\limits_{i=1}^M \big((\chi_1\phi_1+\chi_2\phi_2)(\nu\Sigma_{f,1}\phi_1+\nu\Sigma_{f,2}\phi_2)\big)_{\boldsymbol{x}}}.
		\end{aligned}
		\label{eq-RayleighQ}
	\end{equation}
	
	\subsection{Algorithm}
	In the case of critical in nuclear reactor physics, we only need to find the minimum eigenvalue. Moreover, at this stage, the corresponding physical quantity on the computational domain should be non-negative. Therefore, our neural network design includes a post-processing step to ensure that the output satisfies the non-negativity condition. We achieve this by applying an element-wise squaring operation to the output of the neural network $\mathcal{NN}^2$, ensuring that the computed function values are non-negative. Our algorithm is presented in Algorithm \ref{Alg-1}.
	
	\begin{algorithm}
		\caption{Training process}
		\label{Alg-1}
		\KwIn{The initial guess: $\tilde{\lambda}_0:=\frac{1}{k^{(0)}_{\text{eff}}}, \phi_{1}^{NN(0)},\phi_{2}^{NN(0)}.$ }
		\KwIn{The initial neural network parameters: $\Theta.$ }
		\KwIn{Training points: $\mathbf{X}_{Res},~\mathbf{X}_{Bou},~\mathbf{X}_{Int}$.}
		\KwIn{Training epochs: $N_{epoch}$.}
		\For {$i=1:N_{epoch}$}
		{
			$u=\sum\limits_{j=1}^p\mathbbm{1}_{\Omega_j}u^{NN}_j, v=\sum\limits_{j=1}^p\mathbbm{1}_{\Omega_j}v^{NN}_j \gets \mathcal{NN}^2(\mathbf{X}_{Res};\Theta)$\;
			
			Calculate power $\boldsymbol{p} \gets \int_\Omega (u~\nu\Sigma_{f,1}+v~\nu\Sigma_{f,2})d\mathbf{x}$\;
			Normalize by power $\phi_{1}^{NN(i)},\phi_{2}^{NN(i)} \gets \frac{u}{\boldsymbol{p}}, \frac{v}{\boldsymbol{p}}$\;
			
			Loss of residual$\gets Loss_{Res}(\phi_{1}^{NN(i)},\phi_{2}^{NN(i)},\phi_{1}^{NN(i-1)},\phi_{2}^{NN(i-1)},\tilde{\lambda}_{i-1})$\;
			
			$\phi_{1,Bou}^{NN}=\sum\limits_{j=1}^p\mathbbm{1}_{\Omega_j}u^{NN}_j, \phi_{2,Bou}^{NN}=\sum\limits_{j=1}^p\mathbbm{1}_{\Omega_j}v^{NN}_j\gets \mathcal{NN}^2(X_{Bou};\Theta)$\;
			Loss of boundary$\gets Loss_{Bou}(\phi_{1,Bou}^{NN}, \phi_{2,Bou}^{NN})$\;
			
			$\phi_{1,Int}^{NN}=\sum\limits_{j=1}^p\mathbbm{1}_{\Omega_j}u^{NN}_j, \phi_{2,Int}^{NN}=\sum\limits_{j=1}^p\mathbbm{1}_{\Omega_j}v^{NN}_j\gets \mathcal{NN}^2(X_{Int};\Theta)$\;
			Loss of interface$\gets Loss_{Int0}(\phi_{1,Int}^{NN}, \phi_{2,Int}^{NN}), Loss_{Int1}(\phi_{1,Int}^{NN}, \phi_{2,Int}^{NN})$\;
			
			Update neural networks parameters $\Theta$ with respect to Eq.(\ref{loss-total-1})\;
			
			$\tilde{\lambda}_{i} \gets$ in Rayleigh-Quotient form Eq.(\ref{eq-RayleighQ})\;
			\If{\text{achieve the minimum loss}}{Save the model and corresponding
				$\tilde{\lambda}_s = \tilde{\lambda}_{i}, \phi_1 = \phi_1^{NN(s)}, \phi_2 = \phi_2^{NN(s)}$.}
			$i = i + 1$\;
		}
		\KwResult{The eigenvalue $\frac{1}{k_{\text{eff}}}=\tilde{\lambda}_s-\sigma$ and the eigenfunction $\vec{\phi}=(\phi_1,\phi_2)^T$.}
	\end{algorithm}
	
	\section{Numerical Experiments} \label{Numerical Experiments}
	In this section, we plan to show numerical results for solving few eigenvalue problems. We give statements of problem to be solved by our method. And we will propose corresponding evaluation criteria to verify the validity of numerical results in Section \ref{Error criterion}. Further, we will give details, numerical solutions and related evaluation results of one-dimensional problem \cite{demaziere2007calculation}, two-dimensional and three-dimensional TWIGL problem \cite{hageman1969comparison}, and two-dimensional and three-dimensional IAEA problem \cite{theler2011solution} with different residual loss function applied. 
	
	As mentioned earlier, we have a total of three training point sets $\mathbf{X}_{Res}, \mathbf{X}_{Bou}, \mathbf{X}_{Int}$ for our training points. The boundary point set $\mathbf{X}_{Bou}$ consists of points strictly defined on the boundary, {and different examples will be used as uniform sampling points on boundary according to different resolutions.} The interface point $\mathbf{X}_{Int}$ set and residual point $\mathbf{X}_{Res}$ set are disjoint, meaning that if a coordinate point is subject to an interface condition, it will not be subjected to the equation loss condition, even if it satisfies both conditions.
	
	\subsection{Error criterion}\label{Error criterion}
	
	In this section, we propose our evaluation method for the quality of numerical solutions. Additionally, due to the engineering application of the problem, we will provide acceptance criteria for the numerical solutions in practical problems. We will continue to use the symbols introduced in Section 3 and define our solution as $\vec{\boldsymbol{\phi}}^{NN} = (\phi_1^{NN}, \phi_2^{NN})^T, ~ k_{\text{eff}}^{NN} := \frac{1}{\lambda^{NN}}$.To evaluate the quality of the numerical solutions, we will introduce a high-fidelity solution obtained by FreeFEM++ \cite{hecht2012new} as our reference solution, which is defined as $\vec{\boldsymbol{\phi}}^{FF} = (\phi_1^{FF}, \phi_2^{FF})^T, ~ k_{\text{eff}}^{FF} := \frac{1}{\lambda^{FF}}$.
	
	Next, we will provide different relative error formulation for the solution in terms of the $L^p$-norm:
	
	\begin{equation}
		\mathbf{E}_{R,p}(\phi_i) = \frac{||\phi_i^{NN}-\phi_i^{FF}||_{L^p(\Omega)}}{||\phi_i^{FF}||_{L^p(\Omega)}},~~i = 1, 2,
	\end{equation}
	where the subscript $R$ denotes relative, and the superscript $p$ means $L^p$-norm. Here, $p$ can only take $2$ and $\infty$.
	
	Since the eigenvalues in these cases are scalars and real numbers, the relative error with respect to the minimal eigenvalue is defined as follows:
	\begin{equation}
		\mathbf{E}_R(k_{\text{eff}}) = \frac{|k_{\text{eff}}^{NN}-k_{\text{eff}}^{FF}|}{|k_{\text{eff}}^{FF}|}.
	\end{equation}
	
	\subsection{1-D Problem}
\begin{figure*}[h]
		\centering
		\includegraphics[width=15cm]{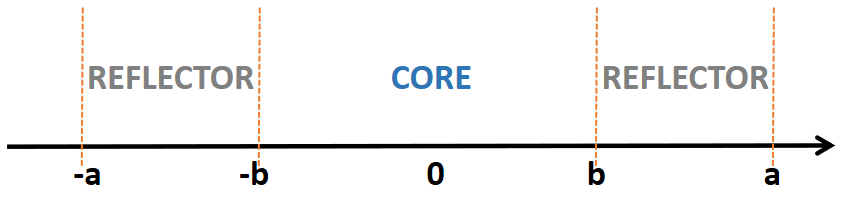}
		\caption{Computational domain of 1-D Swedish Ringhals-4 pressurized water reactor, where $a = 279.5\text{cm}, b = 161.25\text{cm}$.}
		\label{fig-1d-domain}
	\end{figure*}

	This problem originates from the core of the Swedish Ringhals-4 pressurized water reactor. The cross-section data for various materials were obtained from simulation calculations, then transformed using homogenization techniques to make them applicable to the one-dimensional scenario. In this paper, this is considered as the benchmark problem for the one-dimensional case, and the corresponding data is presented in Table \ref{tab-coef-1D}. The computational domain for this problem is depicted in the Figure \ref{fig-1d-domain}, consisting of three regions and two types of materials. Robin boundary conditions are applied at the endpoints of the domain, so we apply \eqref{loss-robin} as the boundary loss function and $c_{bou}=0.5$. 
	
	We solved the problem using FreeFEM++ on a uniform grid with $\Delta x = 0.05$ as a reference solution. The reference eigenvalue is $k^{FF}_{\text{eff}} = 1.0037$. {Similarly, for our neural network, we selected 11,181 uniformly distributed sampling points as our training data, where 11,179 are residual points and 2 are interface points.} We chose 2 residual blocks, with 20 neurons per layer in each residual block. The Adam optimizer with a learning rate of $0.001$ was used for $100,000$ epochs optimization during the training process.
	
	\begin{table}[h]
		\centering
		\caption{Coefficients of different regions for 1-D Swedish Ringhals-4 pressurized water reactor problem}
		\begin{tabular}{|c|c|c|c|c|c|c|c|c|}
			\hline
			& $D_1$ & $D_2$ & $\Sigma_{a,1}$ & $\Sigma_{a,2}$ & $\Sigma_{1\rightarrow2}$ & $\nu\Sigma_{f,1}$ & $\nu\Sigma_{f,2}$\\
			& $(cm)$ & $(cm)$ & $(cm^{-1})$ & $(cm^{-1})$ & $(cm^{-1})$ & $(cm^{-1})$ & $(cm^{-1})$\\
			\hline
			Core & 1.4376 & 0.3723 & 0.0115 & 0.1019 & 0.0151 & 0.0057 & 0.1425\\
			\hline
			Reflector & 1.3116 & 0.2624 & -0.0098 & 0.0284 & 0.0238 & 0.0 & 0.0\\
			\hline
		\end{tabular}
		\label{tab-coef-1D}
	\end{table}
	
\begin{figure*}[h]
		\begin{minipage}{0.32\textwidth}
			\centering
			\includegraphics[width=5.25cm]{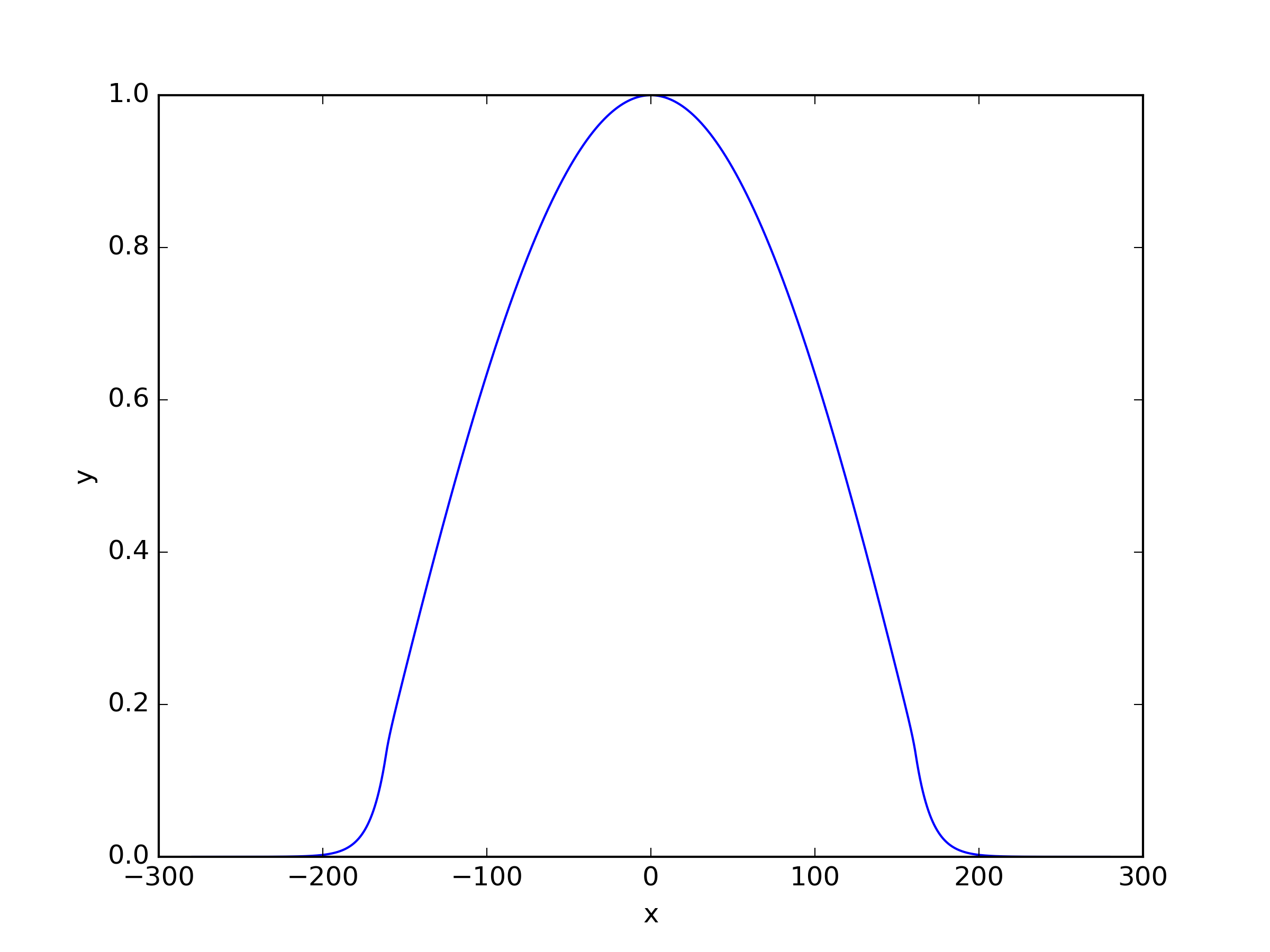}
		\end{minipage}
		\begin{minipage}{0.32\textwidth}
			\centering
			\includegraphics[width=5.25cm]{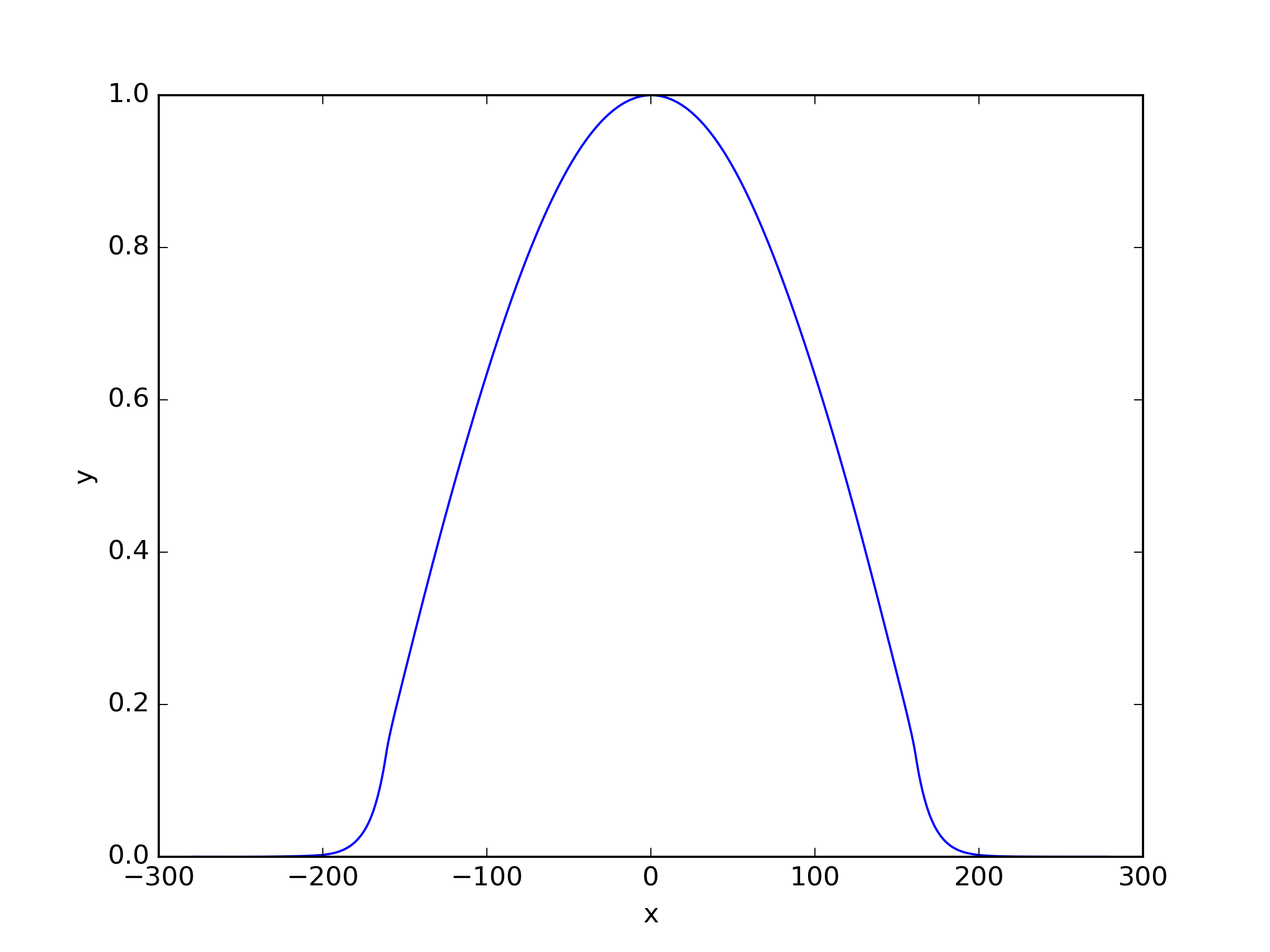}
		\end{minipage}
		\begin{minipage}{0.32\textwidth}
			\centering
			\includegraphics[width=5.25cm]{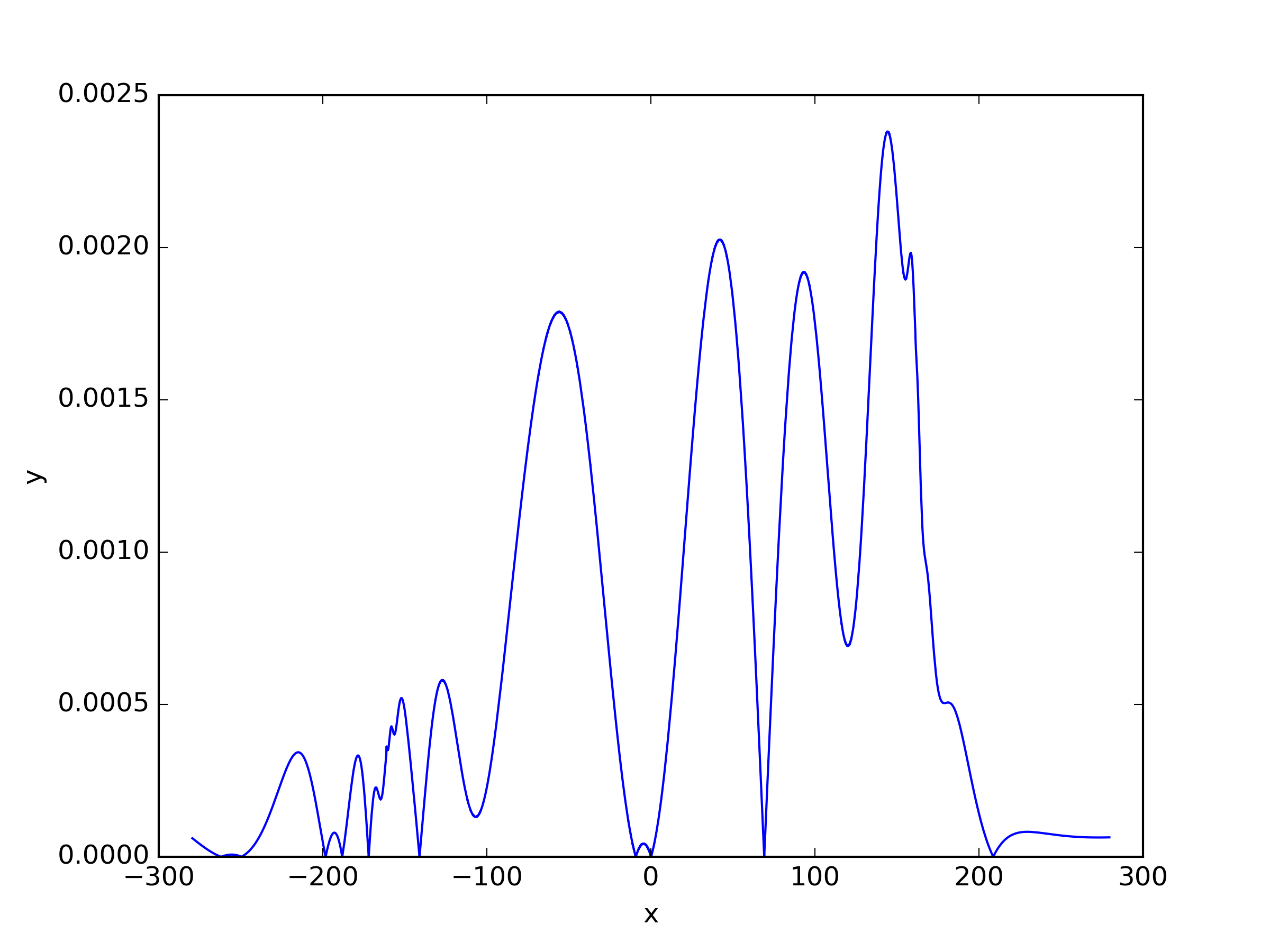}
		\end{minipage}
		
		\begin{minipage}{0.32\textwidth}
			\centering
			\includegraphics[width=5.25cm]{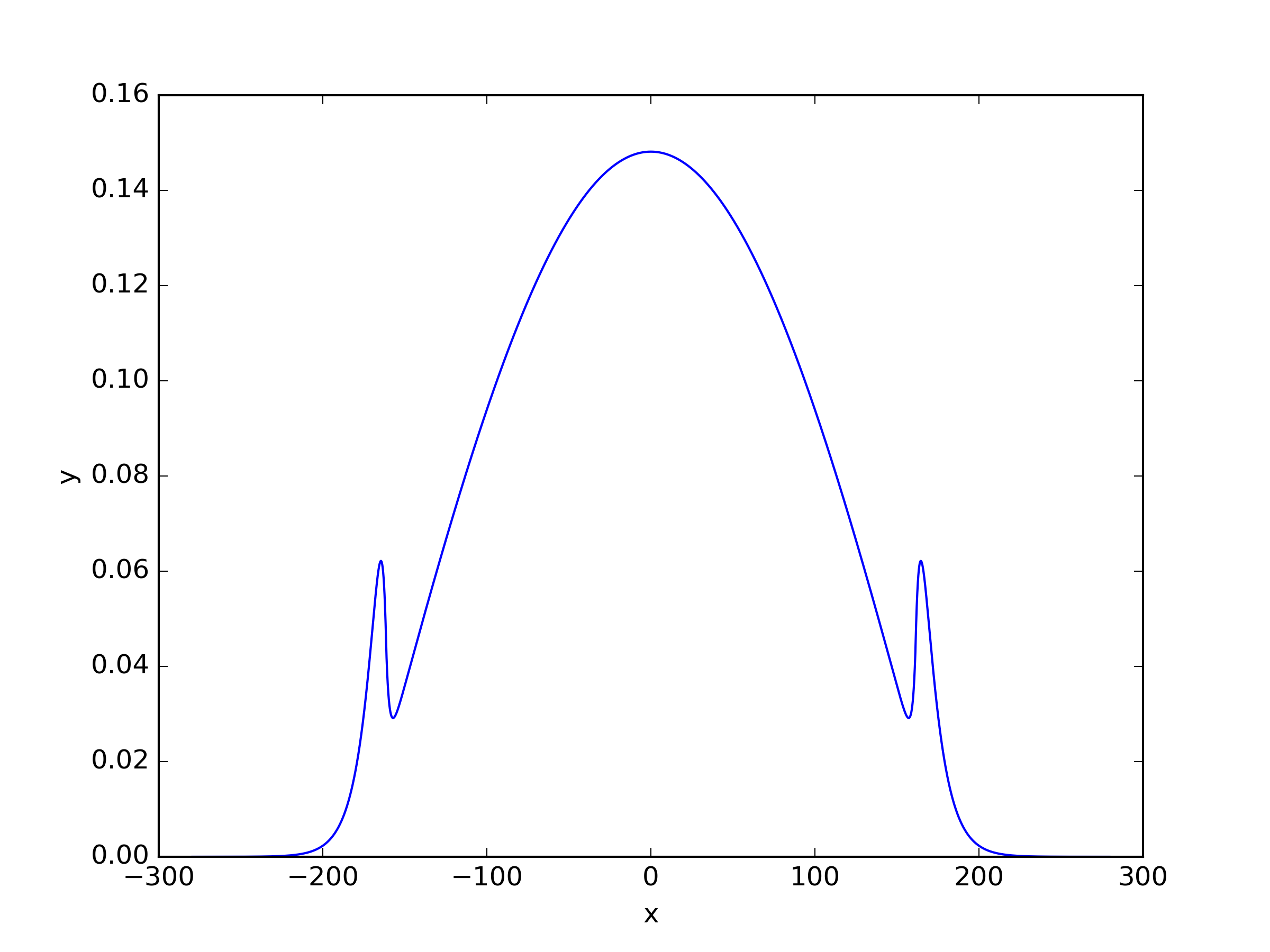}
		\end{minipage}
		\begin{minipage}{0.32\textwidth}
			\centering
			\includegraphics[width=5.25cm]{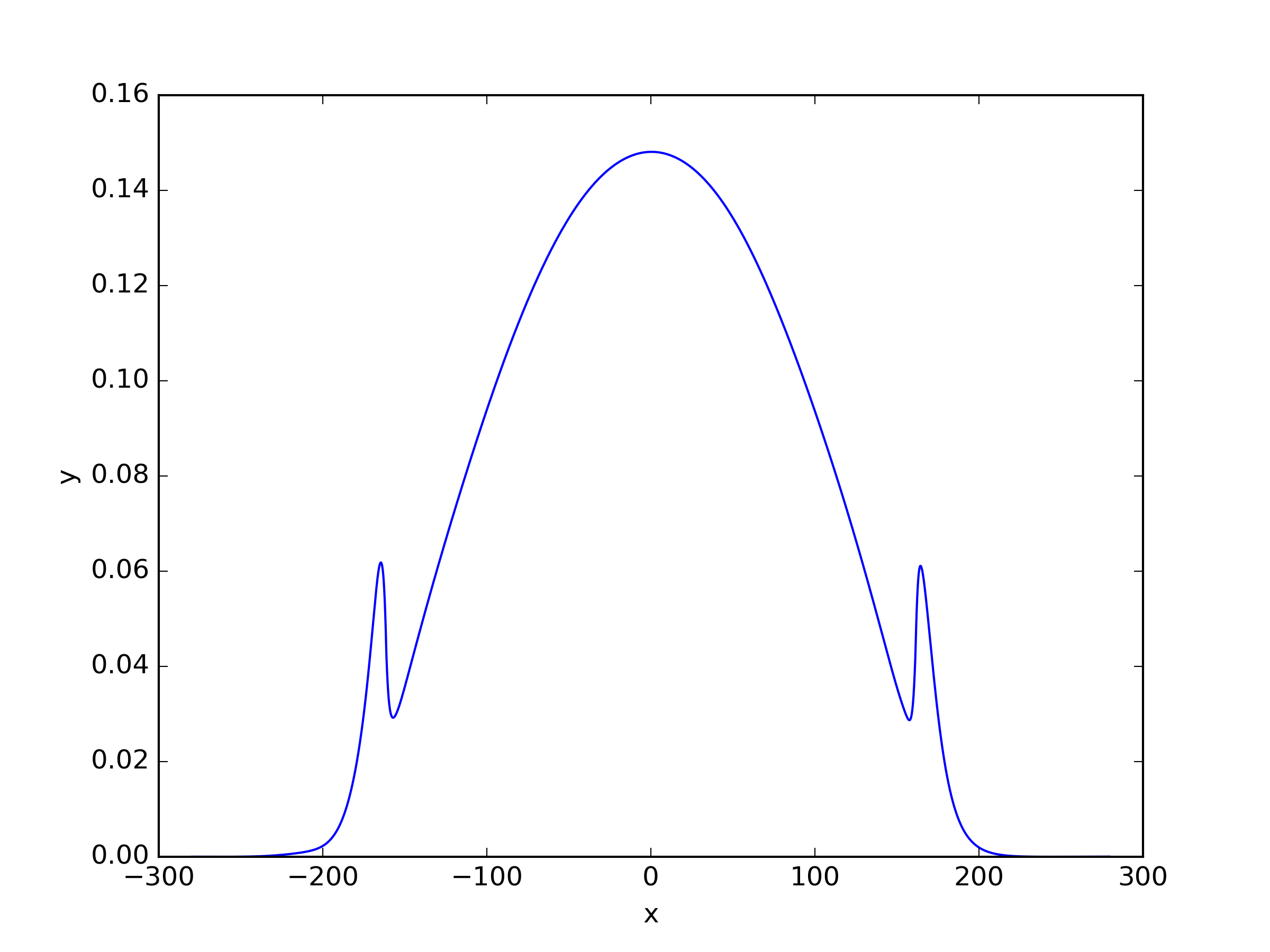}
		\end{minipage}
		\begin{minipage}{0.32\textwidth}
			\centering
			\includegraphics[width=5.25cm]{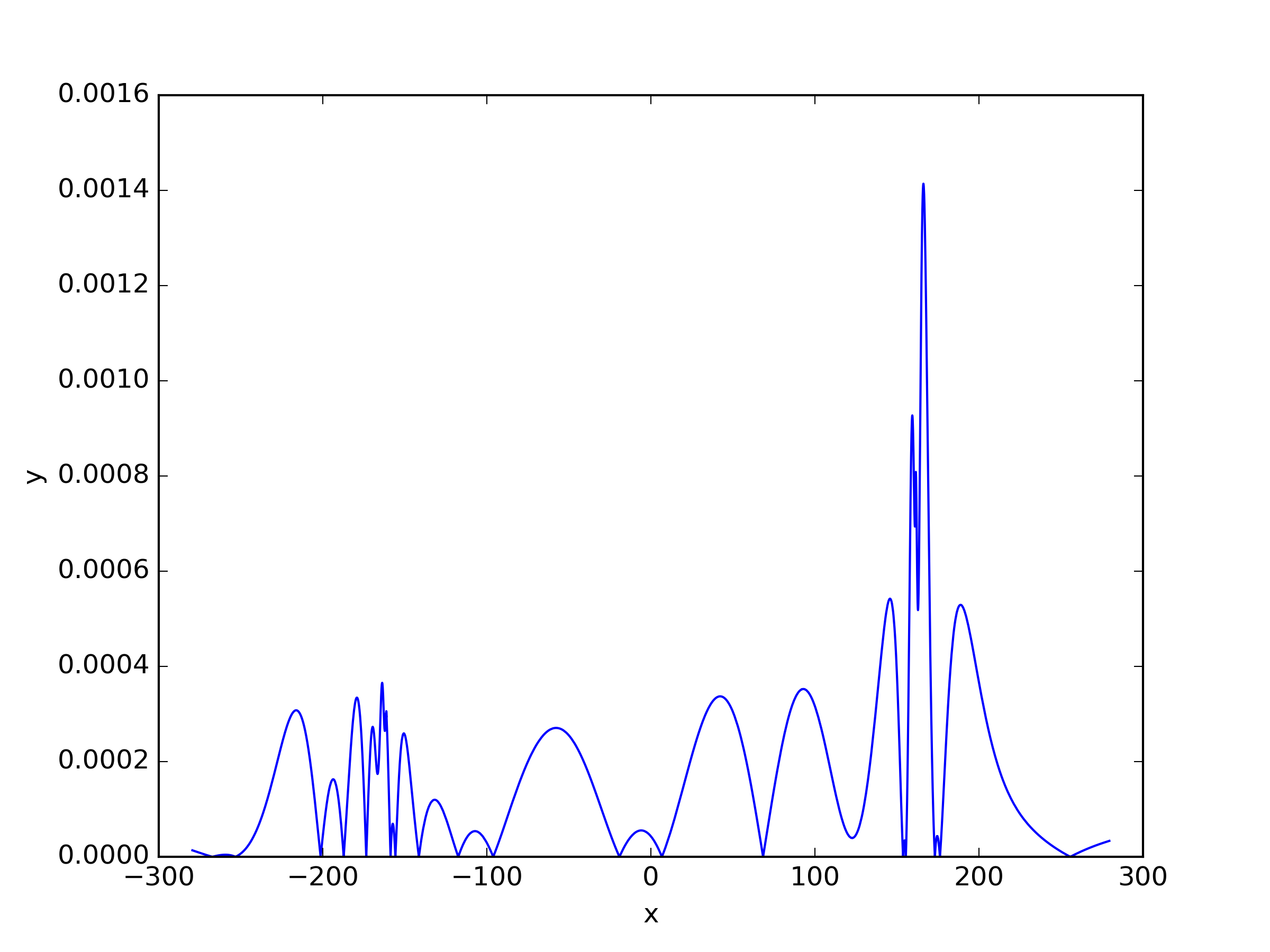}
		\end{minipage}
		\caption{Results of  1-D Swedish Ringhals-4 pressurized water reactor problem: Left: reference solution. Mid: neural network solution. Right: absolute error. First row: information of $\phi_1(x)$. Second row: information of $\phi_2(x)$.}
		\label{fig-slab-1D2R}
	\end{figure*}

	We applied three different residual loss functions to this one-dimensional case, and for $Loss_{De}$ \eqref{loss-res-3} and $Loss_{DI}$ \eqref{loss-res-2}, we added a shift $\sigma = 1$. When the shift term for $Loss_{De}$ and $Loss_{DI}$ is set to 0, it corresponds to the same loss function. The results of the solutions can be seen in Table \ref{tab-1D2G} and Figure \ref{fig-slab-1D2R}. Here $Loss_{IPM}$ refers to inverse power method residual loss function proposed in \cite{yang2023physics}, and we also apply this residual loss function in other problems. The expression of it is shown as \eqref{loss-res-1}. We observed that the inverse power method loss function \eqref{loss-res-1}, which was effective for solving single-group problems, did not perform well for multi-group problems. Surprisingly, $Loss_{De}$ showed the best performance in this case.
	
	\begin{equation}
		\begin{aligned}
			Loss_{IPM} = & \sum\limits_{\boldsymbol{x} \in \mathbf{X}_{res}} \big(-\nabla \cdot D_1 \nabla \phi_1^{NN(i)} + (\Sigma_{a,1}+\Sigma_{1\rightarrow 2})\phi_1^{NN(i)} -\lambda(\nu\Sigma_{f,1}\phi_1^{NN(i-1)}+\nu\Sigma_{f,2}\phi_2^{NN(i-1)}) \big)^2_{\boldsymbol{x}}\\
			& + \big(-\nabla \cdot D_2 \nabla \phi_2^{NN(i)} + \Sigma_{a,2}\phi_2^{NN(i)}-\Sigma_{1\rightarrow 2}\phi_1^{NN(i)}\big)^2_{\boldsymbol{x}},
		\end{aligned}
		\label{loss-res-1}
	\end{equation}
	
	\begin{table*}[h]
		\caption{Recording numerical results of  1-D Swedish Ringhals-4 pressurized water reactor problem with different residual loss function applied.}
		\centering
		\begin{tabular}{cccccccc}
			\hline
			& Shift & $\mathbf{E}_R(k_{\text{eff}})$ & $\mathbf{E}_{R,\infty}(\phi_1)$ & $\mathbf{E}_{R,2}(\phi_1)$ & $\mathbf{E}_{R,\infty}(\phi_2)$ & $\mathbf{E}_{R,2}(\phi_2)$ \\
			\hline
			$Loss_{IPM}$ & 0 & 1.7623e-04 & 4.8424e-02 & 2.9796e-02 & 1.7779e-01 & 5.0444e-02 \\
			\hline
			$Loss_{De}$ & 0 & 1.7743e-04 & 4.5991e-03 & 2.3375e-03 & 2.3234e-02 & 6.8866e-03 \\
			\hline
			$Loss_{DI}$ & 1 & 7.2184e-04 & 2.6505e-02 & 1.6108e-02 & 1.4934e-01 & 3.8318e-02 \\
			\hline
			$Loss_{De}$ & 1 & \pmb{1.0351e-04} & \pmb{2.3804e-03} & \pmb{1.7113e-03} & \pmb{9.5453e-03} & \pmb{2.9083e-03} \\
			\hline
		\end{tabular}
		\label{tab-1D2G}
	\end{table*}
	
	\subsection{2-D Problem}
	\subsubsection{TWIGL}
\begin{figure*}[h]
		\centering
		\includegraphics[width=\linewidth]{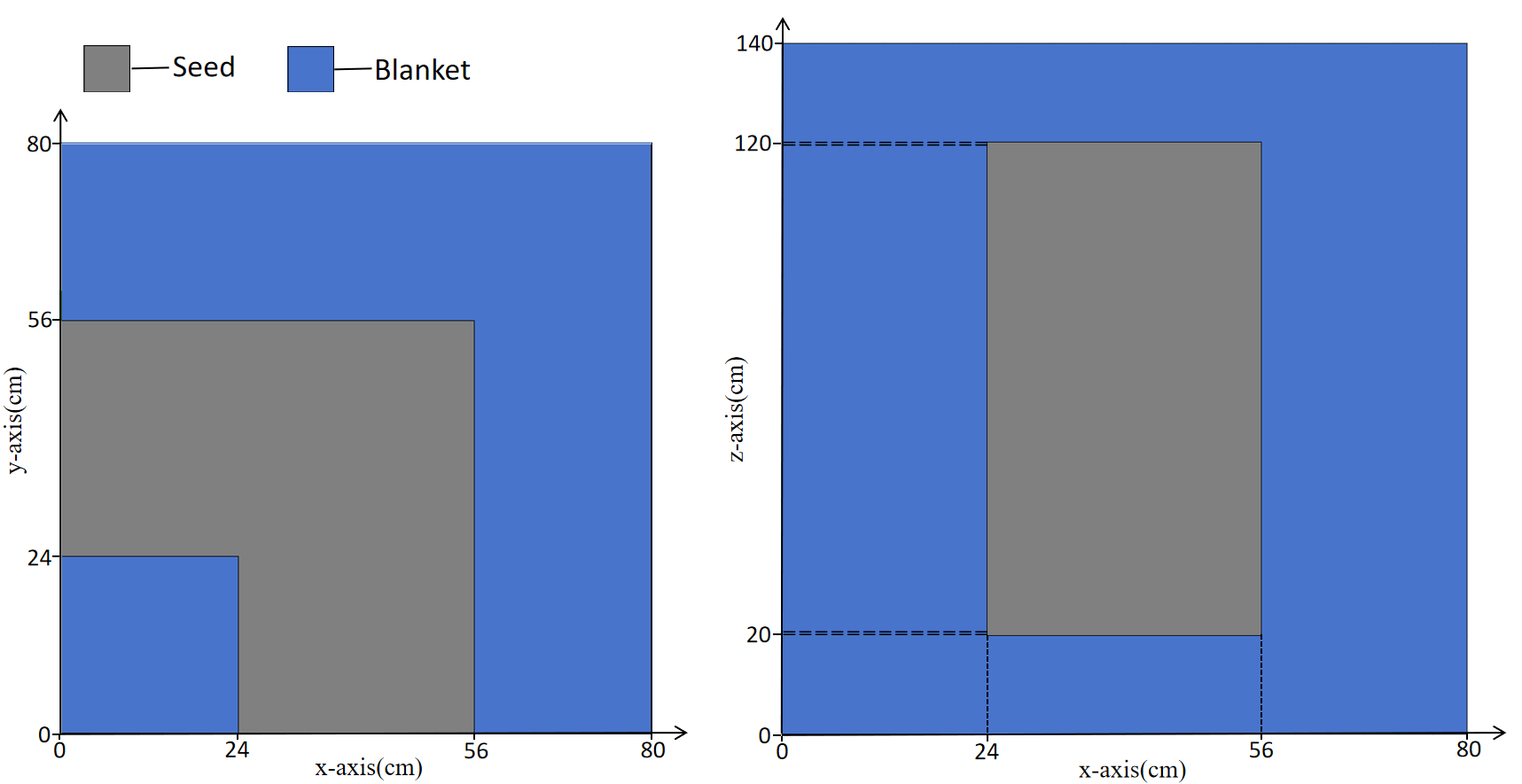}
		\caption{Computational domain of TWIGL problem. Grey and blue region denote seed and blanket materials, respectively. Left: first quadrant of $x-y$ plane at $z$ = 70 $\text{cm}$. Right: $x-z$ plane at $y$ = 0 $\text{cm}$.}
		\label{fig-TWIGL}
	\end{figure*}

	The TWIGL model is a square-shaped nuclear reactor core with dimensions of 160 centimeters for both length and width. For the critical case, the geometric region consists of only two materials: seed and blanket. By exploiting the symmetry of the reactor core, we can reduce the computational domain of the partial differential equation to one-fourth of its original size (Figure \ref{fig-TWIGL}). At the real outer boundary, homogeneous Dirichlet boundary conditions are applied, while homogeneous Neumann boundary conditions are set at the symmetric boundaries. Therefore, we will apply the Neumann boundary loss function \eqref{loss-neumann} at the boundaries x=0 and y=0, and the Dirichlet boundary loss function \eqref{loss-diri} at the boundaries x=80 and y=80. The corresponding coefficients are presented in the Table \ref{tab-coef-2D-TWIGL}. 
	
	We solved the problem using the finite element method with 320 points on the boundaries with intervals $\Delta x = 1$ and $\Delta y = 1$, and then constructed the internal element mesh using the Delaunay algorithm. This critical eigenvalue solved by FreeFEM is $k_{\text{eff}}^{FF} = 0.9133$. {For our neural network, we have selected a total of 6,399 residual points and 162 interface points.} We chose 4 residual blocks, with 32 neurons per layer in each residual block. The Adam optimizer with a learning rate of $0.001$ was used for $100,000$ epochs optimization during the training process.
	
	\begin{table*}[h]
		\caption{Coefficients of different regions for TWIGL problem}
		\centering
		\begin{tabular}{|c|c|c|c|c|c|c|c|c|}
			\hline
			Region & $D_1$ & $D_2$ & $\Sigma_{a,1}$ & $\Sigma_{a,2}$ & $\Sigma_{1\rightarrow2}$ & $\nu\Sigma_{f,1}$ & $\nu\Sigma_{f,2}$\\
			& $(cm)$ & $(cm)$ & $(cm^{-1})$ & $(cm^{-1})$ & $(cm^{-1})$ & $(cm^{-1})$ & $(cm^{-1})$\\
			\hline
			Seed & 1.4 & 0.4 & 0.01 & 0.15 & 0.01 & 0.007 & 0.2\\
			\hline
			Blanket & 1.3 & 0.5 & 0.008 & 0.05 & 0.01 & 0.003 & 0.06\\
			\hline
		\end{tabular}
		\label{tab-coef-2D-TWIGL}
	\end{table*}
	
	\begin{figure*}[h]
		\begin{minipage}{0.32\textwidth}
			\centering
			\includegraphics[width=5.5cm]{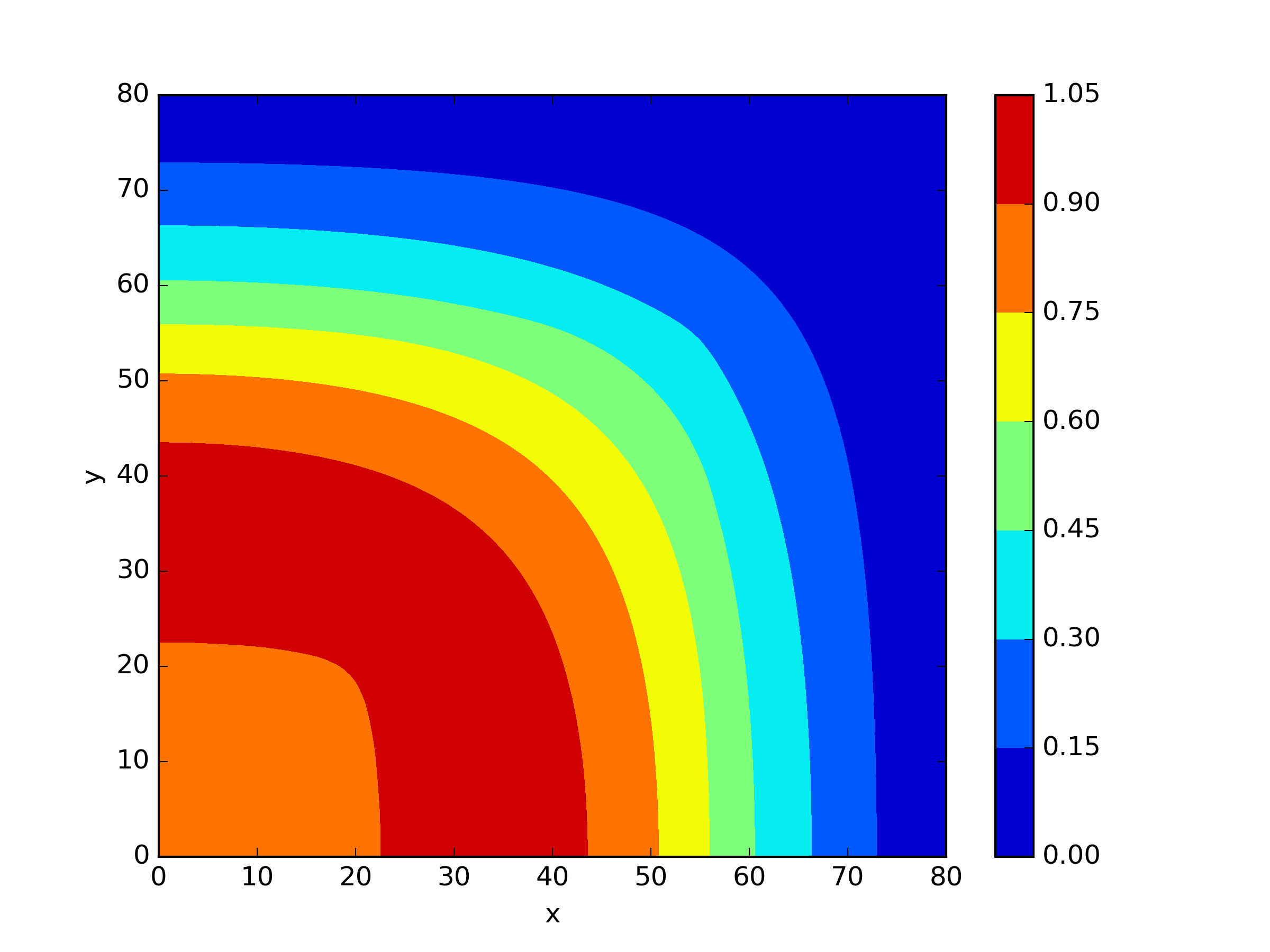}
		\end{minipage}
		\begin{minipage}{0.32\textwidth}
			\centering
			\includegraphics[width=5.5cm]{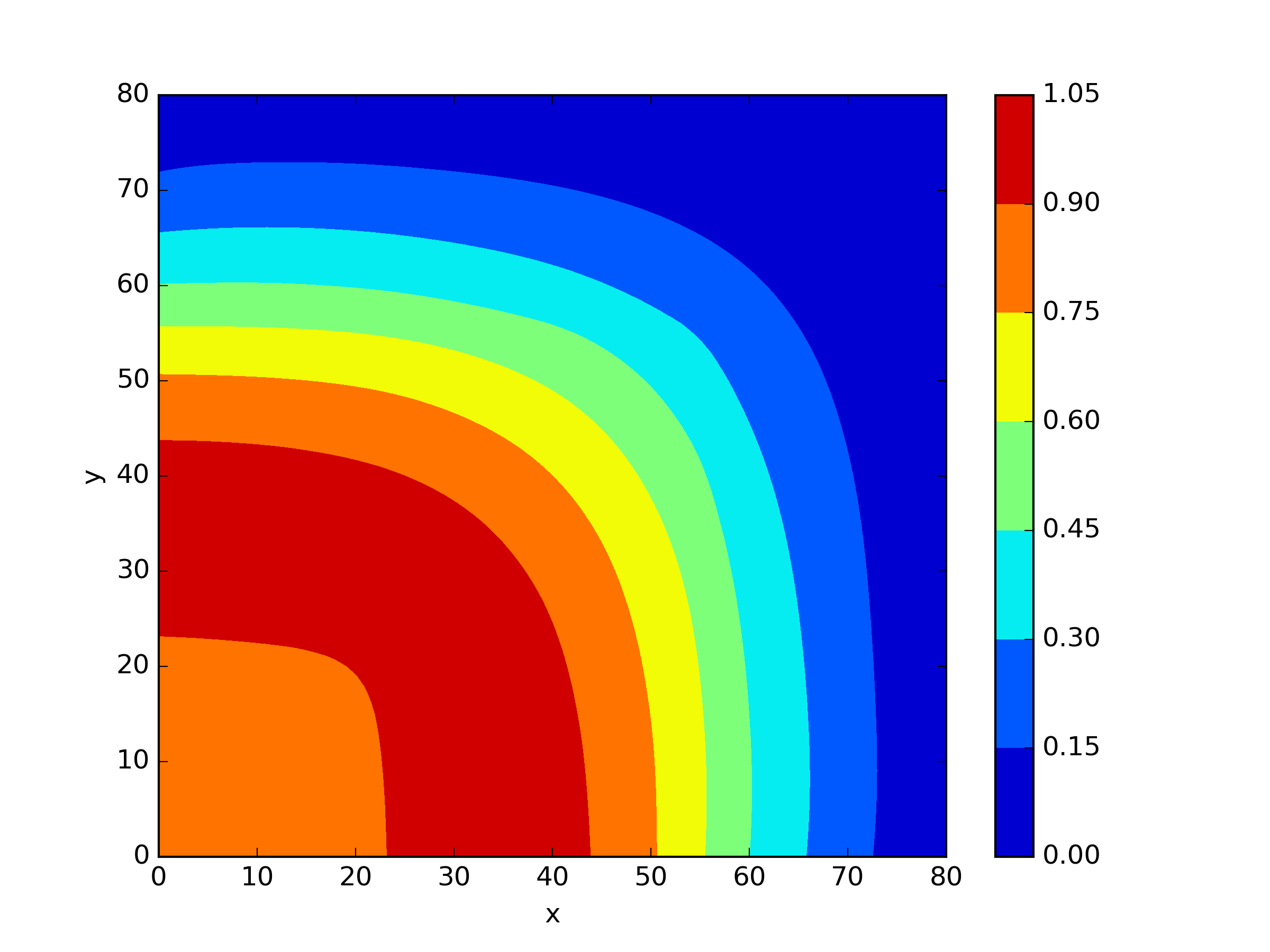}
		\end{minipage}
		\begin{minipage}{0.32\textwidth}
			\centering
			\includegraphics[width=5.5cm]{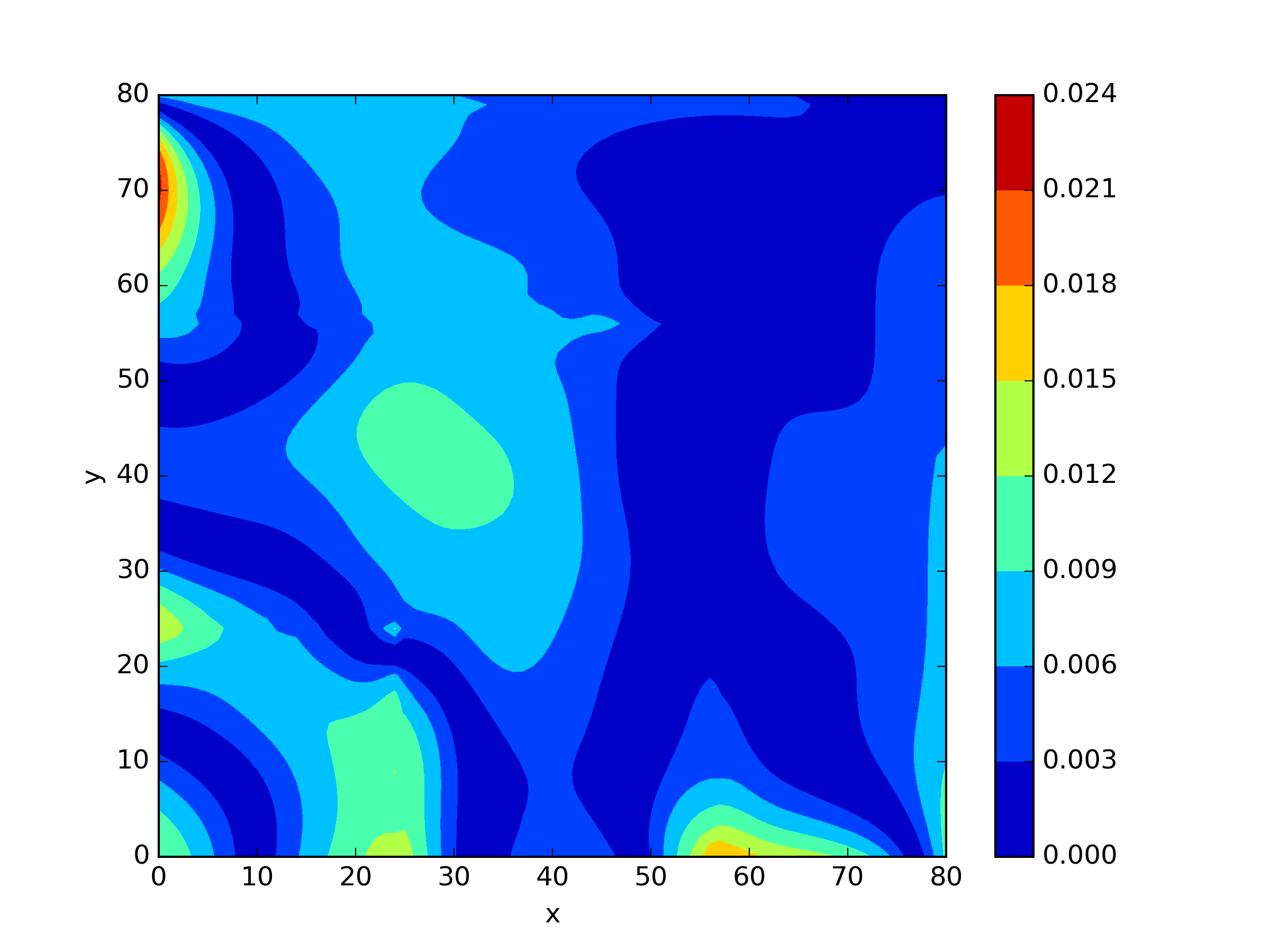}
		\end{minipage}
		
		\begin{minipage}{0.32\textwidth}
			\centering
			\includegraphics[width=5.5cm]{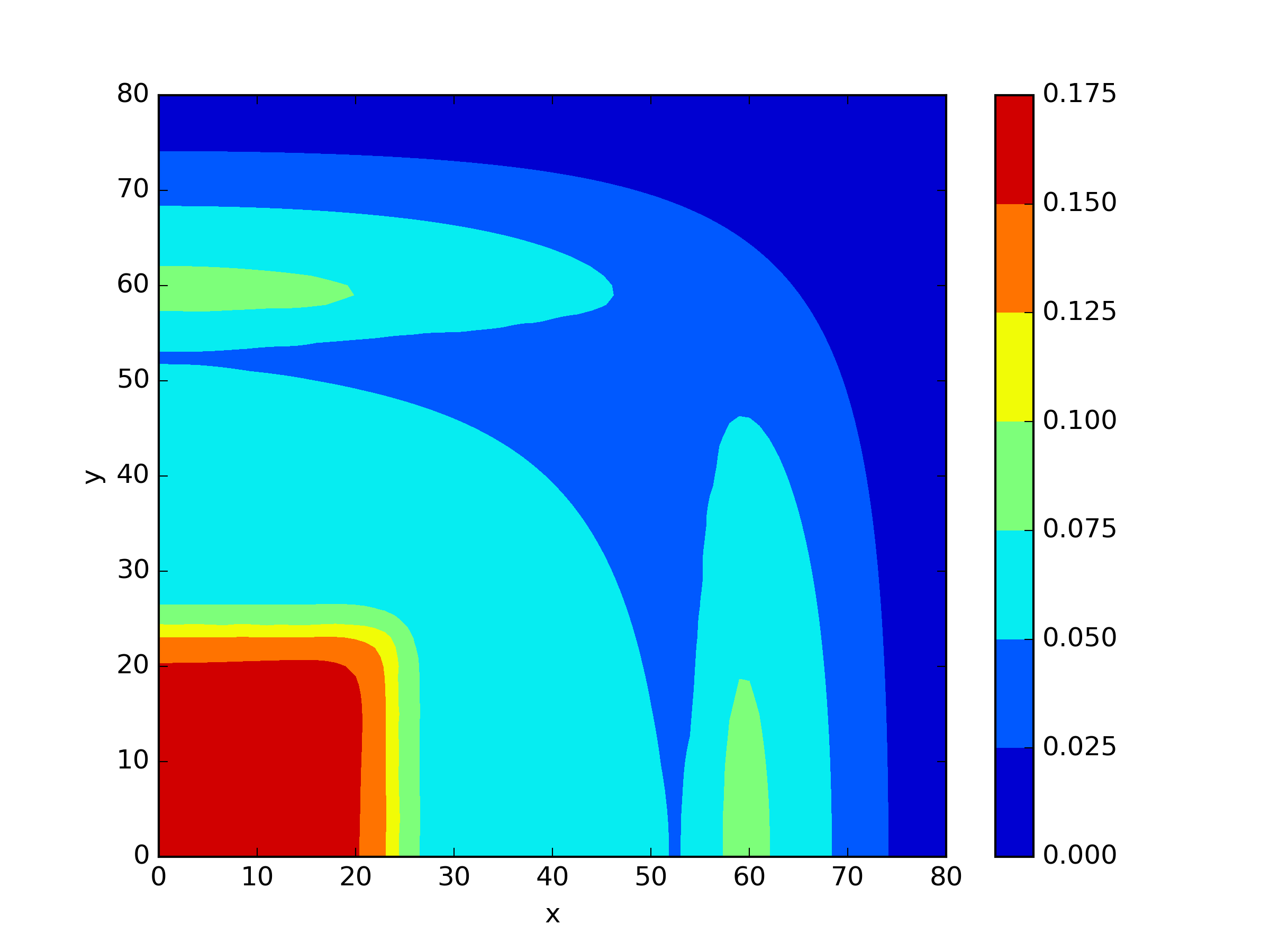}
		\end{minipage}
		\begin{minipage}{0.32\textwidth}
			\centering
			\includegraphics[width=5.5cm]{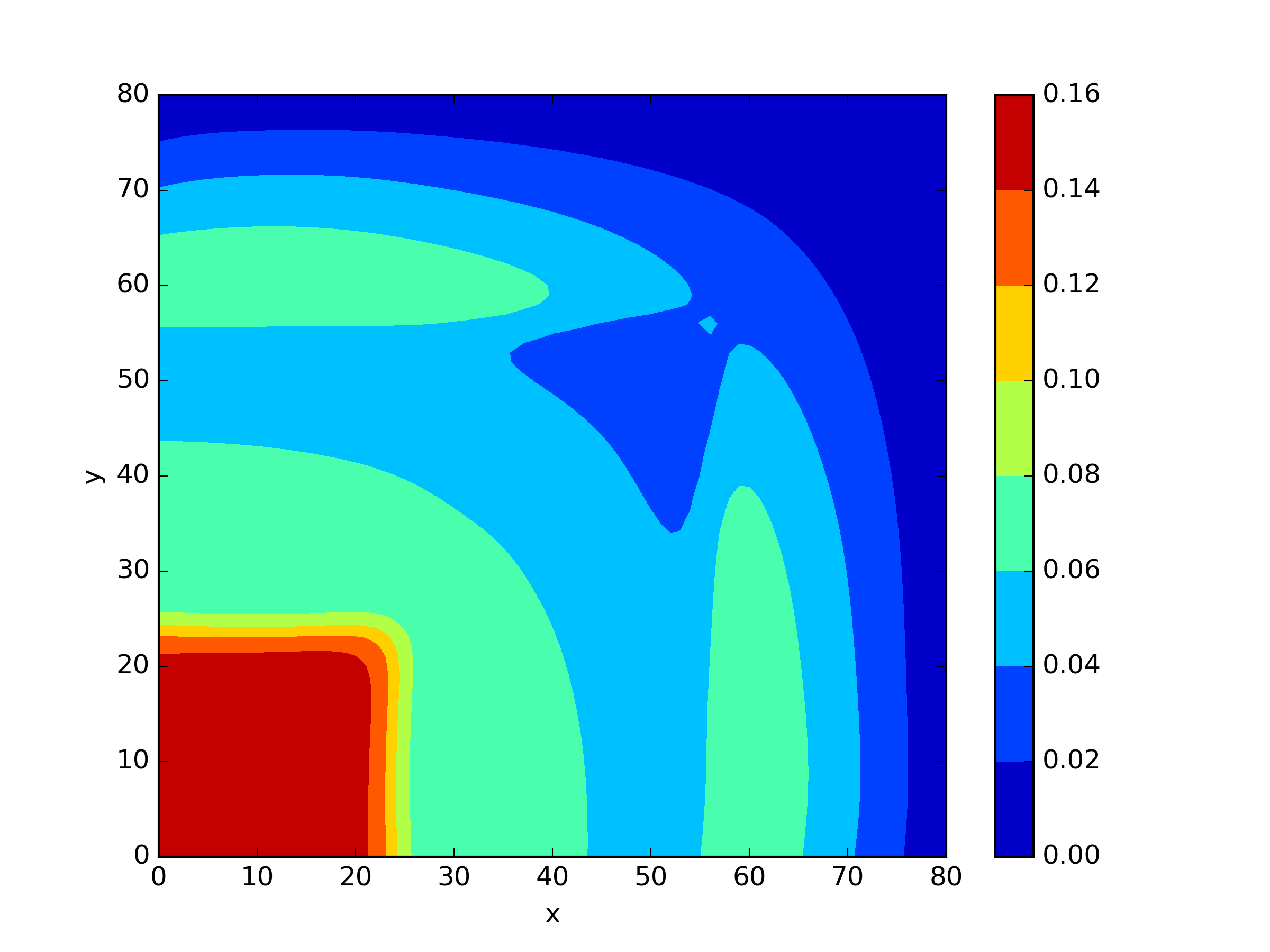}
		\end{minipage}
		\begin{minipage}{0.32\textwidth}
			\centering
			\includegraphics[width=5.5cm]{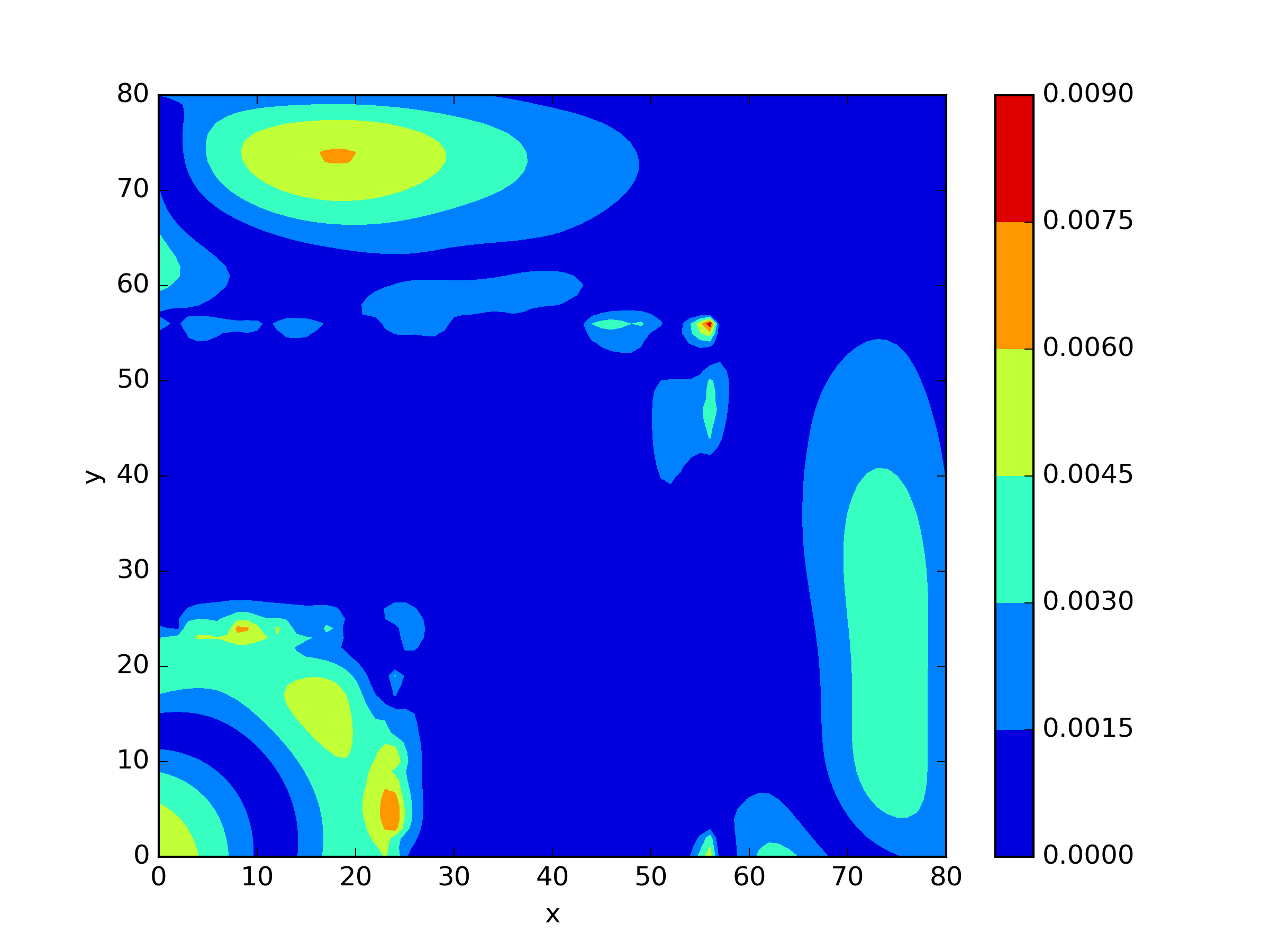}
		\end{minipage}
		\caption{Results of 2-D TWIGL problem: Left: reference solution. Mid: neural network solution. Right: absolute error. First row: information of $\phi_1(x, y)$. Second row: information of $\phi_2(x, y)$.}
		\label{fig-TWIGL-2D2R}
	\end{figure*}
	
	We applied $Loss_{De}, Loss_{DI}$ and $Loss_{IPM}$ mentioned earlier to this scenario. For $Loss_{IPM}$, we only tested the approach without shift. For  $Loss_{DI}$ and $Loss_{De}$, we attempted to use a shift value of 1. The results are presented in Table \ref{tab-TWIGL-2D2G} and Figure \ref{fig-TWIGL-2D2R}. We found that $Loss_{IPM}$ did not exhibit any significant disadvantages in solving this problem. However, it can be observed that the methods with a shift in $Loss_{DI}$ and $Loss_{De}$ achieved slightly higher accuracy compared to the method without shift.
	
	\begin{table*}[h]
		\caption{Recording numerical results for 2-D TWIGL problem with different residual loss function applied.}
		\centering
		\begin{tabular}{cccccccc}
			\hline
			& Shift & $\mathbf{E}_R(k_{\text{eff}})$ & $\mathbf{E}_{R,\infty}(\phi_1)$ & $\mathbf{E}_{R,2}(\phi_1)$ & $\mathbf{E}_{R,\infty}(\phi_2)$ & $\mathbf{E}_{R,2}(\phi_2)$ \\
			\hline
			$Loss_{IPM}$ &  0 & 4.6363e-03 & 3.4499e-02 & 1.3948e-02 & \pmb{5.1534e-02} & \pmb{2.4732e-02} \\
			\hline
			$Loss_{De}$ &  0 & 4.3026e-03 & 3.4864e-02 & 1.6327e-02 & 5.6210e-02 & 2.7221e-02 \\
			\hline
			$Loss_{DI}$ &  1 & \pmb{3.3199e-03} & \pmb{2.1937e-02} & \pmb{8.9844e-03} & 5.3010e-02 & 2.9896e-02 \\
			\hline
			$Loss_{De}$ &  1 & 3.9902e-03 & 3.1621e-02 & 1.4112e-02 & 5.4007e-02 & 2.6699e-02 \\
			\hline
		\end{tabular}
		\label{tab-TWIGL-2D2G}
	\end{table*}
	
\begin{figure*}[h]
		\centering
		\begin{minipage}[b]{0.32\linewidth}
			\subfloat[Relative error of flux]{
				\includegraphics[width=5.15cm]{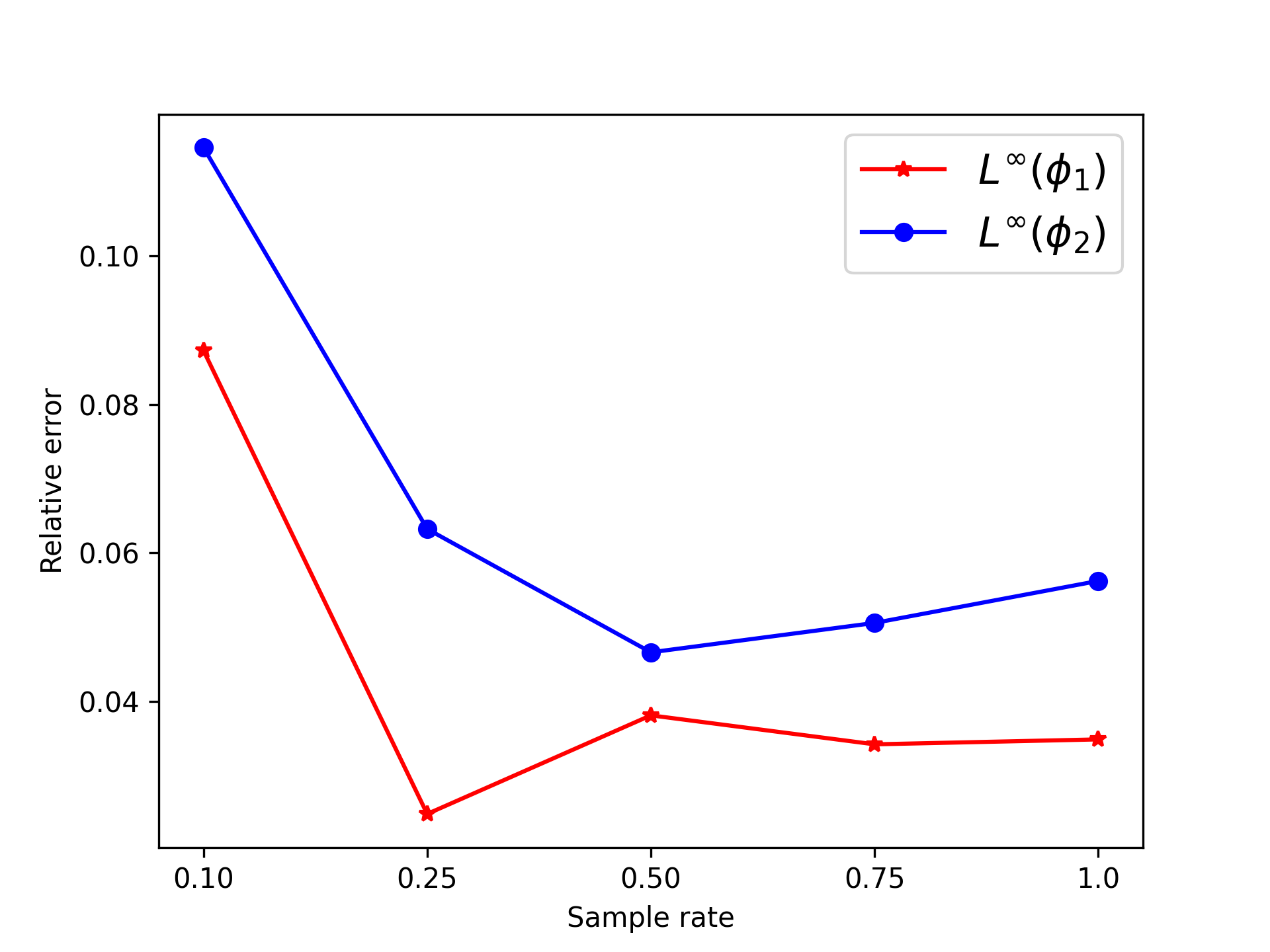}}
		\end{minipage}
		\begin{minipage}[b]{0.32\linewidth}
			\subfloat[Relative error of flux]{
				\includegraphics[width=5.15cm]{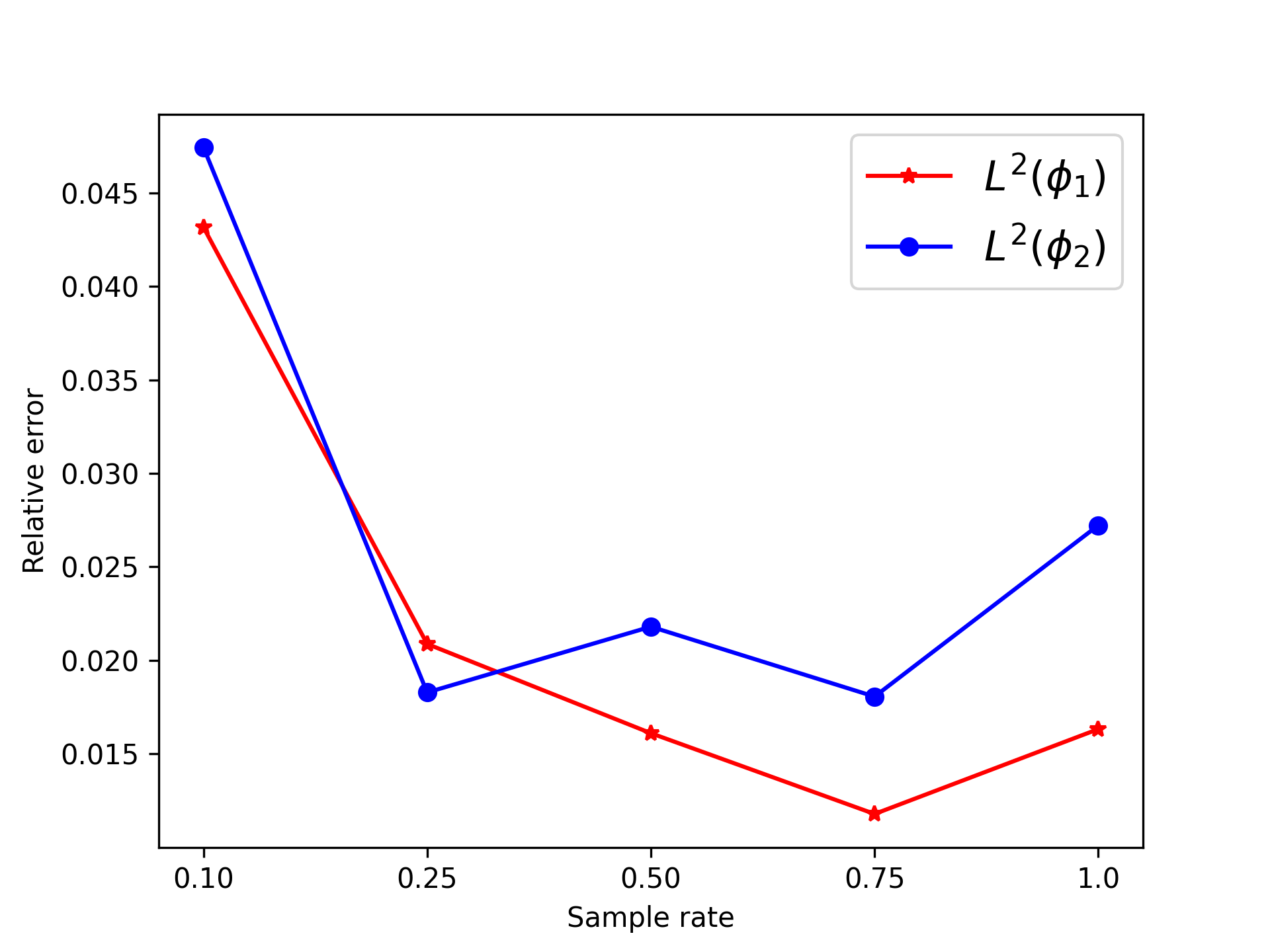}}
		\end{minipage}
		\begin{minipage}[b]{0.32\linewidth}
			\subfloat[Relative error of eigenvalue]{
				\includegraphics[width=5.15cm]{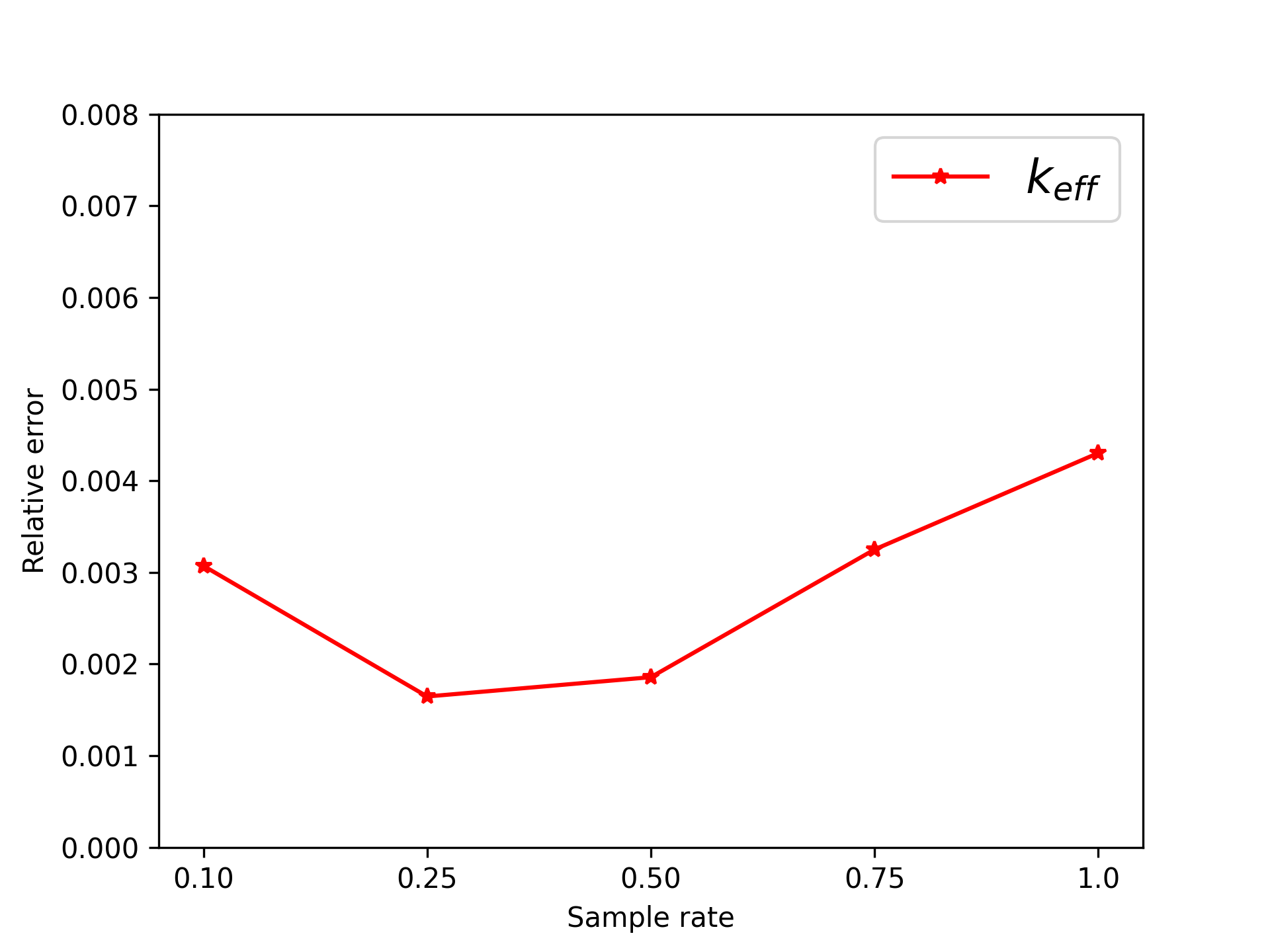}}
		\end{minipage}
		\caption{2-D TWIGL problem: relative error of neural network solution at different sample rate.}
		\label{pic-TWIGL-samplerate}
	\end{figure*}

\begin{figure*}[h]
		\centering
		\begin{minipage}[b]{0.49\linewidth}
			\subfloat[Relative $L_{\infty}$ error of $\phi_1$]{
				\includegraphics[width=7cm]{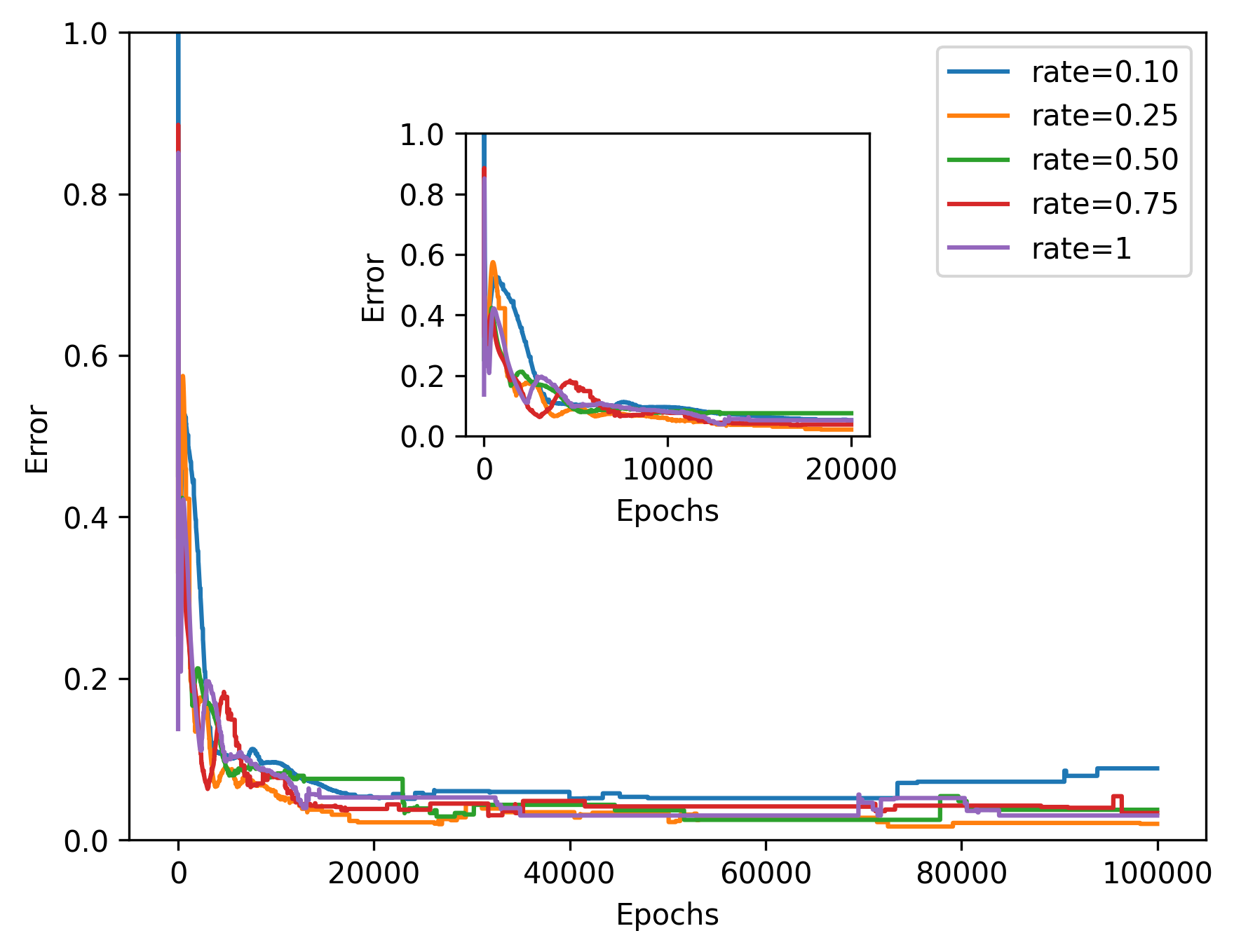}}
		\end{minipage}
		\begin{minipage}[b]{0.49\linewidth}
			\subfloat[Relative $L_{2}$ error of $\phi_1$]{
				\includegraphics[width=7cm]{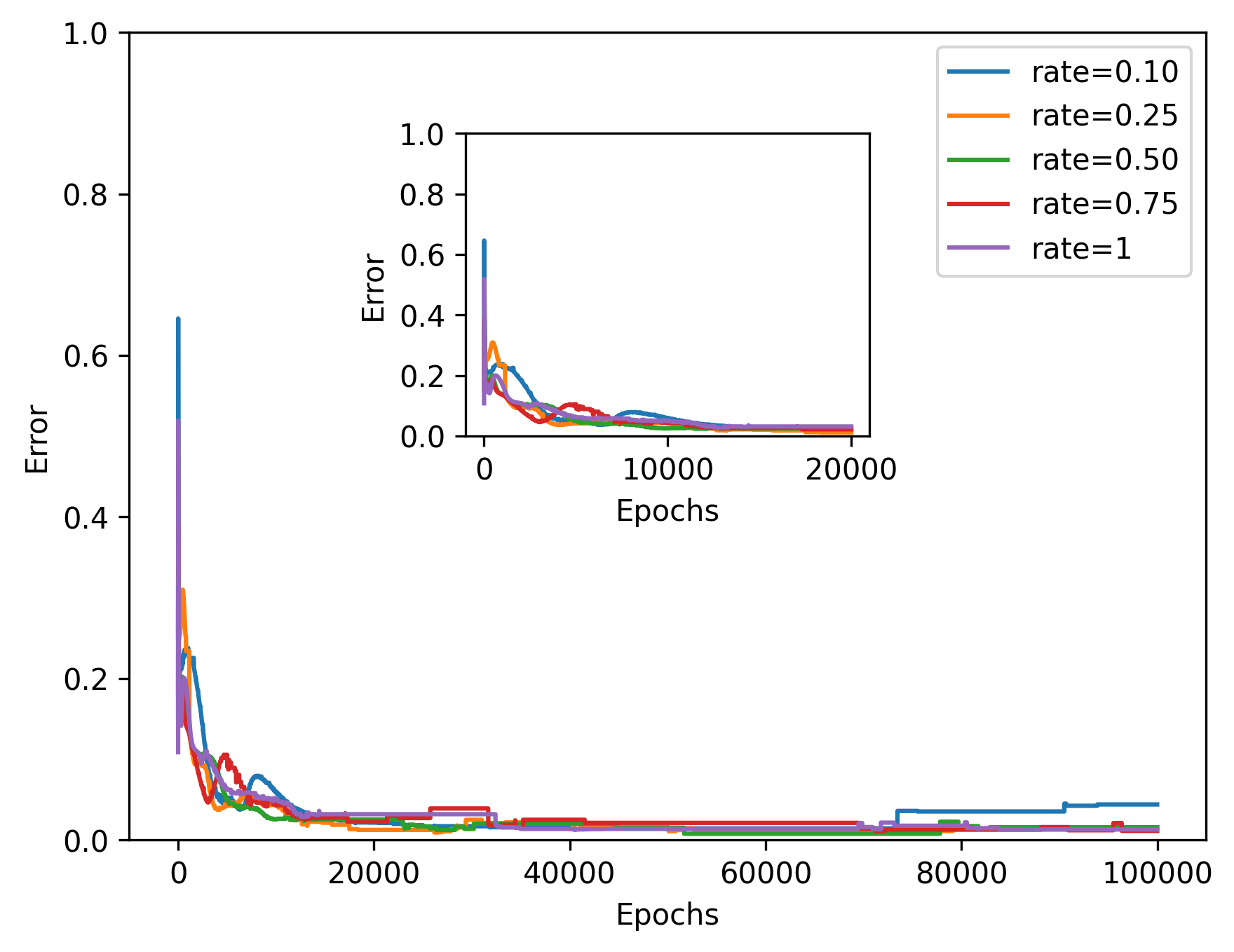}}
		\end{minipage}
		\begin{minipage}[b]{0.49\linewidth}
			\subfloat[Relative $L_{\infty}$ error of $\phi_2$]{
				\includegraphics[width=7cm]{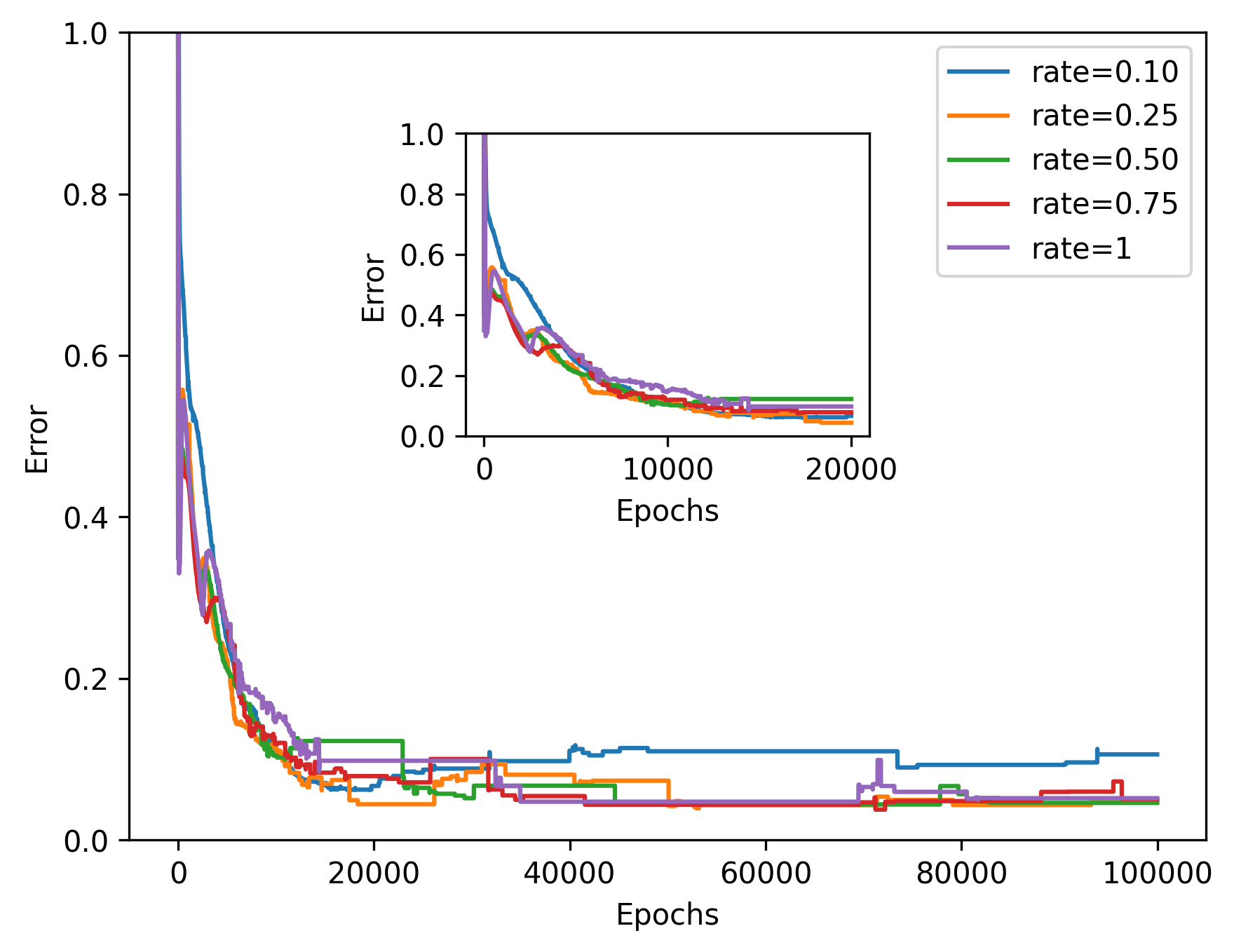}}
		\end{minipage}
		\begin{minipage}[b]{0.49\linewidth}
			\subfloat[Relative $L_{2}$ error of $\phi_2$]{
				\includegraphics[width=7cm]{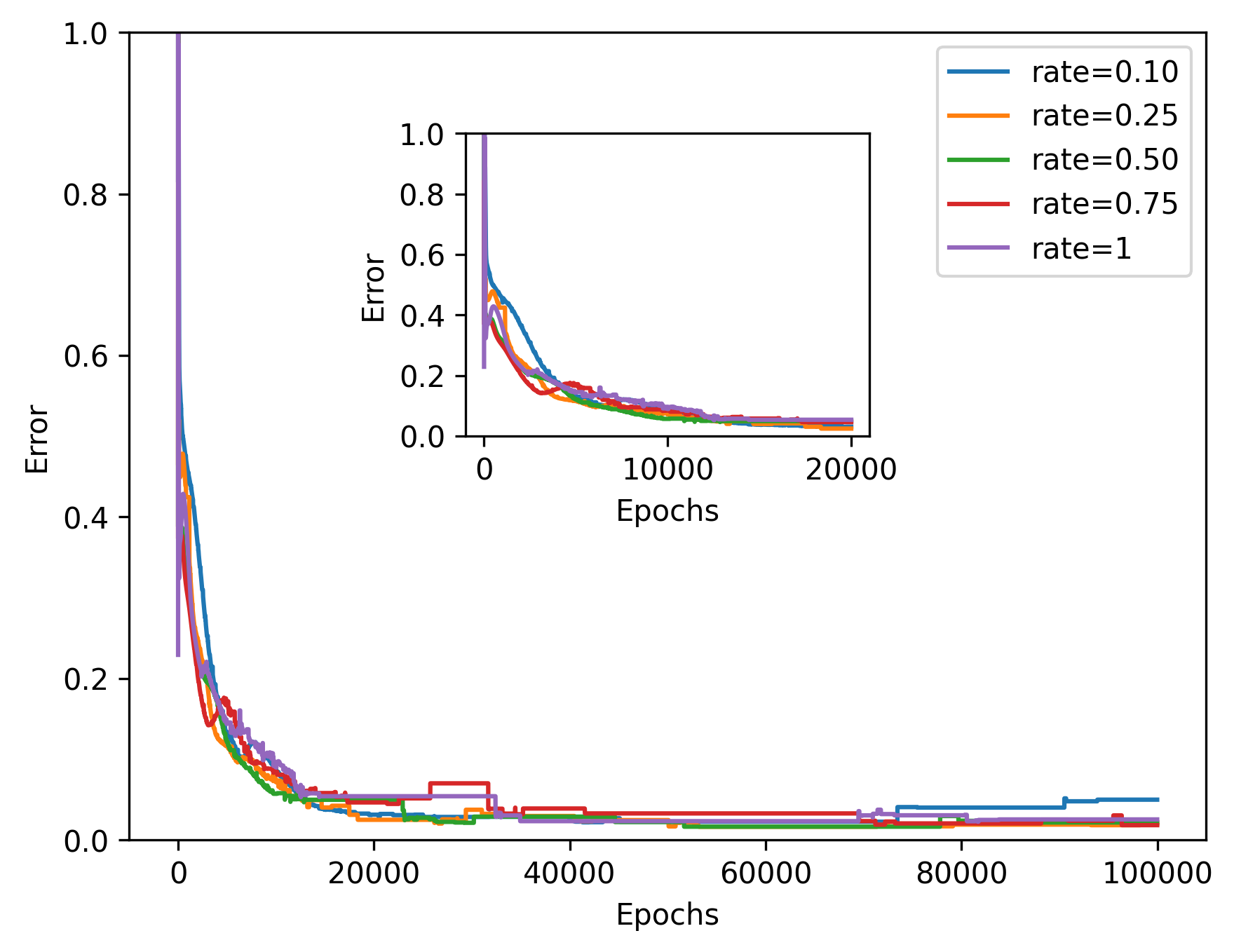}}
		\end{minipage}
		\caption{2-D TWIGL problem: the first row shows relative error of $\phi_1$ during training process; the second row shows relative error of $\phi_2$ during training process. Black, red, green, blue, and purple colors represent the relative error curves at different sampling rates of 0.1, 0.25, 0.5, 0.75, and 1, respectively.}
		\label{pic-TWIGL-samplerate-train}
	\end{figure*}

	In the two-dimensional problem of TWIGL, we conducted a series of experiments regarding the number of sampling points and the solution accuracy. In these experiments, we initially used 6561 sampling points as the complete set of sampling points. In the mentioned paper, the sampling rates of 0.1, 0.25, 0.5, and 0.75 for residual points refer to using only $10\%, 25\%, 50\%$, and $75\%$ of the residual points as training data, respectively. For example, a sampling rate of 0.5 means that only half of the residual points, which amounts to 3,281 points, are used for residual loss function training. We found that a significant reduction in the number of sampling points, such as a sampling rate of 0.1, leads to a noticeable decrease in accuracy (Figure \ref{pic-TWIGL-samplerate}). Reducing the sampling rate will invariably result in slower convergence of the neural network (Figure \ref{pic-TWIGL-samplerate-train}).
	
\begin{figure*}[h]
		\centering
		\begin{minipage}[b]{0.49\linewidth}
			\subfloat[Relative $L_{\infty}$ error of $\phi_1$]{
				\includegraphics[width=8cm]{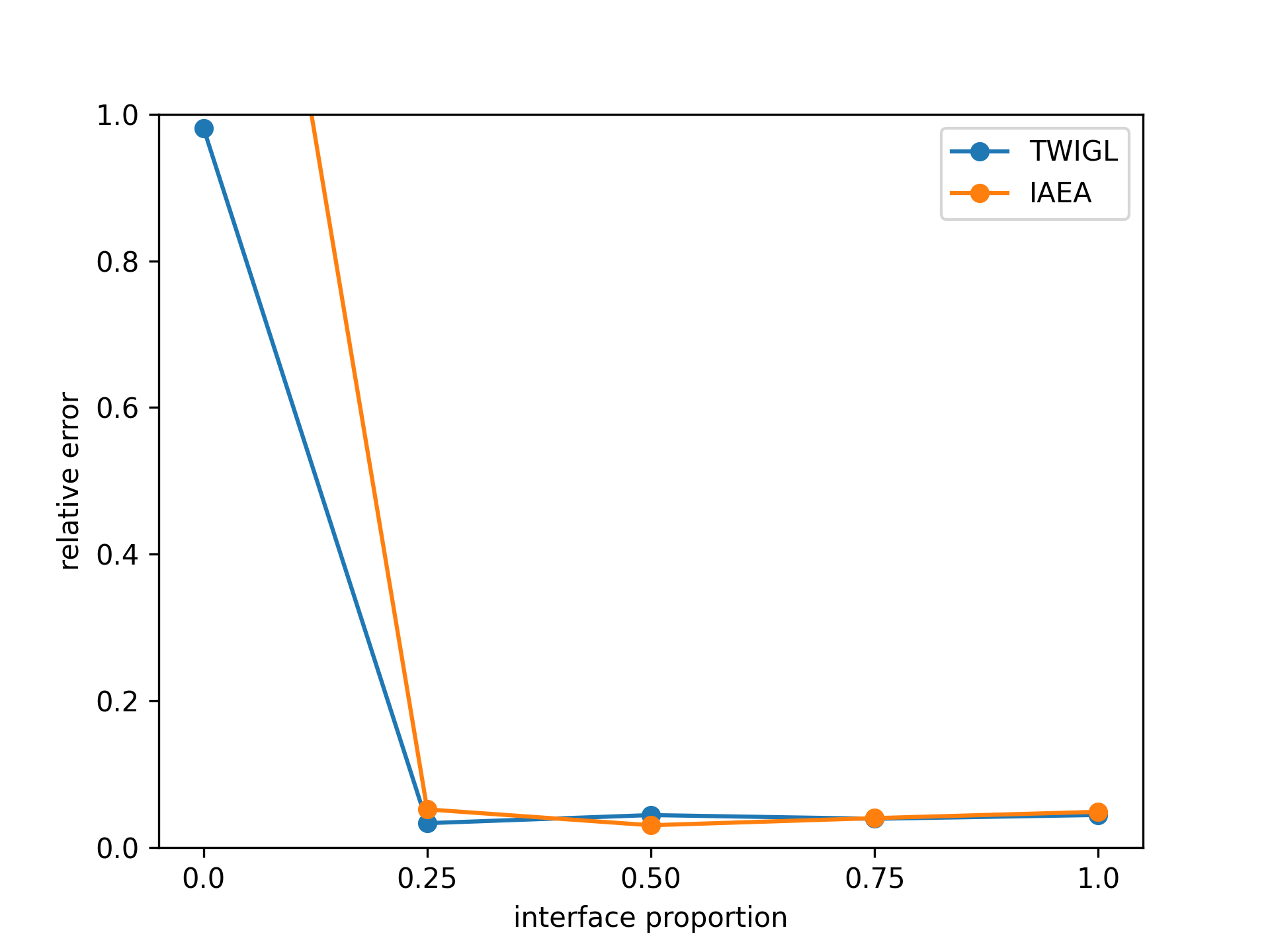}}
		\end{minipage}
		\begin{minipage}[b]{0.49\linewidth}
			\subfloat[Relative $L_{2}$ error of $\phi_1$]{
				\includegraphics[width=8cm]{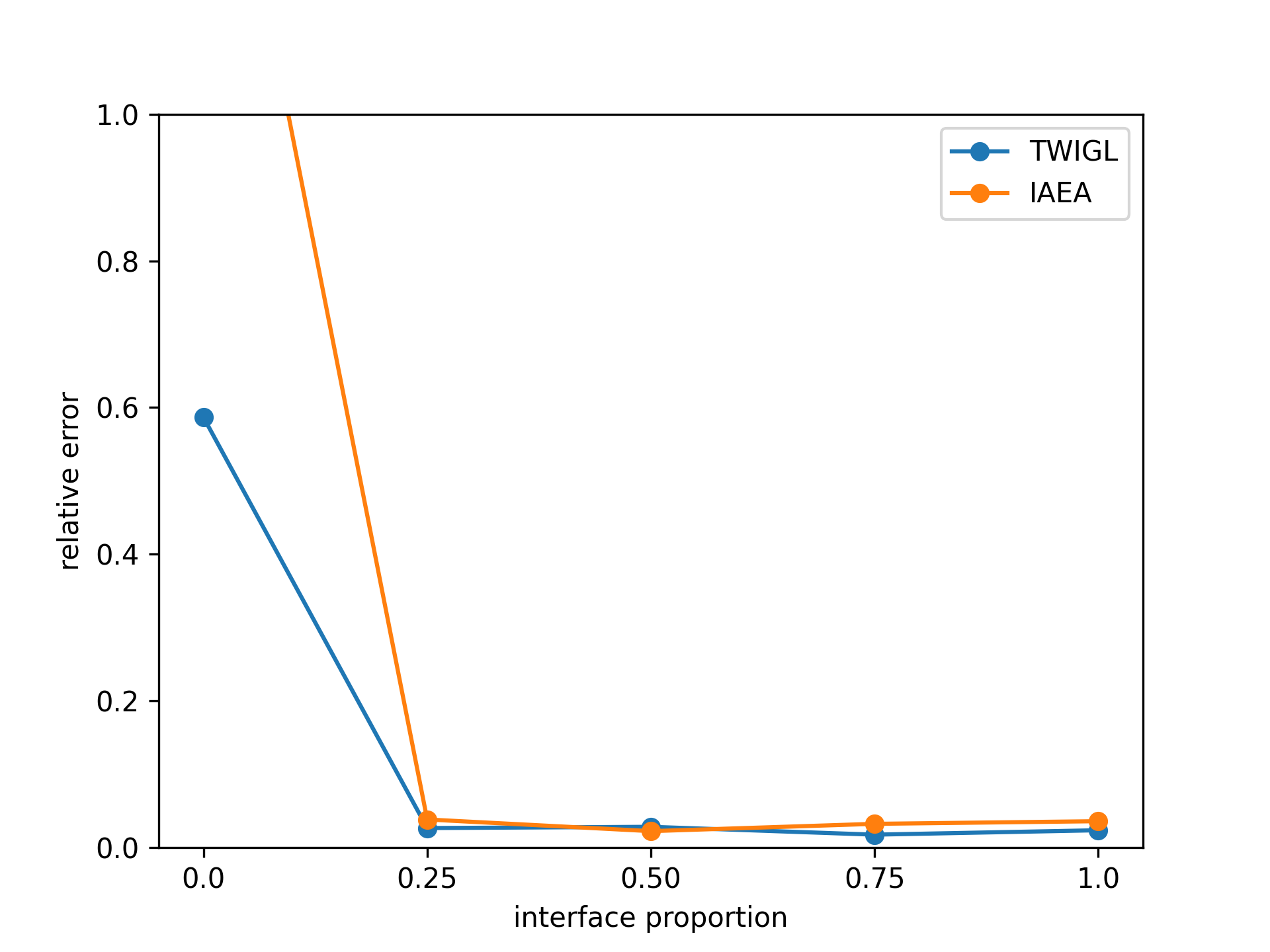}}
		\end{minipage}
		\begin{minipage}[b]{0.49\linewidth}
			\subfloat[Relative $L_{\infty}$ error of $\phi_2$]{
				\includegraphics[width=8cm]{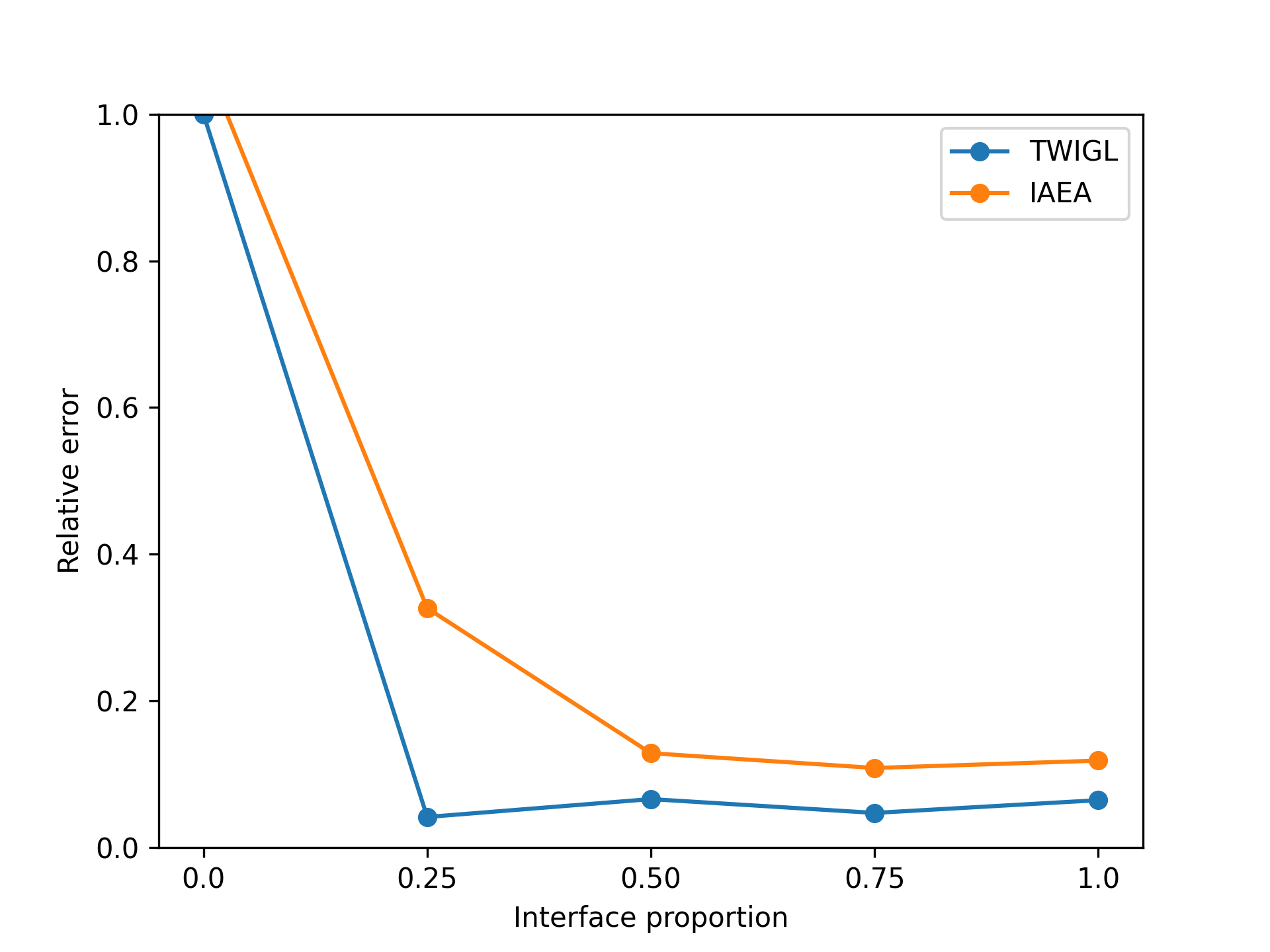}}
		\end{minipage}
		\begin{minipage}[b]{0.49\linewidth}
			\subfloat[Relative $L_{2}$ error of $\phi_2$]{
				\includegraphics[width=8cm]{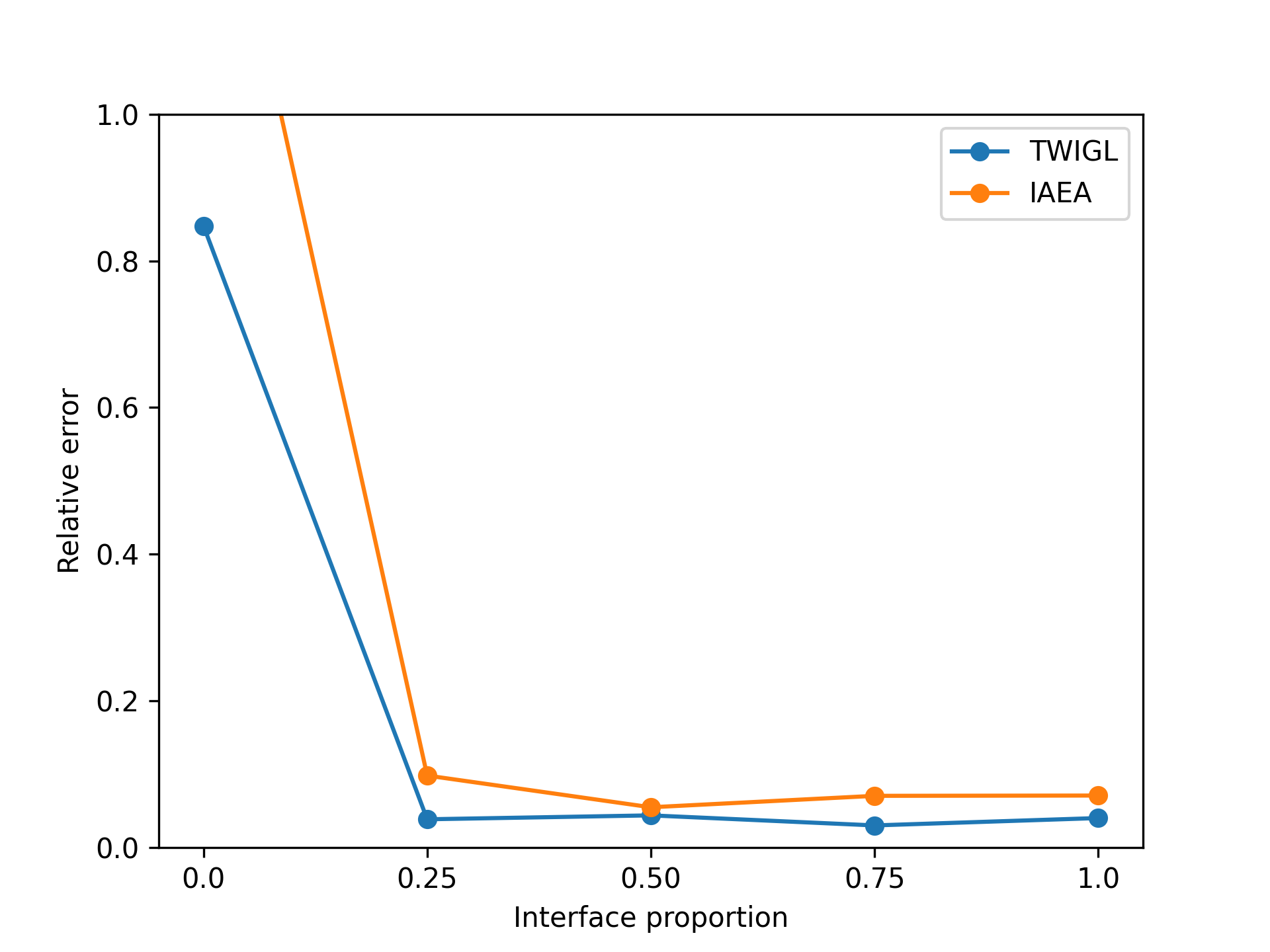}}
		\end{minipage}
		\caption{Experiment results of different proportion of interface points used in training process. The blue and orange solid line denote results from 2-D TWIGL problem and 2-D IAEA problem, respectively. The first row shows relative error of $\phi_1$ during training process; the second row shows relative error of $\phi_2$ during training process.}
		\label{pic-interface}
	\end{figure*}

\begin{figure*}[h]
		\centering
		\includegraphics[width=8cm]{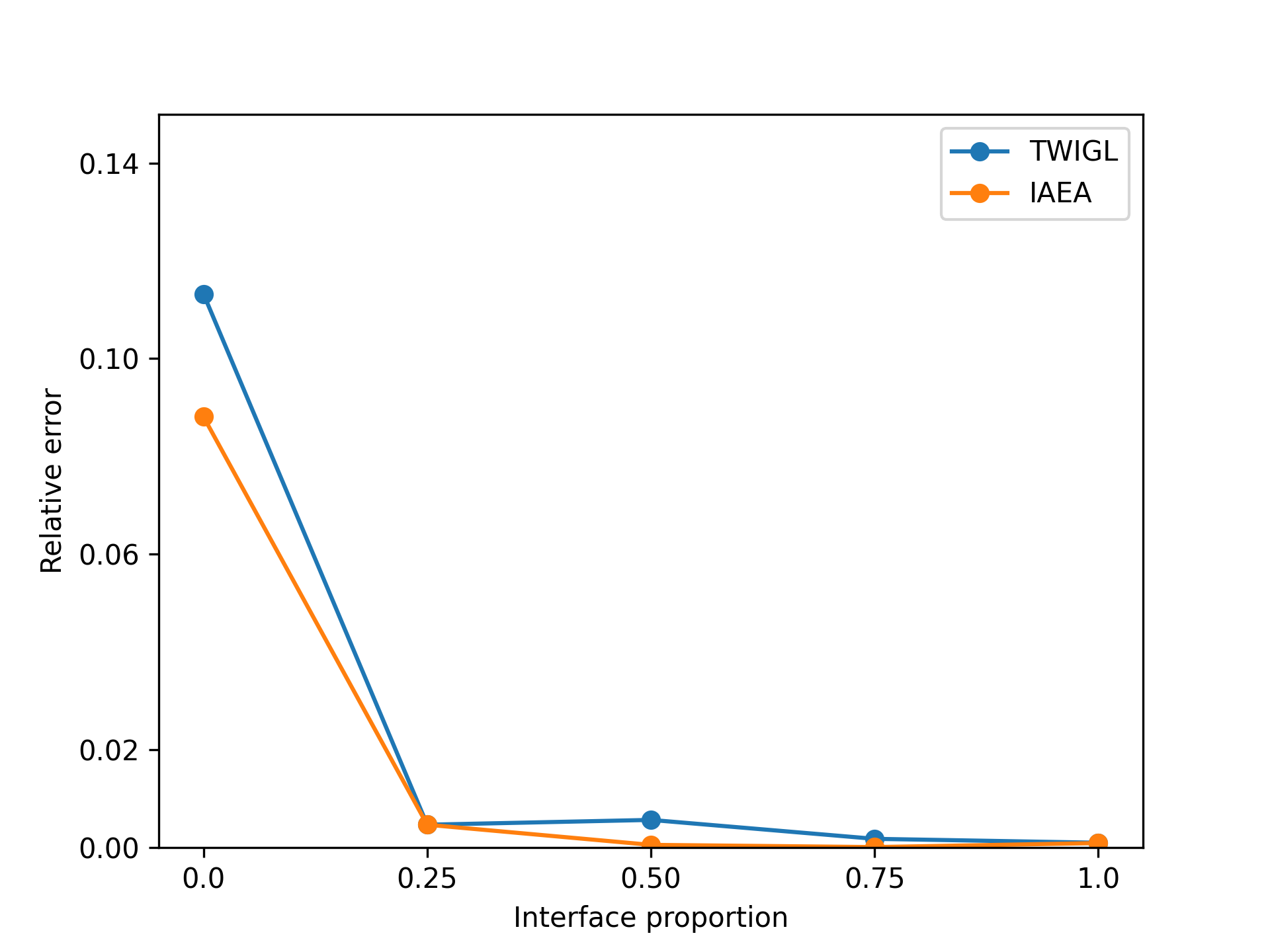}
		\caption{The relative error of $k_{\text{eff}}$ under different proportion of interface points used in training process, where blue represents 2-D TWIGL problem and yellow represents 2-D IAEA problem.}
		\label{pic-int-keff}
	\end{figure*}

	Additionally, we conducted another series of experiments to examine the impact of introducing interface conditions on the numerical experiments. In this experiment, while keeping the number of training points constant, we reduced the number of points corresponding to the interface conditions. As a result, the training points for the equation loss term increased. We tested the numerical examples using models ranging from fully incorporating interface conditions to not including them at all. The results of this experiment indicate that neural network completely fails to solve the problem when interface conditions are not introduced (Figure \ref{pic-interface} and \ref{pic-int-keff}).
	
	\subsubsection{IAEA}
\begin{figure*}[h]
		\centering
		\includegraphics[width=\linewidth]{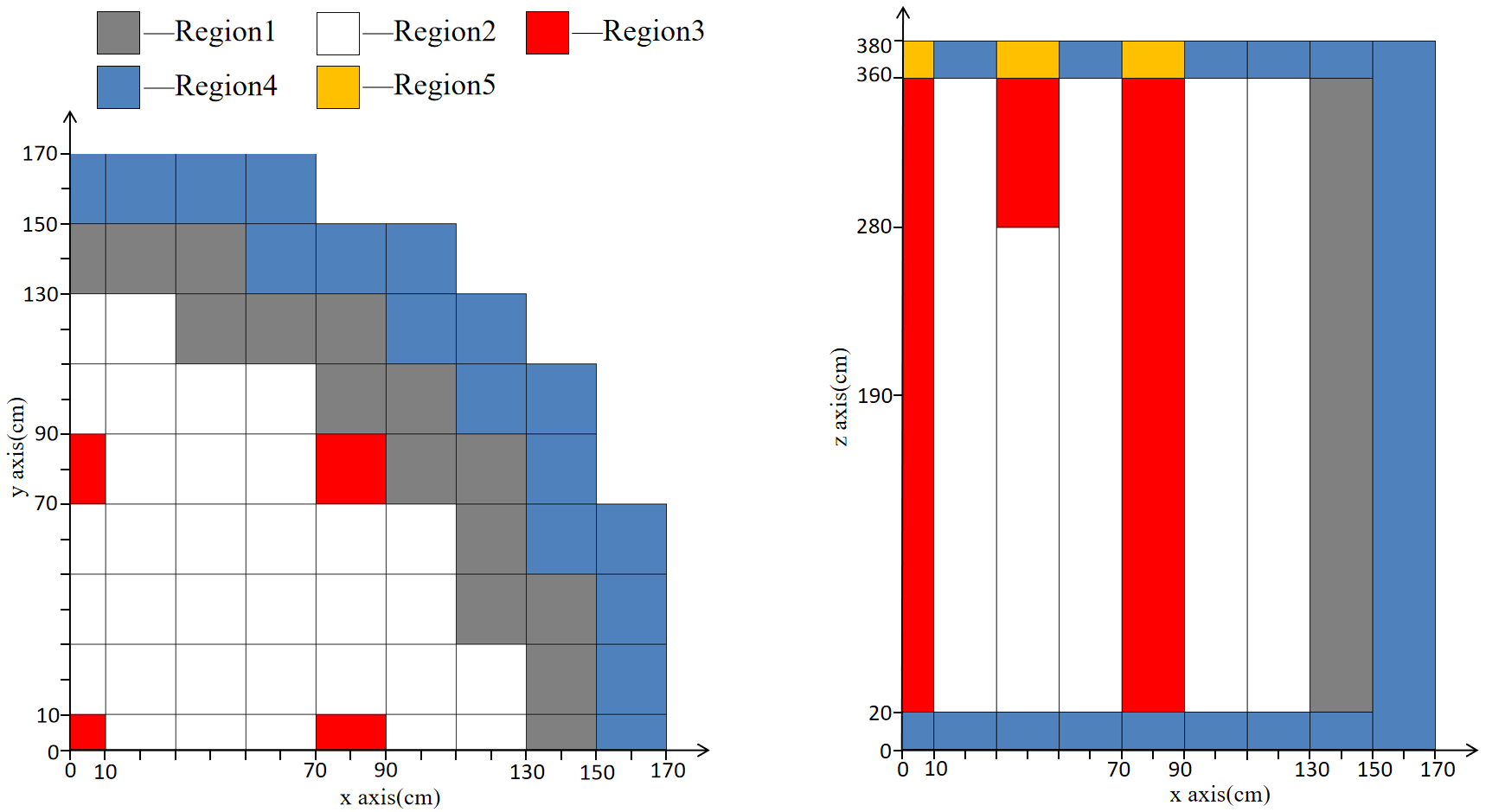}
		\caption{Three-dimensional IAEA problem calculation domain. Different colors represent different regions. In this case, there are five regions in the domain. Left: first quadrant of $x-y$ plane at $z$ = 190 $\text{cm}$. Right: $x-z$ plane at $y$ = 0 $\text{cm}$.}
		\label{fig-IAEA}
	\end{figure*}

	The IAEA benchmark problem is a classic benchmark case widely used for the development of reactor physics computational programs, as published by the Argonne National Laboratory. It refers to a three-dimensional neutron diffusion eigenvalue problem. The computational domain of the problem is illustrated in Figure \ref{fig-IAEA}, and it also includes a two-dimensional problem, which is a cross-section at $z=190cm$. The material parameters are presented in Table \ref{tab-IAEA-coef-2D}. The computational domain includes fuel, control rods, reflector, and their corresponding mixed regions, with a total of four regions corresponding to four sets of coefficients.
	We apply the Neumann loss function \eqref{loss-neumann} at the boundaries $x=0$ and $y=0$, and impose Robin boundary conditions \eqref{loss-robin} at the external boundaries where $c_{bou} = 0.4692$. The critical eigenvalue of this problem is $k_{\text{eff}}^{FF} = 1.0296$, which is solved by FreeFEM++ with intervals $\Delta x =1$ and $\Delta y =1$. 
	
\begin{figure*}[h]
		\begin{minipage}{0.32\textwidth}
			\centering
			\includegraphics[width=5.5cm]{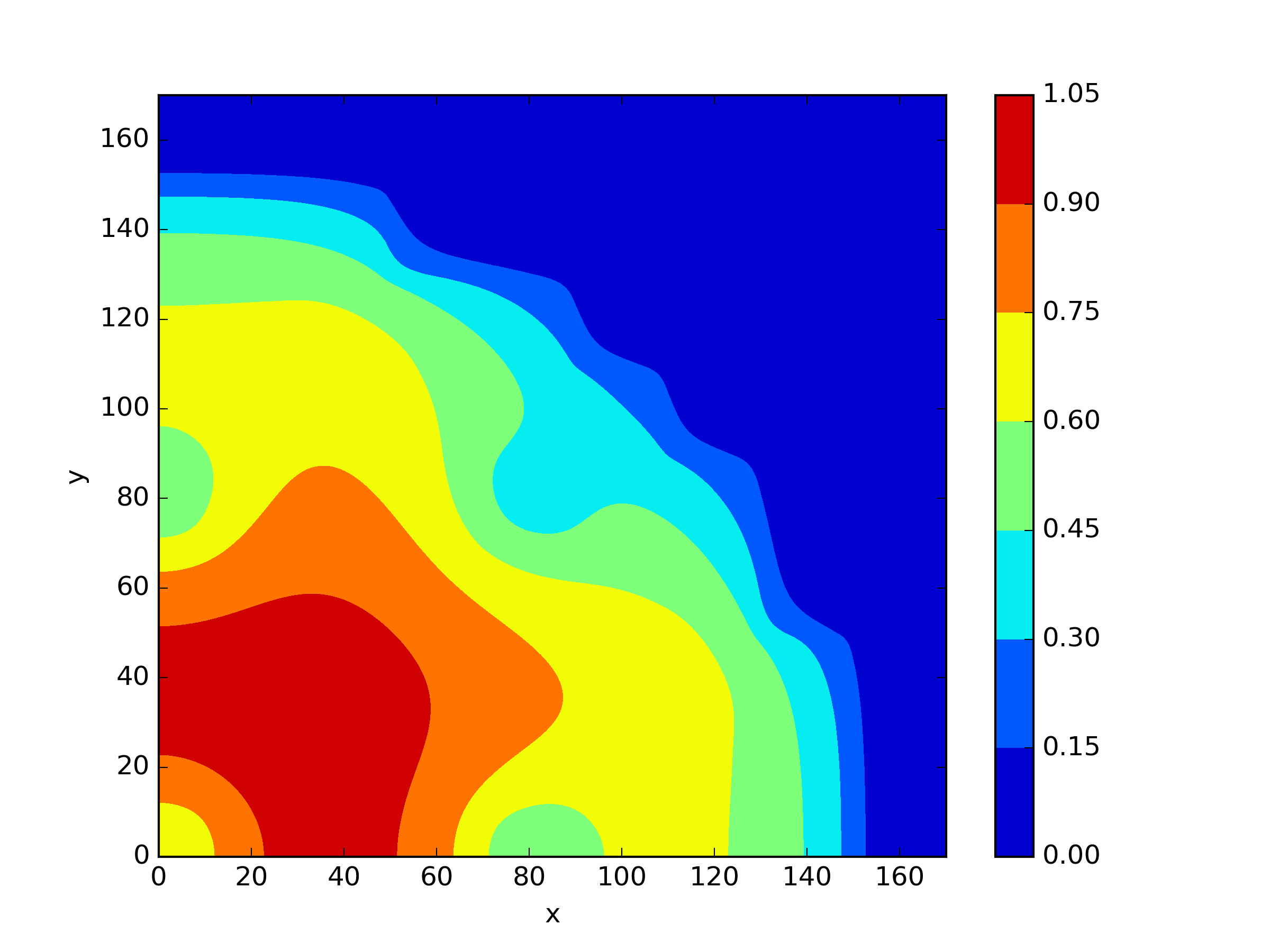}
		\end{minipage}
		\begin{minipage}{0.32\textwidth}
			\centering
			\includegraphics[width=5.5cm]{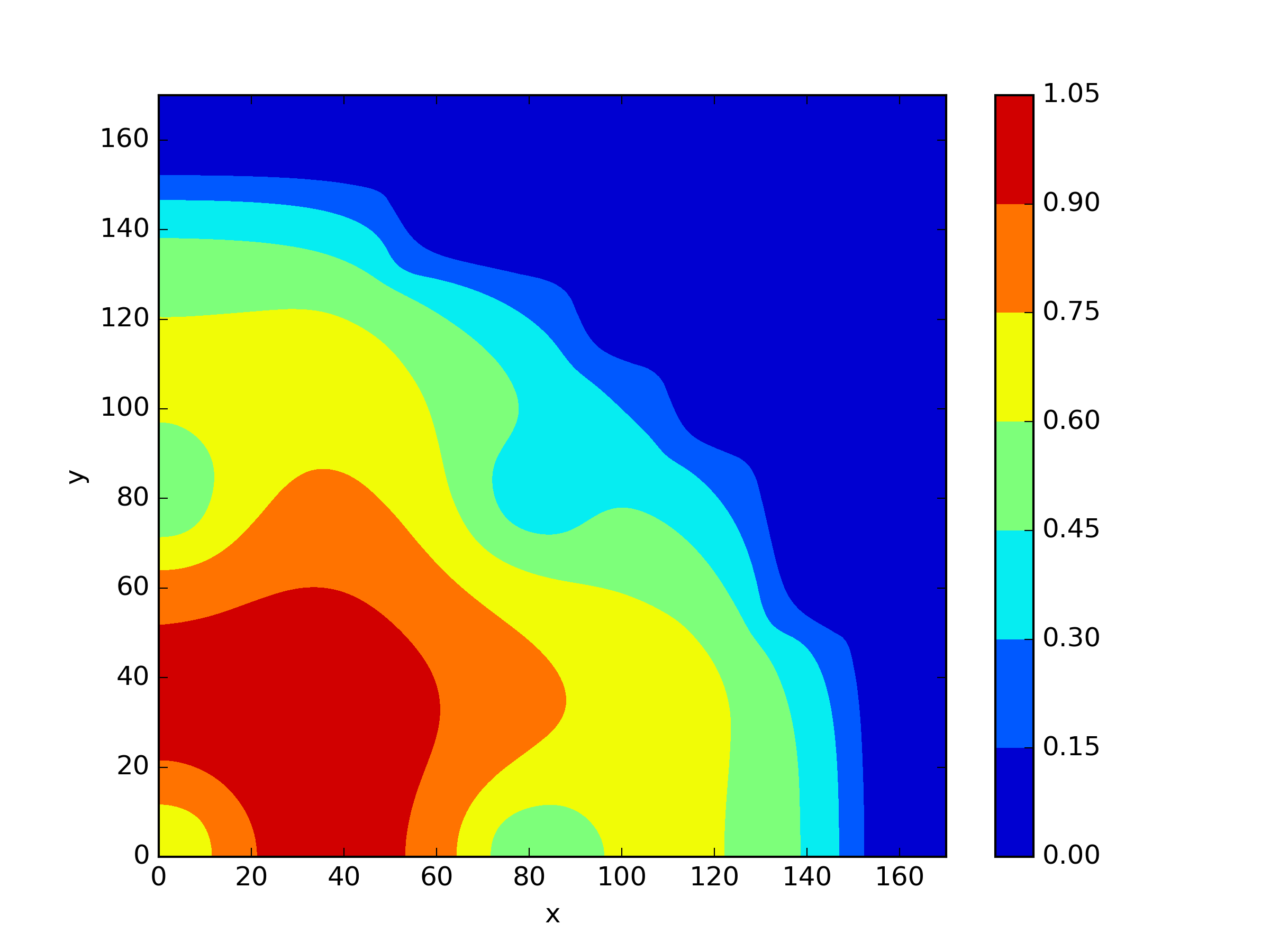}
		\end{minipage}
		\begin{minipage}{0.32\textwidth}
			\centering
			\includegraphics[width=5.5cm]{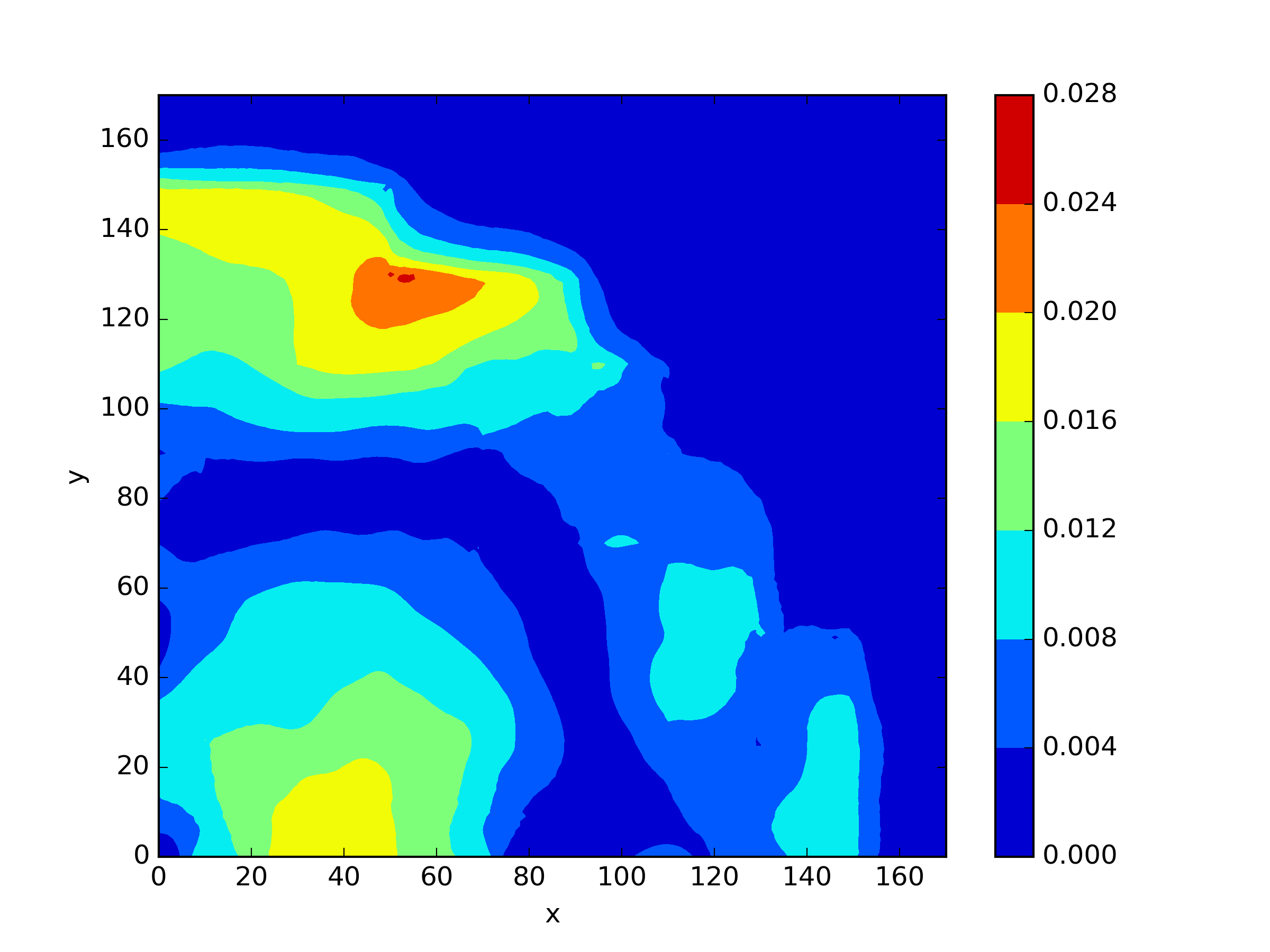}
		\end{minipage}
		
		\begin{minipage}{0.32\textwidth}
			\centering
			\includegraphics[width=5.5cm]{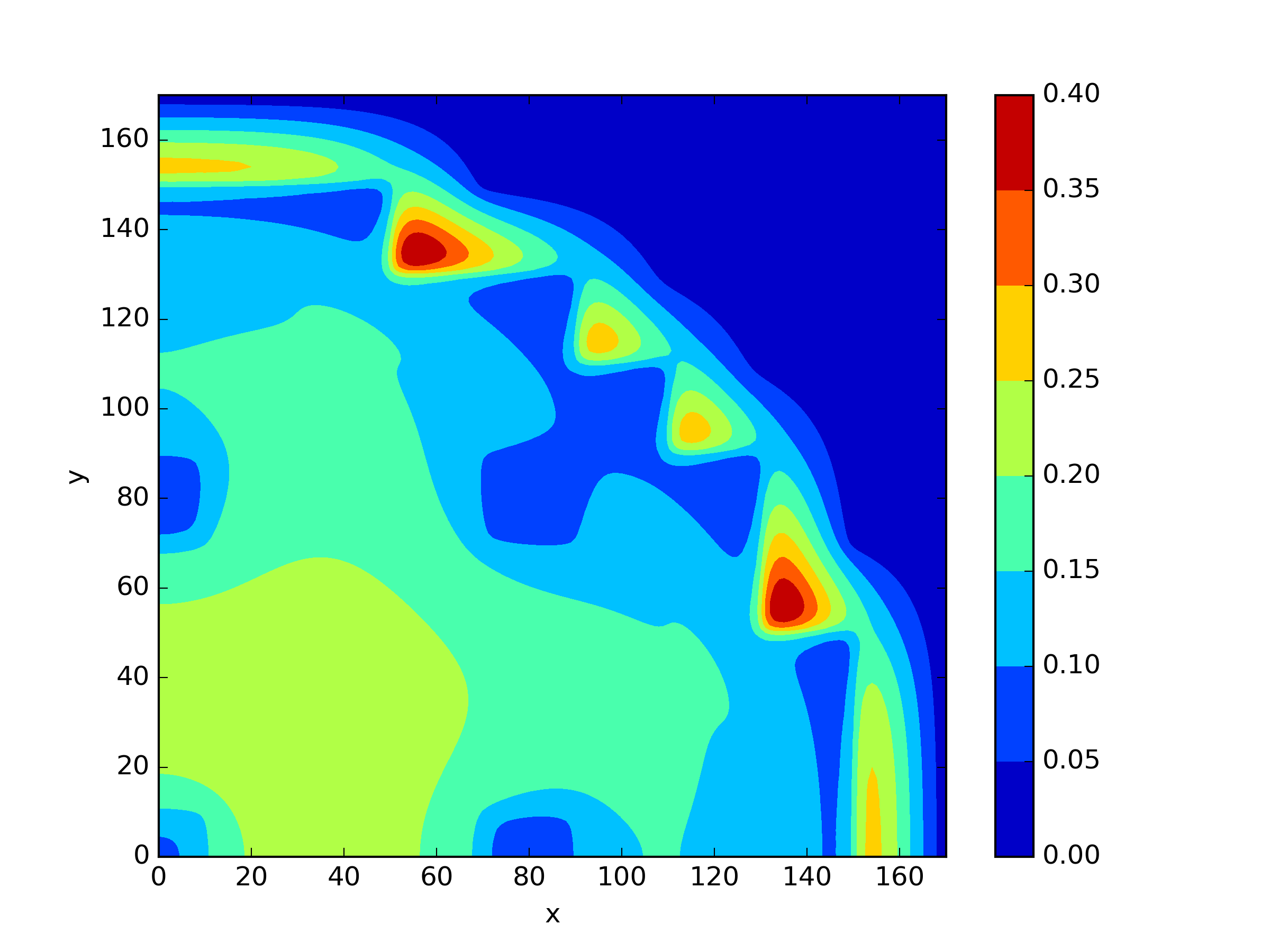}
		\end{minipage}
		\begin{minipage}{0.32\textwidth}
			\centering
			\includegraphics[width=5.5cm]{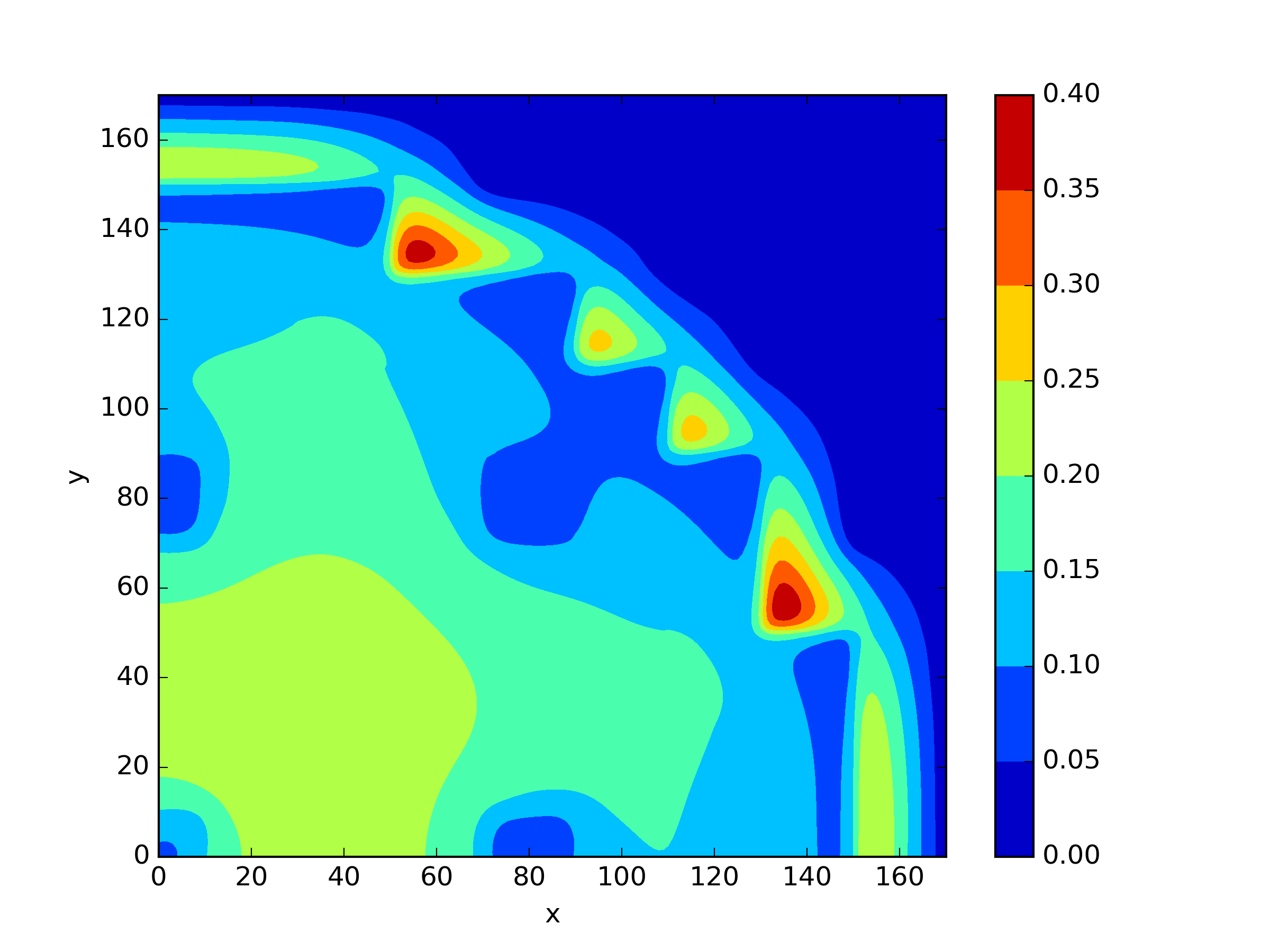}
		\end{minipage}
		\begin{minipage}{0.32\textwidth}
			\centering
			\includegraphics[width=5.5cm]{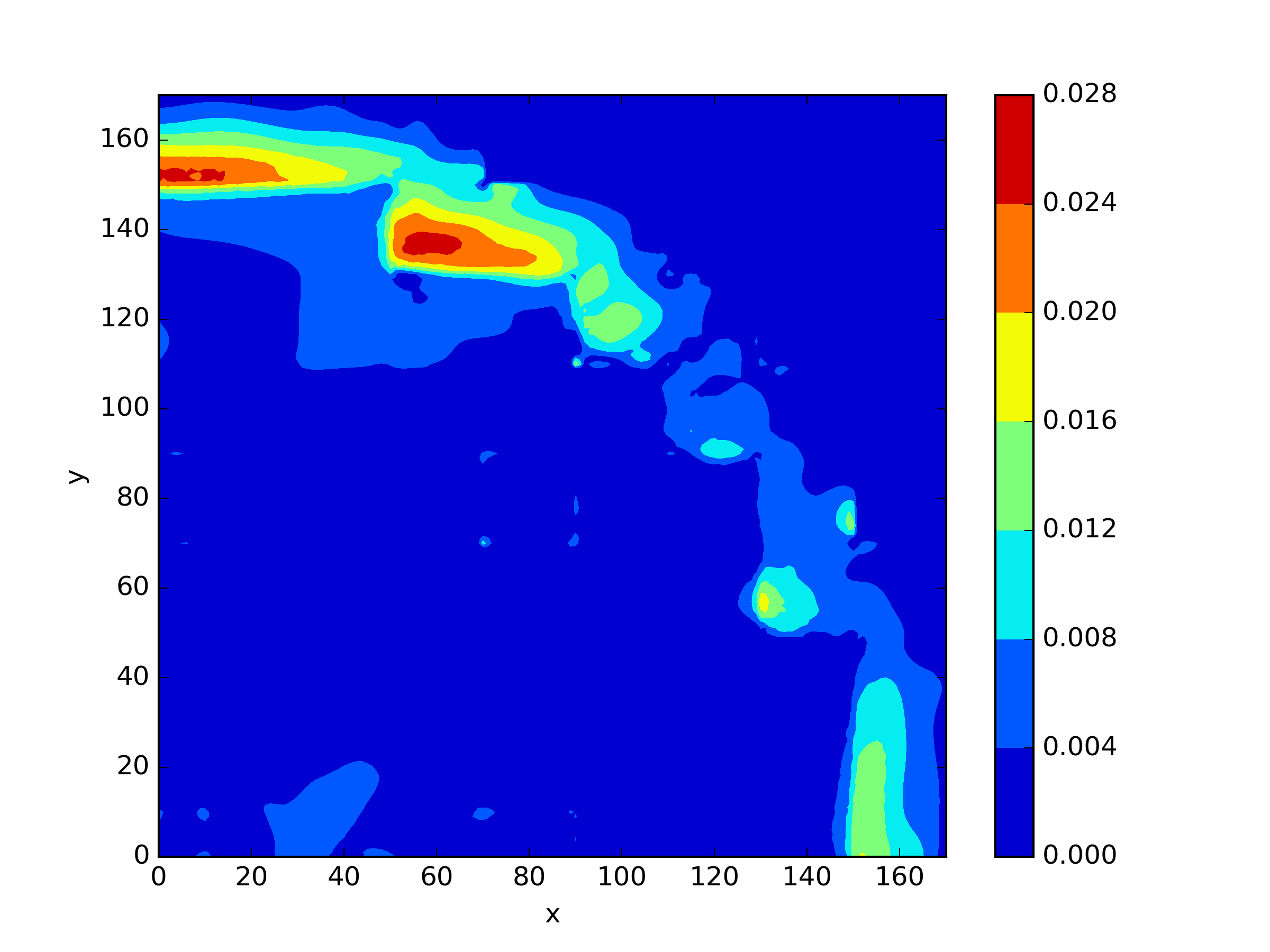}
		\end{minipage}
		\caption{Results of 2-D IAEA problem: Left: reference solution. Mid: neural network solution. Right: absolute error. First row: information of $\phi_1(x, y)$. Second row: information of $\phi_2(x, y)$.}
		\label{fig-IAEA-2D2R}
	\end{figure*}

	Indeed, the IAEA problem is larger in scale and involves more materials compared to the TWIGL problem, which increases the difficulty of solving it using neural networks. We chose 4 residual blocks, with 40 neurons per layer in each residual block. The Adam optimizer with a learning rate of $0.001$ was used for $100,000$ epochs optimization during the training process. {For the two-dimensional IAEA problem, we have selected 23,738 residual points and 703 interface points.} We applied the three residual loss function mentioned earlier to this problem as well (Table \ref{tab-IAEA-2D2G} and Figure \ref{fig-IAEA-2D2R}). In this scenario, the drawback of $Loss_{IPM}$ is amplified, as evident from the table. The errors of $Loss_{IPM}$ are significantly larger because $Loss_{IPM}$ discards the information of $\phi_1$ during the iteration process, resulting in a decrease in the solution accuracy, while $Loss_{De}$ and $Loss_{DI}$ continue to maintain a good level of accuracy. Furthermore, the results with shift $\sigma = 1$ show slightly better performance.
	
	\begin{table*}[h]
		\caption{Coefficients of different regions for IAEA problem}
		\centering
		\begin{tabular}{|c|c|c|c|c|c|c|c|c|}
			\hline
			& $D_1$ & $D_2$ & $\Sigma_{a,1}$ & $\Sigma_{a,2}$ & $\Sigma_{1\rightarrow2}$ & $\nu\Sigma_{f,1}$ & $\nu\Sigma_{f,2}$\\
			& $(cm)$ & $(cm)$ & $(cm^{-1})$ & $(cm^{-1})$ & $(cm^{-1})$ & $(cm^{-1})$ & $(cm^{-1})$\\
			\hline
			Region-1 & 1.5 & 0.4 & 0.01 & 0.085 & 0.02 & 0.0 & 0.135\\
			\hline
			Region-2 & 1.5 & 0.4 & 0.01 & 0.08 & 0.02 & 0.0 & 0.135\\
			\hline
			Region-3 & 1.5 & 0.4 & 0.01 & 0.13 & 0.02 & 0.0 & 0.135\\
			\hline
			Region-4 & 2.0 & 0.3 & 0.0 & 0.01 & 0.04 & 0.0 & 0.0\\
			\hline
			Region-5 & 2.0 & 0.3 & 0.0 & 0.055 & 0.04 & 0.0 & 0.0\\
			\hline
		\end{tabular}
		\label{tab-IAEA-coef-2D}
	\end{table*}
	
	\begin{table*}[h]
		\caption{Recording numerical results for 2-D IAEA problem with different residual loss function applied.}
		\centering
		\begin{tabular}{ccccccc}
			\hline
			& Shift & $\mathbf{E}_R(k_{\text{eff}})$ & $\mathbf{E}_{R,\infty}(\phi_1)$ & $\mathbf{E}_{R,2}(\phi_1)$ & $\mathbf{E}_{R,\infty}(\phi_2)$ & $\mathbf{E}_{R,2}(\phi_2)$ \\
			\hline
			$Loss_{IPM}$ & 0 & \pmb{3.5993e-04} & 5.8825e-02 & 3.4499e-02 & 1.5972e-01 & 7.1278e-02 \\
			\hline
			$Loss_{De}$ & 0 & 7.8179e-04 & 2.9658e-02 & 2.2165e-02 & 8.3564e-02 & 5.0352e-02 \\
			\hline
			$Loss_{DI}$ & 1 & 5.1034e-04 & 4.0134e-02 & 1.7337e-02 & 7.9622e-02 & 4.0727e-02 \\
			\hline
			$Loss_{De}$ & 1 & 6.9851e-04 & \pmb{2.6378e-02} & \pmb{1.5618e-02} & \pmb{6.5194e-02} & \pmb{3.5830e-02} \\
			\hline
		\end{tabular}
		\label{tab-IAEA-2D2G}
	\end{table*}
	
\begin{figure*}[h]
		\centering
		\begin{minipage}[b]{0.32\linewidth}
			\subfloat[Relative error of flux]{
				\includegraphics[width=5.15cm]{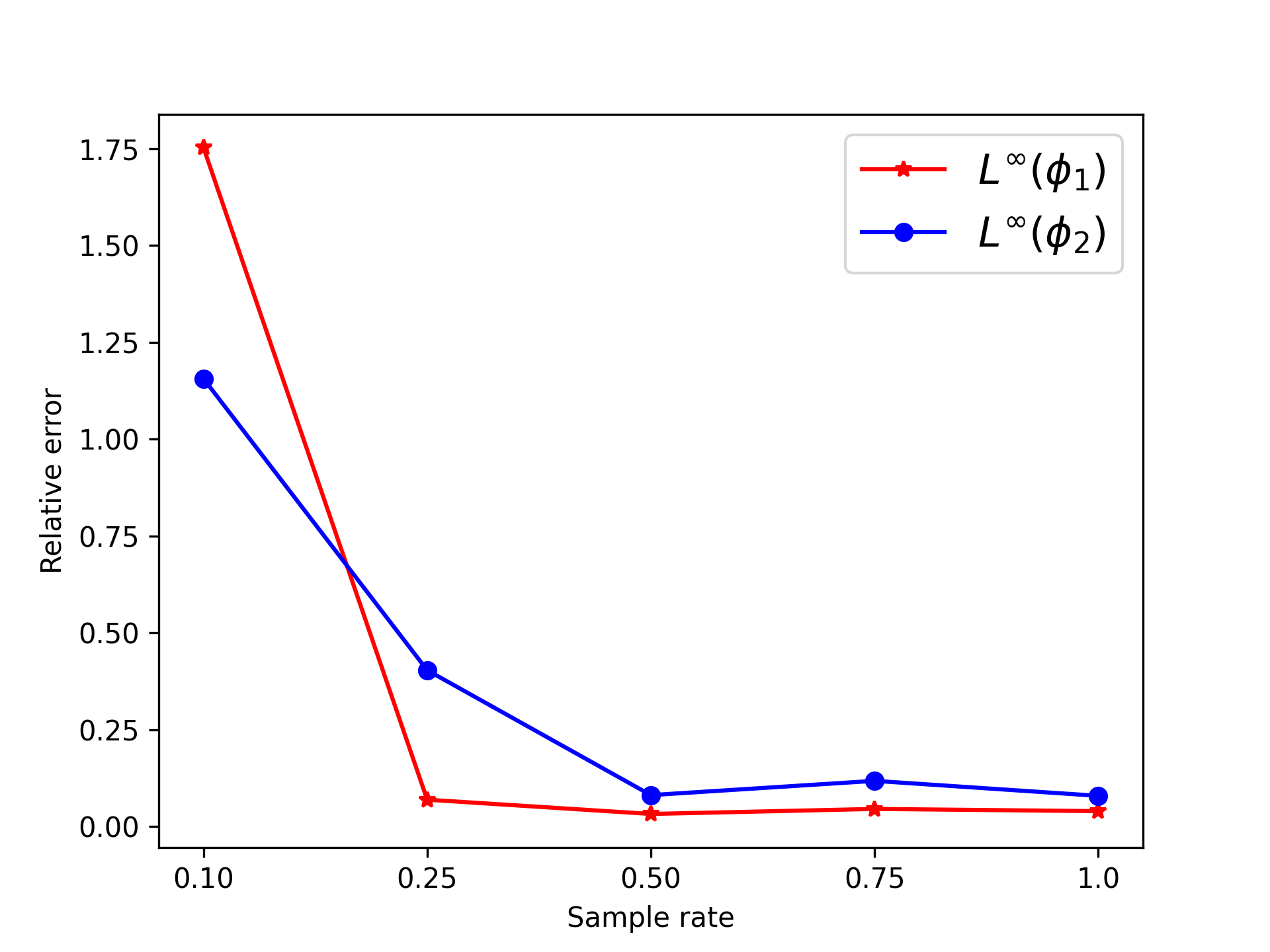}}
		\end{minipage}
		\begin{minipage}[b]{0.32\linewidth}
			\subfloat[Relative error of flux]{
				\includegraphics[width=5.15cm]{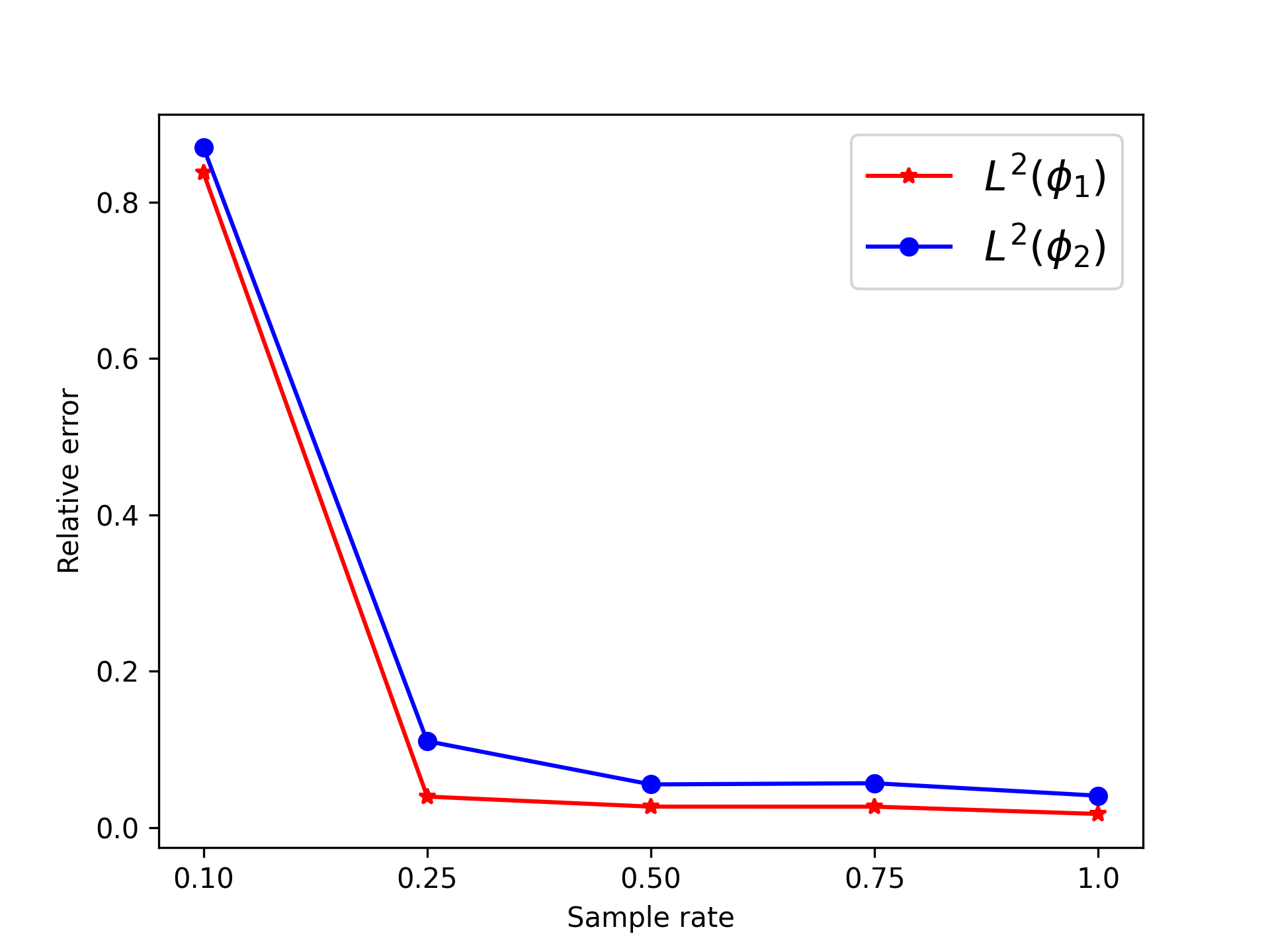}}
		\end{minipage}
		\begin{minipage}[b]{0.32\linewidth}
			\subfloat[Relative error of eigenvalue]{
				\includegraphics[width=5.15cm]{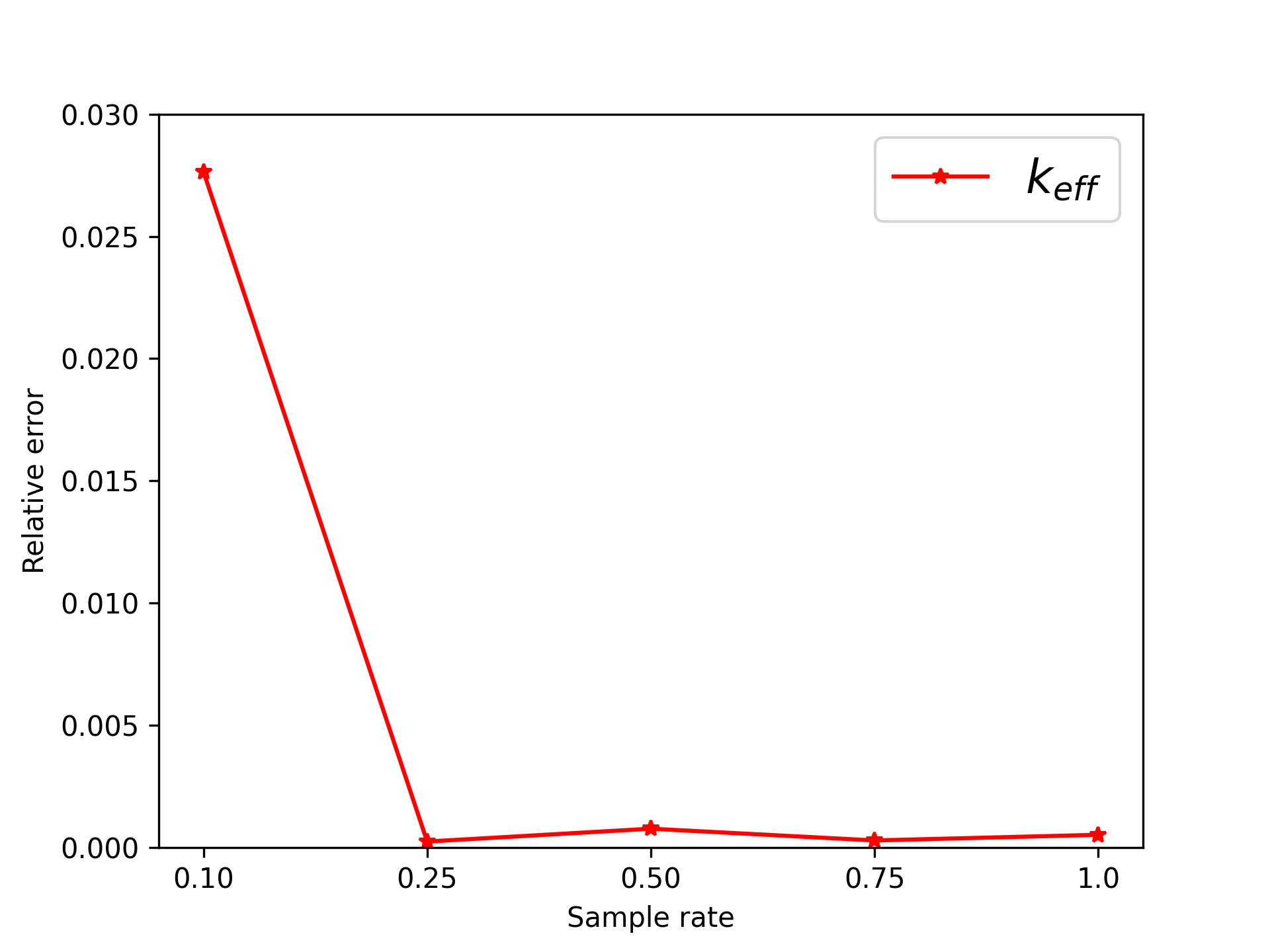}}
		\end{minipage}
		\caption{2-D IAEA problem: relative error of neural network solution at different sample rate.}
		\label{pic-IAEA-2D2G-samplerate}
	\end{figure*}

\begin{figure*}[h]
		\centering
		\begin{minipage}[b]{0.49\linewidth}
			\subfloat[Relative $L^{\infty}$ error of $\phi_1$]{
				\includegraphics[width=7cm]{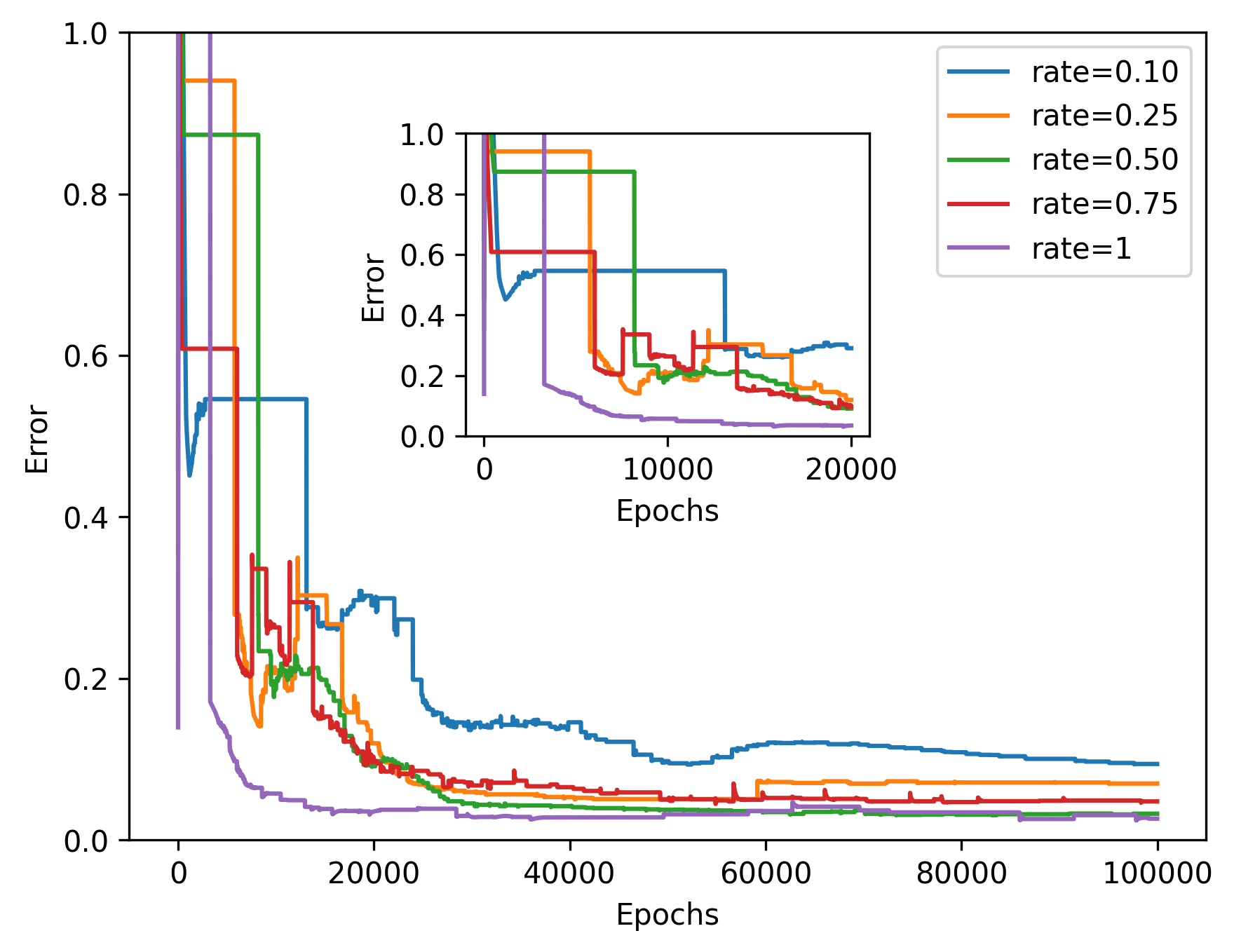}}
		\end{minipage}
		\begin{minipage}[b]{0.49\linewidth}
			\subfloat[Relative $L^{2}$ error of $\phi_1$]{
				\includegraphics[width=7cm]{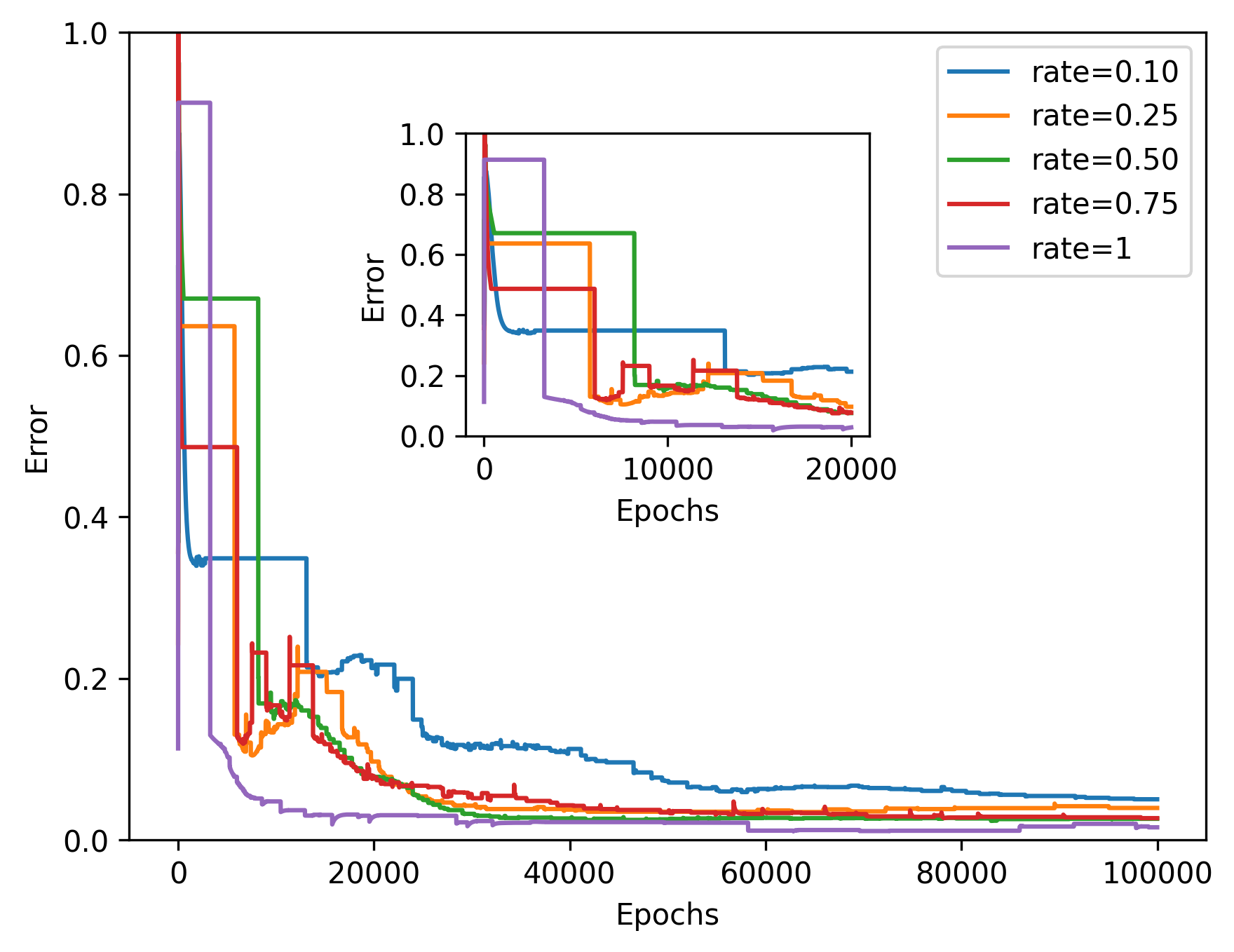}}
		\end{minipage}
		\begin{minipage}[b]{0.49\linewidth}
			\subfloat[Relative $L^{\infty}$ error of $\phi_2$]{
				\includegraphics[width=7cm]{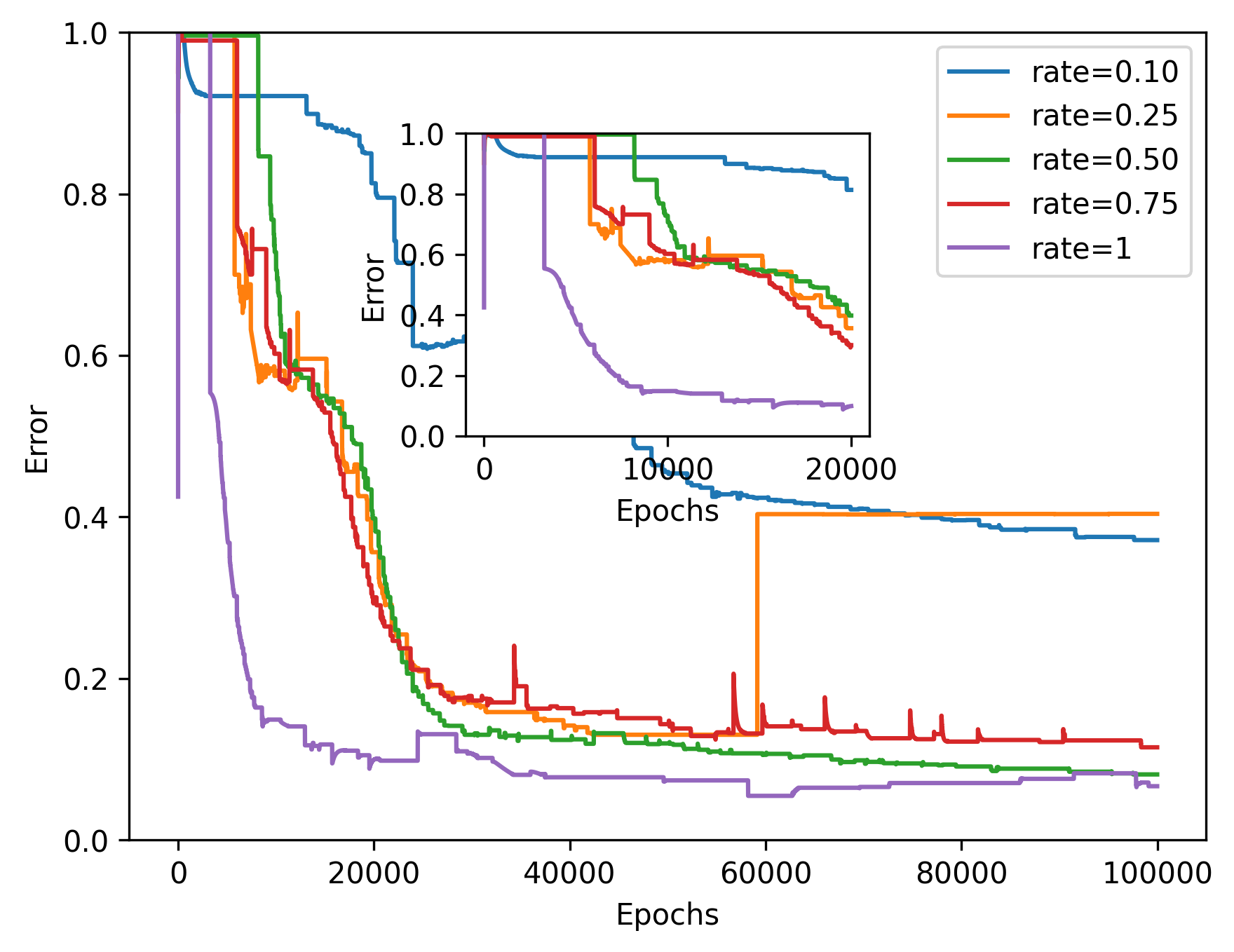}}
		\end{minipage}
		\begin{minipage}[b]{0.49\linewidth}
			\subfloat[Relative $L^{2}$ error of $\phi_2$]{
				\includegraphics[width=7cm]{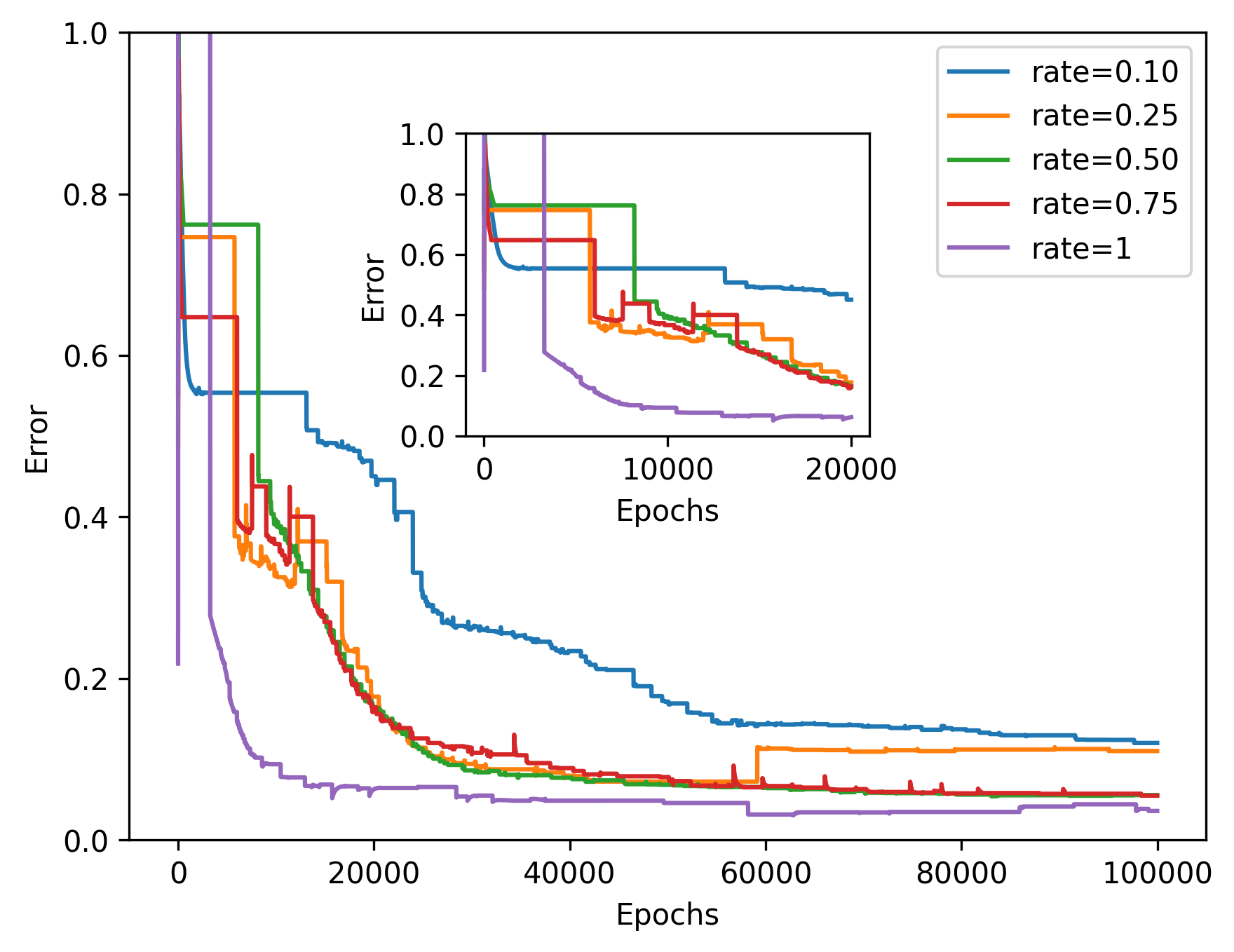}}
		\end{minipage}
		\caption{2-D IAEA: the first row shows relative error of $\phi_1$ during training process; the second row shows relative error of $\phi_2$ during training process. Black, red, green, blue, and purple colors represent the relative error curves at different sampling rates of 0.1, 0.25, 0.5, 0.75, and 1, respectively.}
		\label{pic-IAEA-2D2G-train}
	\end{figure*}

	Similarly, we investigated the impact of the number of sampling points on the solution accuracy for the two-dimensional problem of IAEA. We conducted tests using five different levels of training point numbers. The experimental results (Figure \ref{pic-IAEA-2D2G-samplerate} and Figure \ref{pic-IAEA-2D2G-train}), in the context of this more complex problem, clearly demonstrate the importance of the number of sampling points for problem solving. Having a greater number of sampling points yields higher accuracy in the results.
	
	To investigate the effectiveness of applying interface conditions to the IAEA problem, we conducted an experiment related to interface conditions. Based on the experimental results (Figure \ref{pic-interface} and Figure \ref{pic-int-keff}), it is evident that when interface conditions are not applied at all, the problem remains unsolvable. Moreover, for more complex problems, the implementation of interface conditions becomes even more crucial.
	
	\subsubsection{3-Region TWIGL}
\begin{figure*}
		\centering
		\includegraphics[width=10cm]{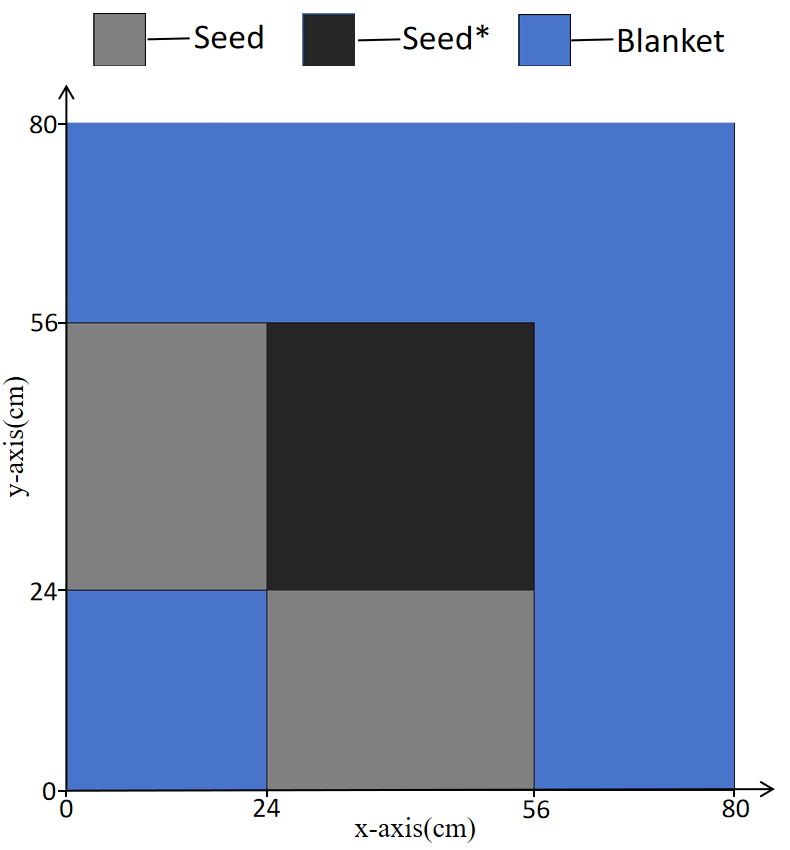}
		\caption{Computational domain of TWIGL 3-R problem. Seed* is a new material (Table \ref{tab-coef-2D-TWIGL-3R}).}
		\label{fig-TWIGL-3R}
	\end{figure*}

\begin{figure*}[h]
		\centering
		\begin{minipage}[b]{0.49\linewidth}
			\subfloat[Relative $L^{\infty}$ error of $\phi_1$]{
				\includegraphics[width=8cm]{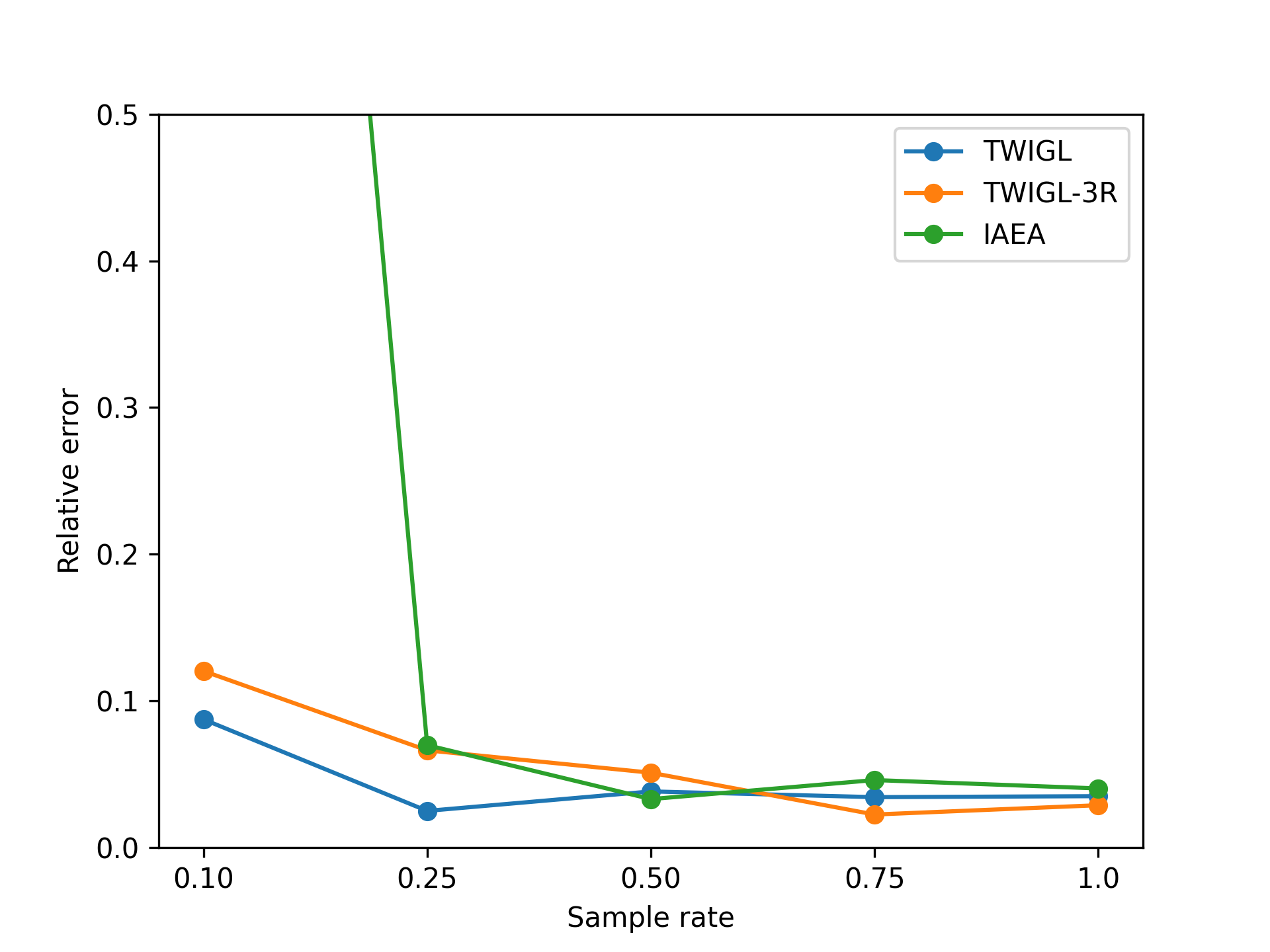}}
		\end{minipage}
		\begin{minipage}[b]{0.49\linewidth}
			\subfloat[Relative $L^{2}$ error of $\phi_1$]{
				\includegraphics[width=8cm]{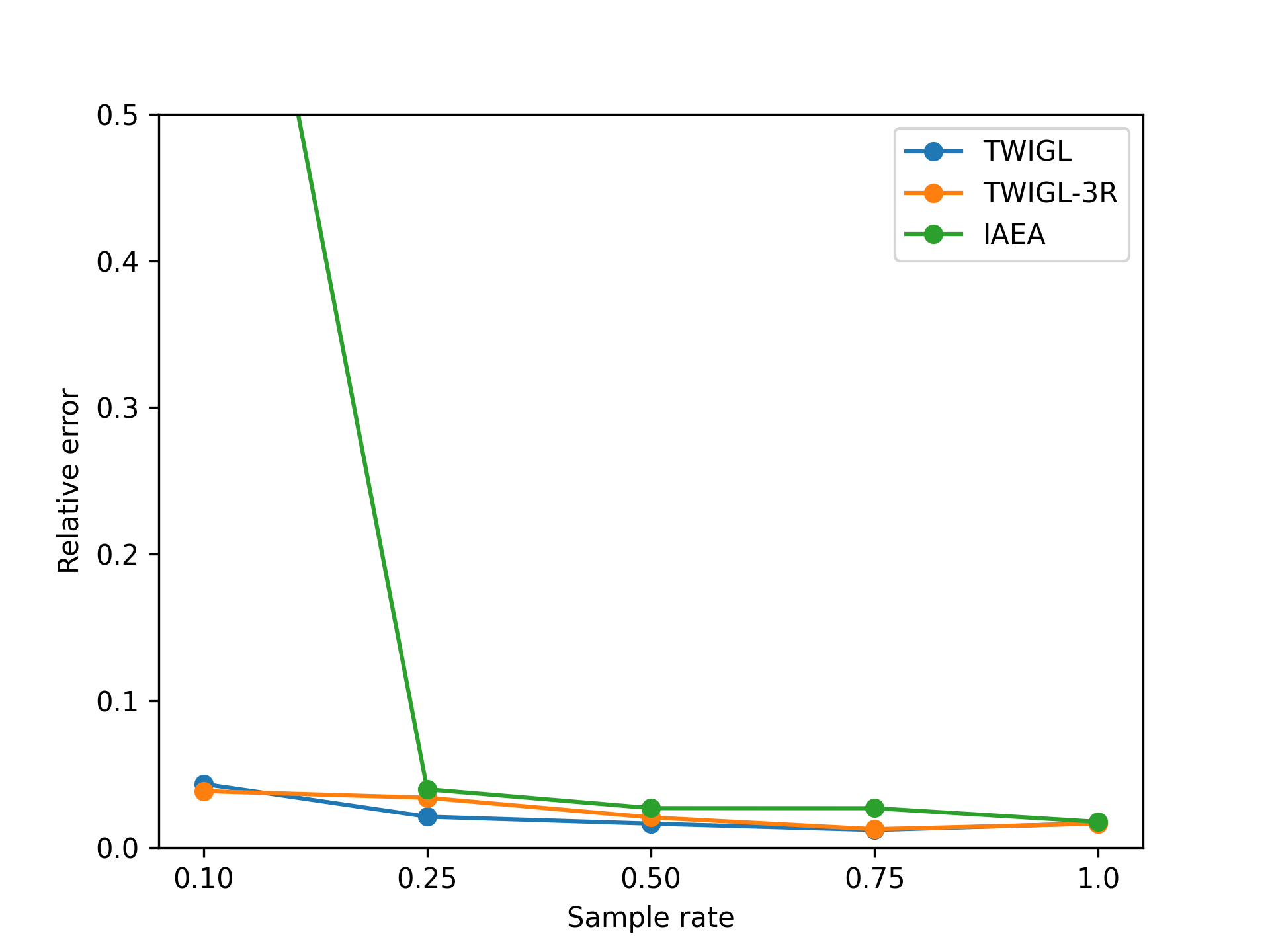}}
		\end{minipage}
		\begin{minipage}[b]{0.49\linewidth}
			\subfloat[Relative $L^{\infty}$ error of $\phi_2$]{
				\includegraphics[width=8cm]{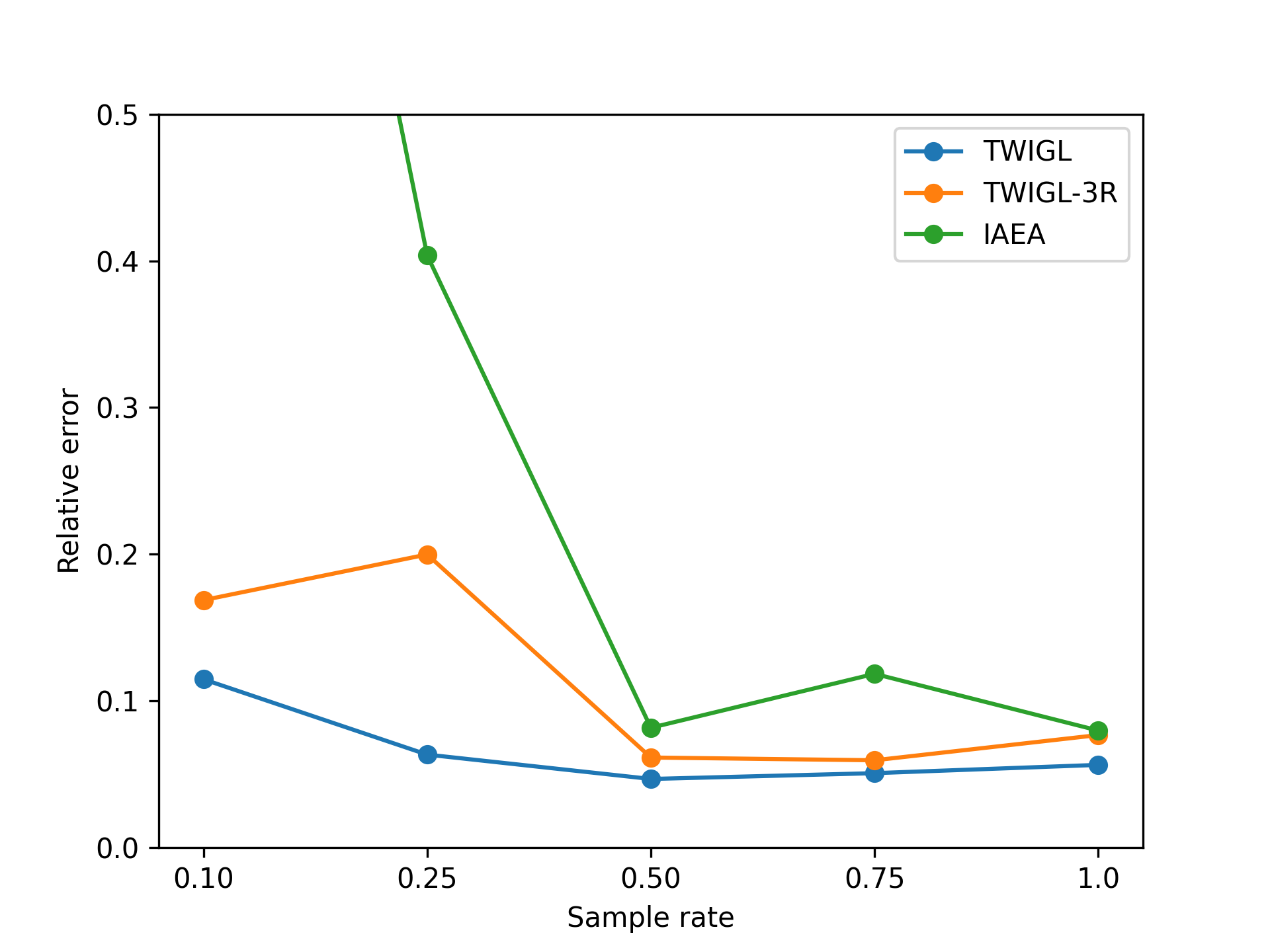}}
		\end{minipage}
		\begin{minipage}[b]{0.49\linewidth}
			\subfloat[Relative $L^{2}$ error of $\phi_2$]{
				\includegraphics[width=8cm]{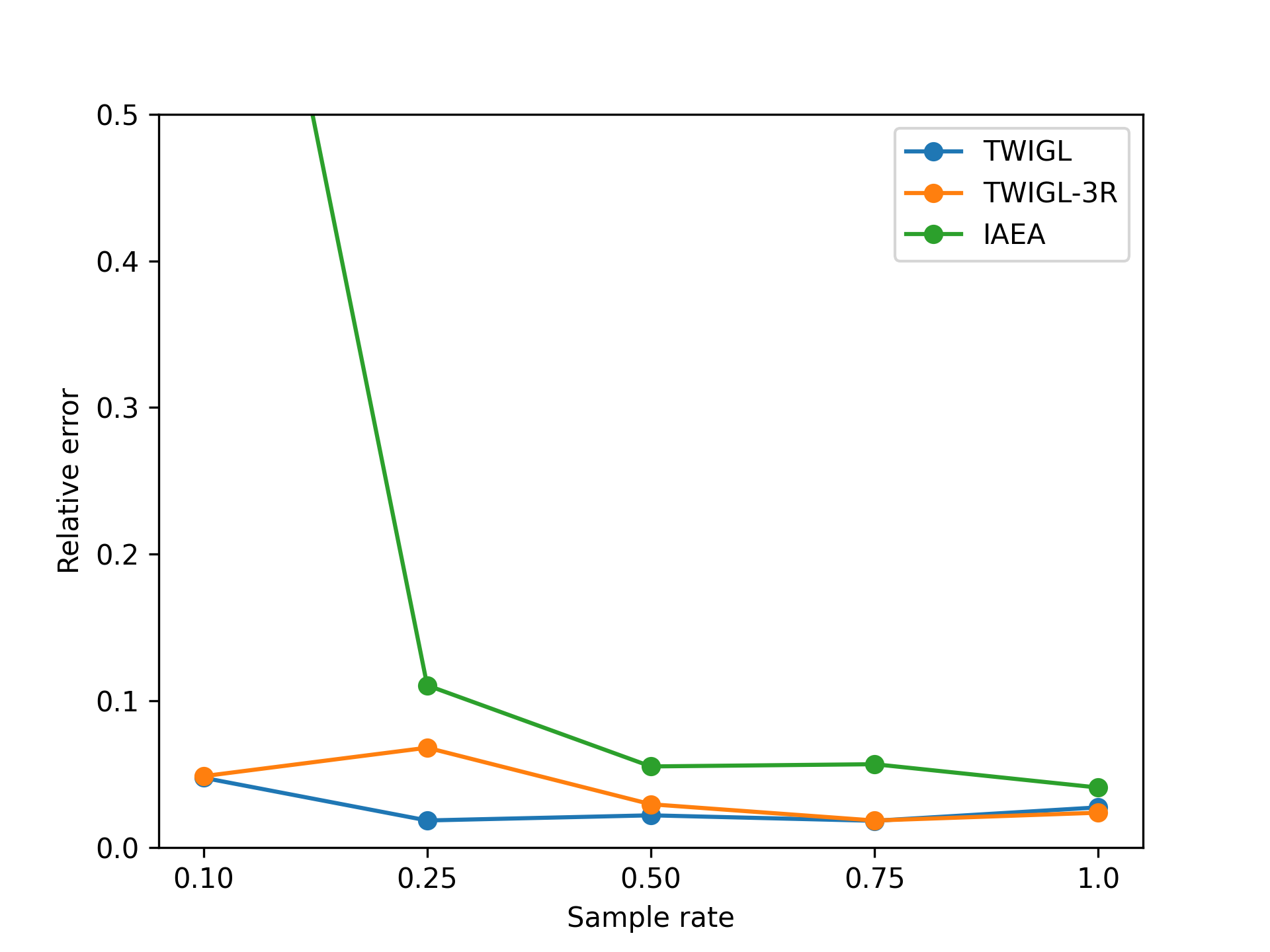}}
		\end{minipage}
		\caption{Relative error curves of different cases under different sampling rates. Blue represents the 2-D TWIGL problem, orange represents the 2-D 3-region TWIGL problem, and green represents the 2-D IAEA problem. The first row shows relative error of $\phi_1$ during training process; the second row shows relative error of $\phi_2$ during training process.}
		\label{pic-3cases}
	\end{figure*}

	In this section, we also considered a slightly more complex version of the TWIGL problem. We extended the computational domain from two materials to three materials, resulting in more intricate interface conditions (Figure \ref{fig-TWIGL-3R}). The boundary conditions for this case are the same as the previous TWIGL problem. Subsequently, we performed neural network computations for this problem at different sampling rates, comparing them with the performance of the best-performing algorithm for TWIGL, as well as two other two-dimensional problems. {The set of training points is consistent with the TWIGL problem.} The results are illustrated in Figure \ref{pic-3cases}.
	
	\begin{table*}[htb]
		\caption{Coefficient of different regions for TWIGL 3-R problem}
		\centering
		\begin{tabular}{|c|c|c|c|c|c|c|c|c|}
			\hline
			Region & $D_1$ & $D_2$ & $\Sigma_{a,1}$ & $\Sigma_{a,2}$ & $\Sigma_{1\rightarrow2}$ & $\nu\Sigma_{f,1}$ & $\nu\Sigma_{f,2}$\\
			& $(cm)$ & $(cm)$ & $(cm^{-1})$ & $(cm^{-1})$ & $(cm^{-1})$ & $(cm^{-1})$ & $(cm^{-1})$\\
			\hline
			Seed* & 1.35 & 0.45 & 0.009 & 0.1 & 0.01 & 0.005 & 0.13\\
			\hline
		\end{tabular}
		\label{tab-coef-2D-TWIGL-3R}
	\end{table*}
	
	From the Figure \ref{pic-3cases}, it is evident that as we progress from TWIGL to 3-R TWIGL and then to the IAEA problem, the average error increases. This implies that the size and complexity of the problem significantly affect the accuracy of our neural network solution. However, it is worth noting that regardless of the case, there are always instances where the results satisfy a threshold of less than 0.08 or even 0.05, which are the engineering acceptance criteria.
	
	\subsection{3-D Problem}
	In the case of the 3-D TWIGL and IAEA problems, we initially employed FreeFEM++ as our reference solution. For discretization, we utilized a grid with intervals of 1$cm$ in the x and y directions, and intervals of 10$cm$ in the z direction. For these two three-dimensional problems, the boundary conditions from the two-dimensional problems naturally extend to the three-dimensional problems. The top and bottom surfaces of the three-dimensional computational domain are considered as the boundaries corresponding to the external interfaces. Quadratic polynomial elements were employed for the numerical approximation. Critical eigenvalue of three-dimensional TWIGL problem is $k_{\text{eff}}^{FF} = 0.8787$ and three-dimensional IAEA problem is $k_{\text{eff}}^{FF} = 1.0293$. 
	
	Our neural network was trained using a sampling approach where we considered all integer coordinate points within the computational domain. For the 3-D TWIGL problem, we selected 4 residual blocks, each containing 80 neurons in a single layer. As for the 3-D IAEA problem, we opted for 6 residual blocks, each with 60 neurons in a single layer. For both problems, The optimization algorithm remained Adam with a learning rate of $0.001$ and we have selected \eqref{loss-res-2} as residual loss function. {For three-dimensional problems, the TWIGL problem comprises 9,027,249 residual training points and 21,384 interface points; the IAEA problem includes 93,903,717 residual points and 284,772 interface points.}
	
	\begin{table*}[ht]
		\caption{Relative error of 3-D problem solved by neural network.}
		\centering
		\begin{tabular}{cccccc}
			\hline
			& $\mathbf{E}_R(k_{\text{eff}})$ & $\mathbf{E}_{R,2}(\phi_1)$ &  $\mathbf{E}_{R,2}(\phi_2)$\\
			\hline
			TWIGL-3D & 6.2617e-03 & 5.5424e-02 & 8.1373e-02 \\
			IAEA-3D & 3.4704e-03 & 1.2472e-01 & 2.6125e-01 \\
			\hline
		\end{tabular}
		\label{tab-3D}
	\end{table*}
	
	\begin{figure*}[h]
		\centering
		\begin{minipage}{0.3\textwidth}
			\centering
			\includegraphics[width=4.7cm]{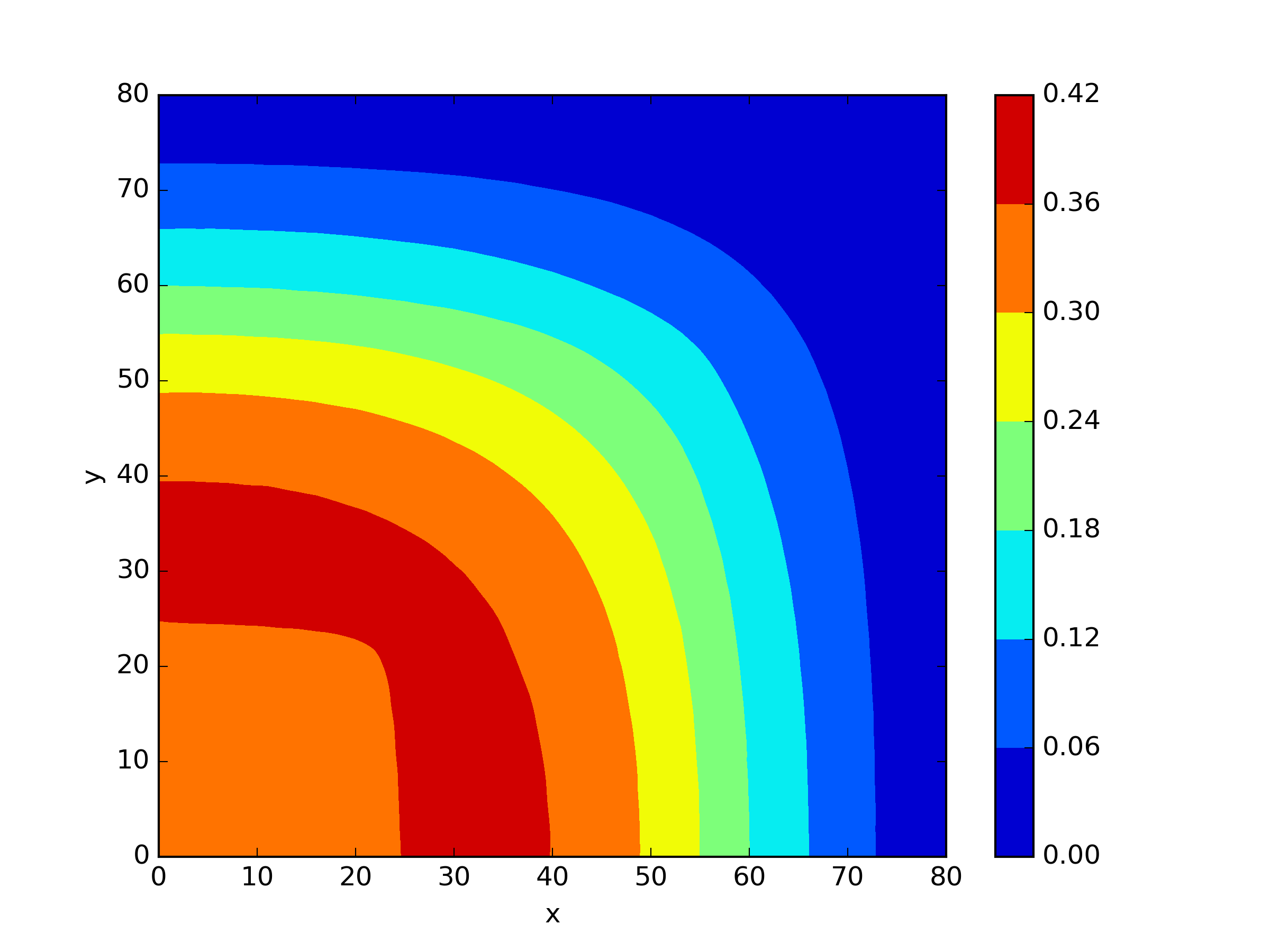}
		\end{minipage}
		\begin{minipage}{0.3\textwidth}
			\centering
			\includegraphics[width=4.7cm]{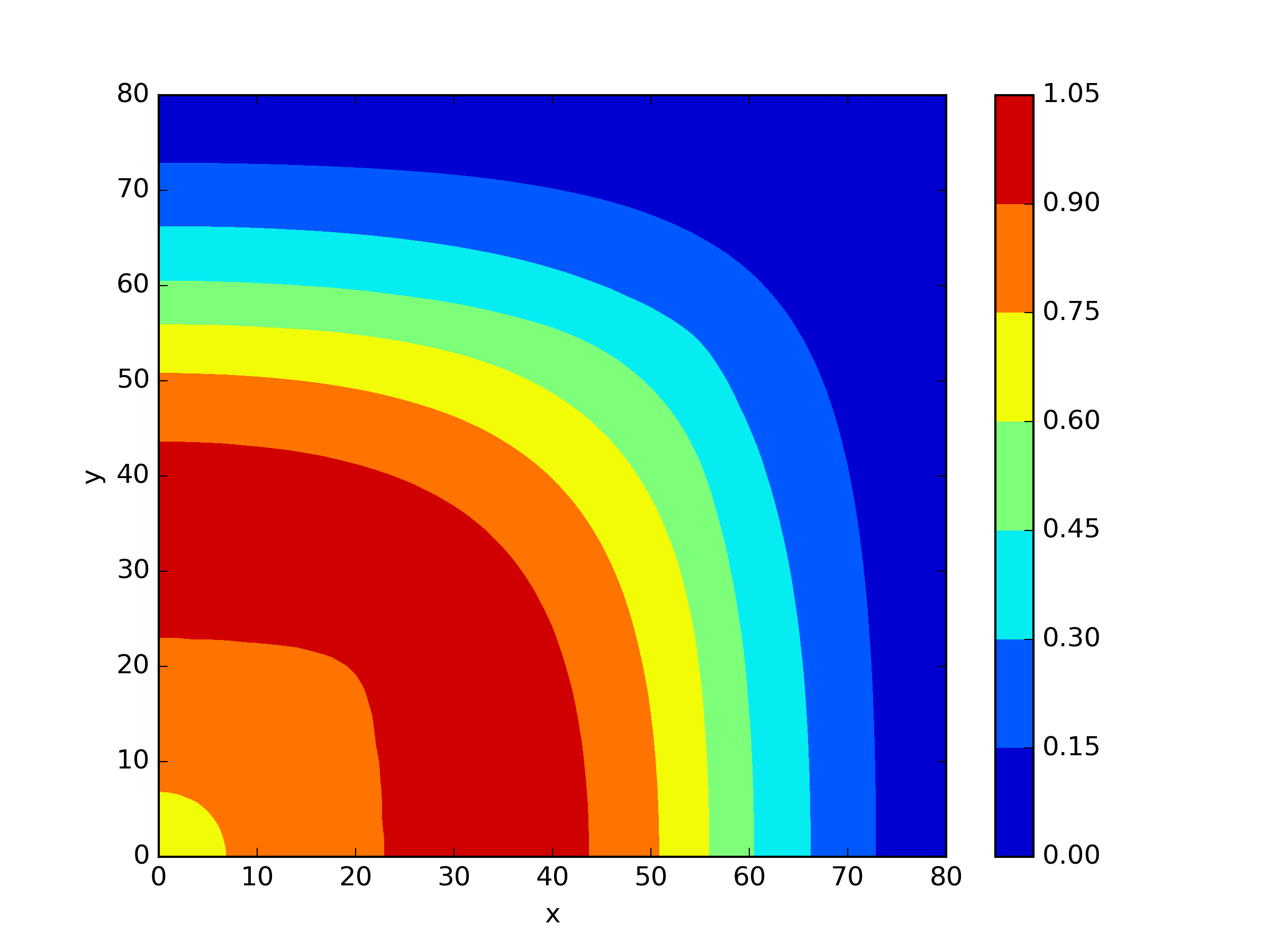}
		\end{minipage}
		\begin{minipage}{0.3\textwidth}
			\centering
			\includegraphics[width=4.7cm]{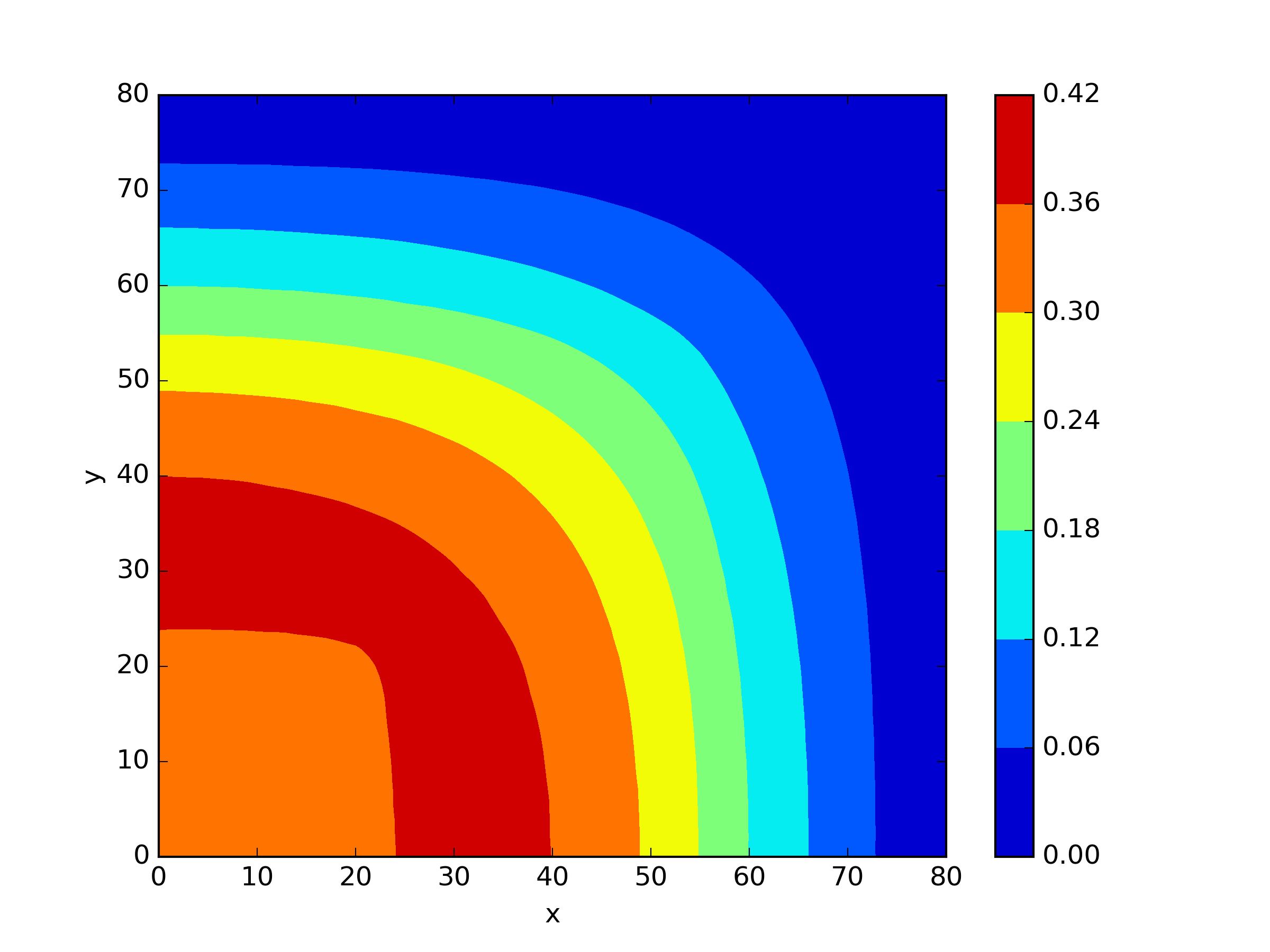}
		\end{minipage}
		
		\begin{minipage}{0.3\textwidth}
			\centering
			\includegraphics[width=4.7cm]{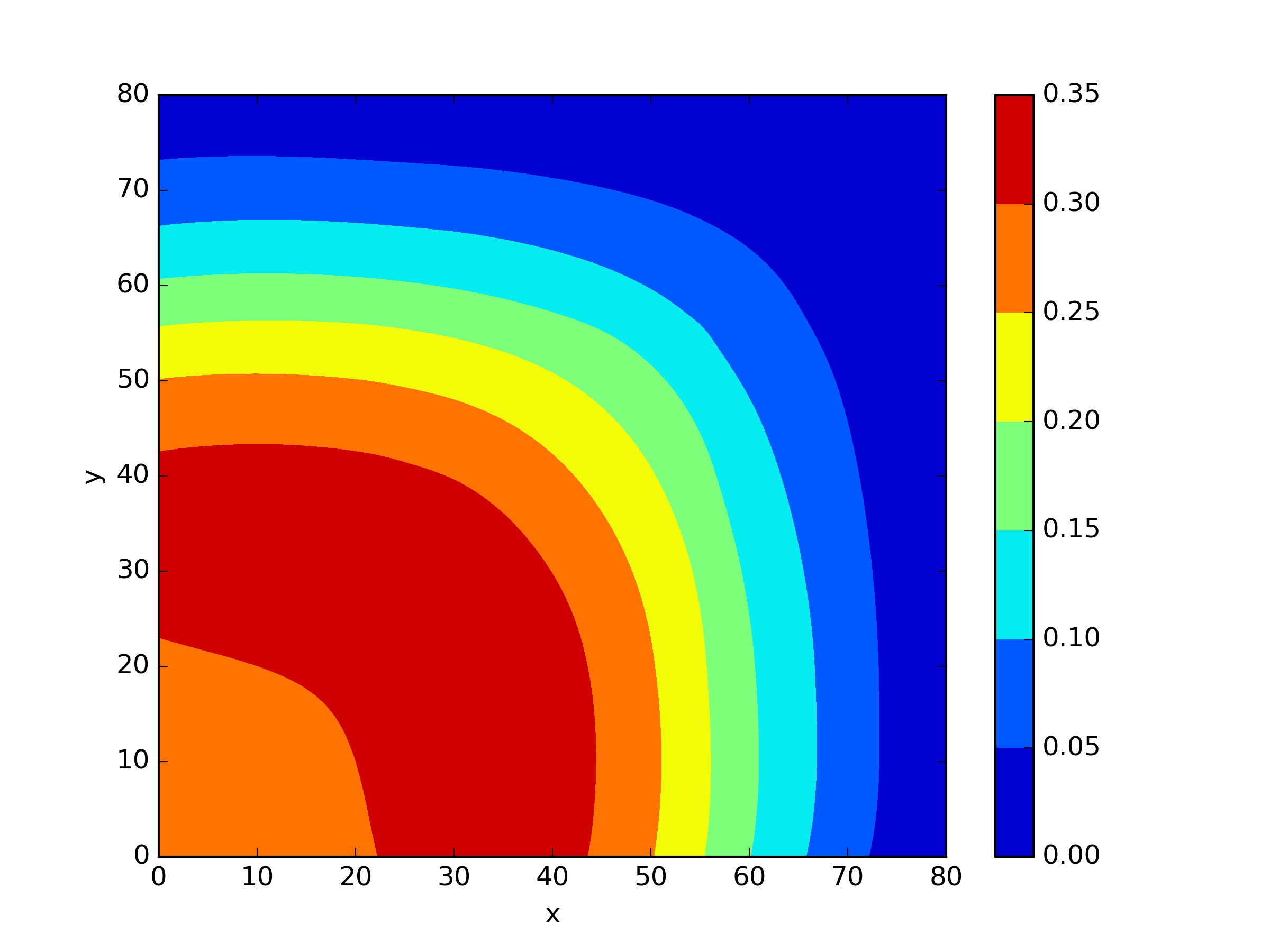}
		\end{minipage}
		\begin{minipage}{0.3\textwidth}
			\centering
			\includegraphics[width=4.7cm]{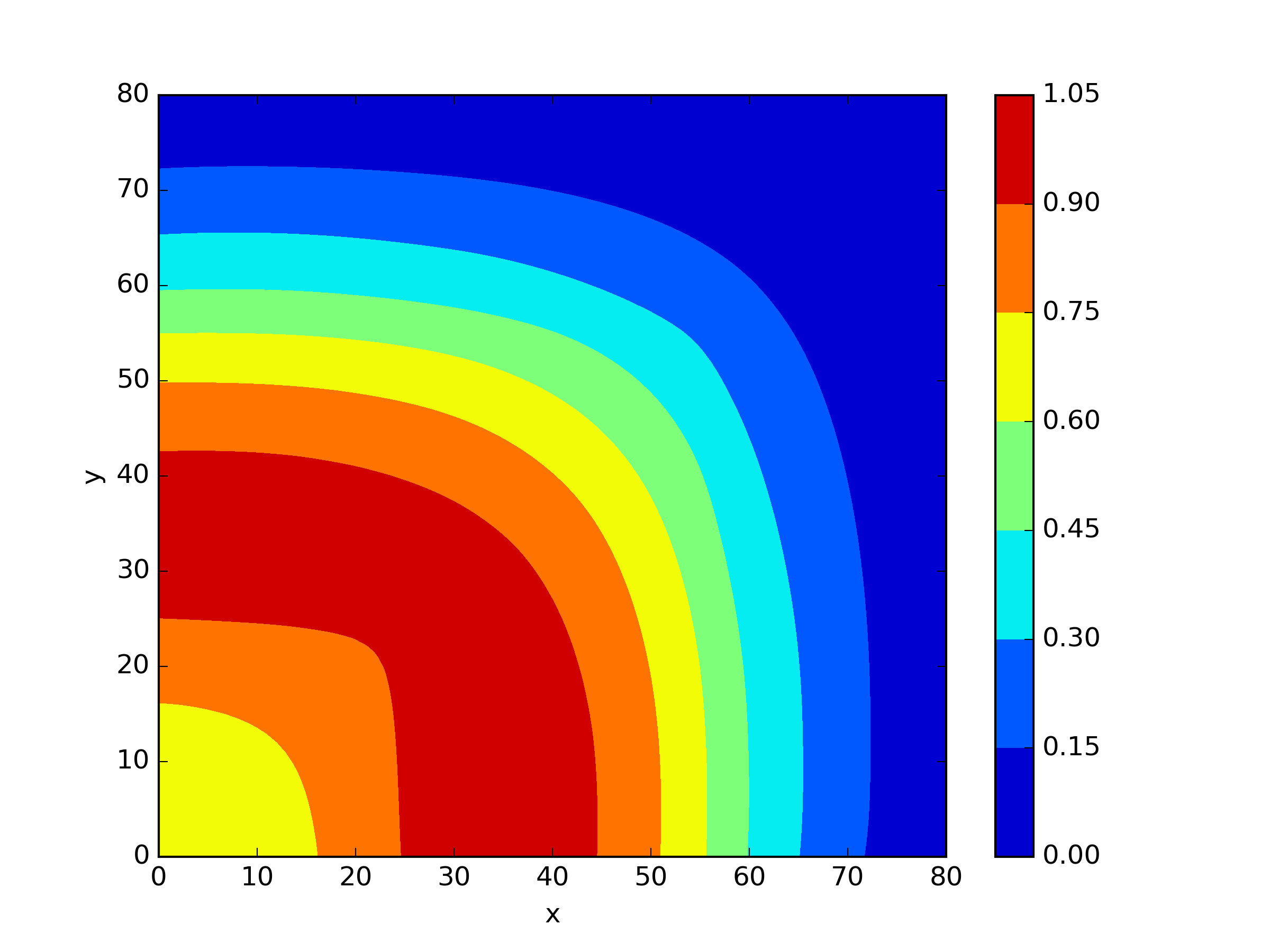}
		\end{minipage}
		\begin{minipage}{0.3\textwidth}
			\centering
			\includegraphics[width=4.7cm]{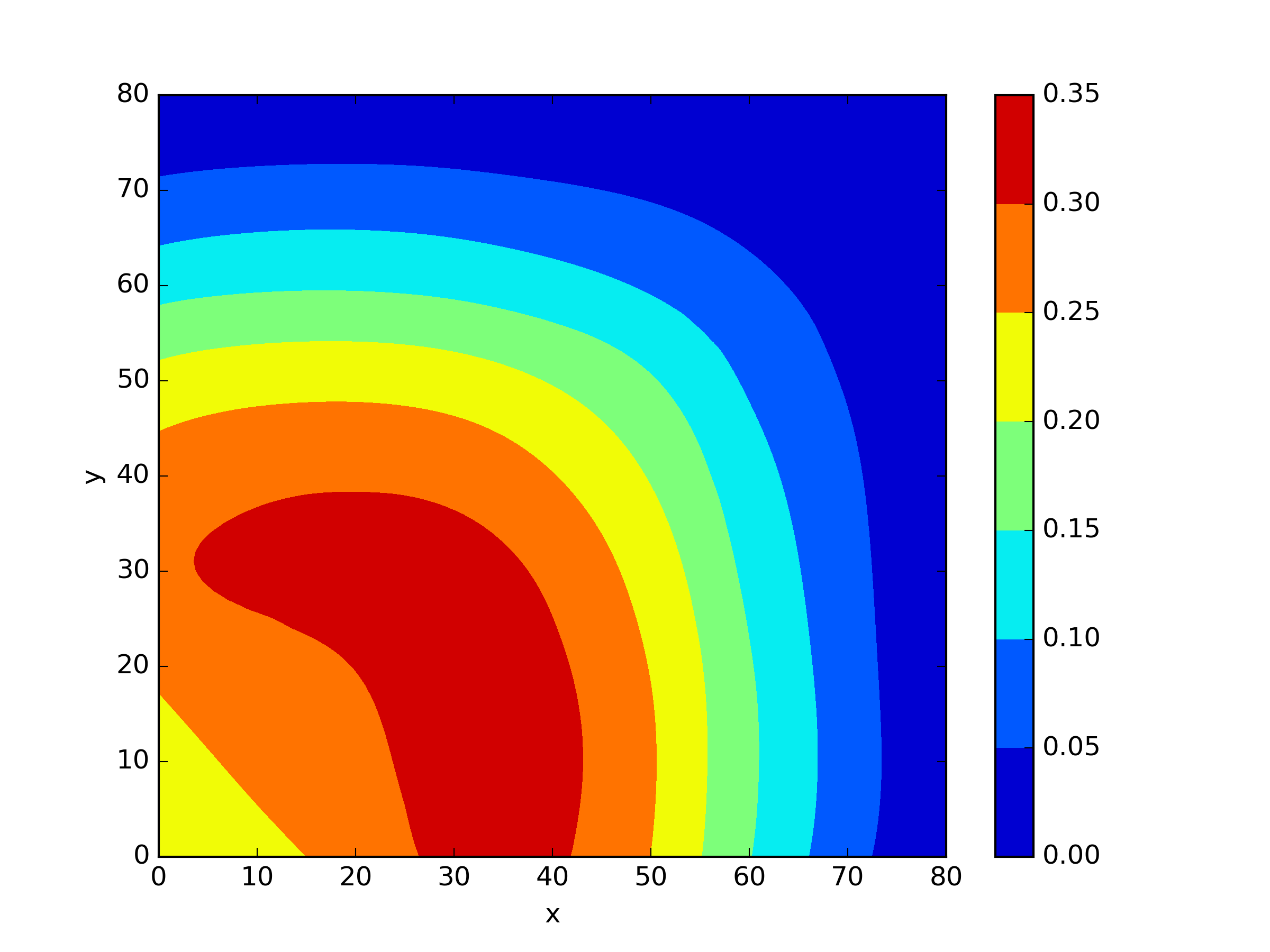}
		\end{minipage}
		
		\begin{minipage}{0.3\textwidth}
			\centering
			\includegraphics[width=4.7cm]{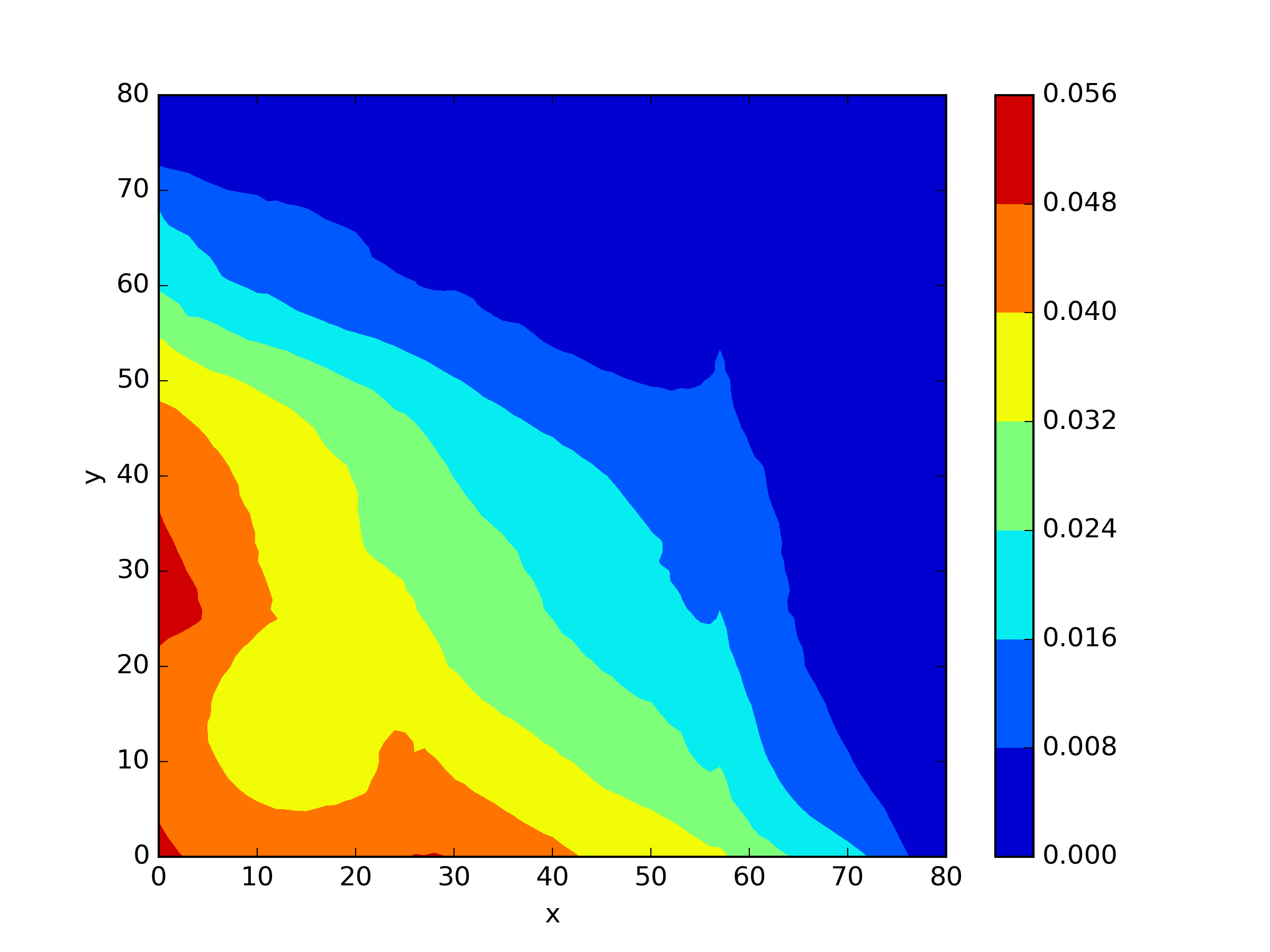}
		\end{minipage}
		\begin{minipage}{0.3\textwidth}
			\centering
			\includegraphics[width=4.7cm]{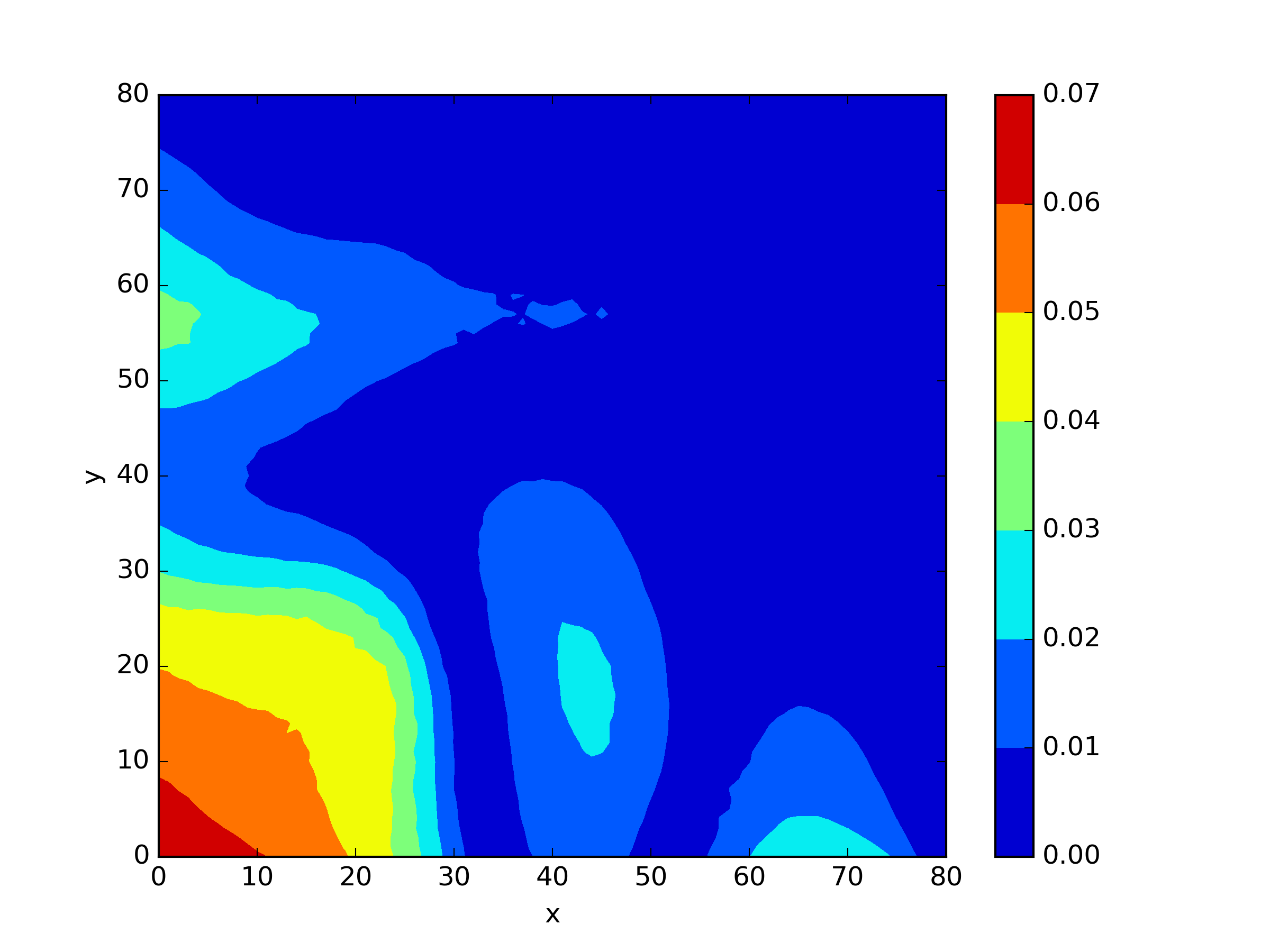}
		\end{minipage}
		\begin{minipage}{0.3\textwidth}
			\centering
			\includegraphics[width=4.7cm]{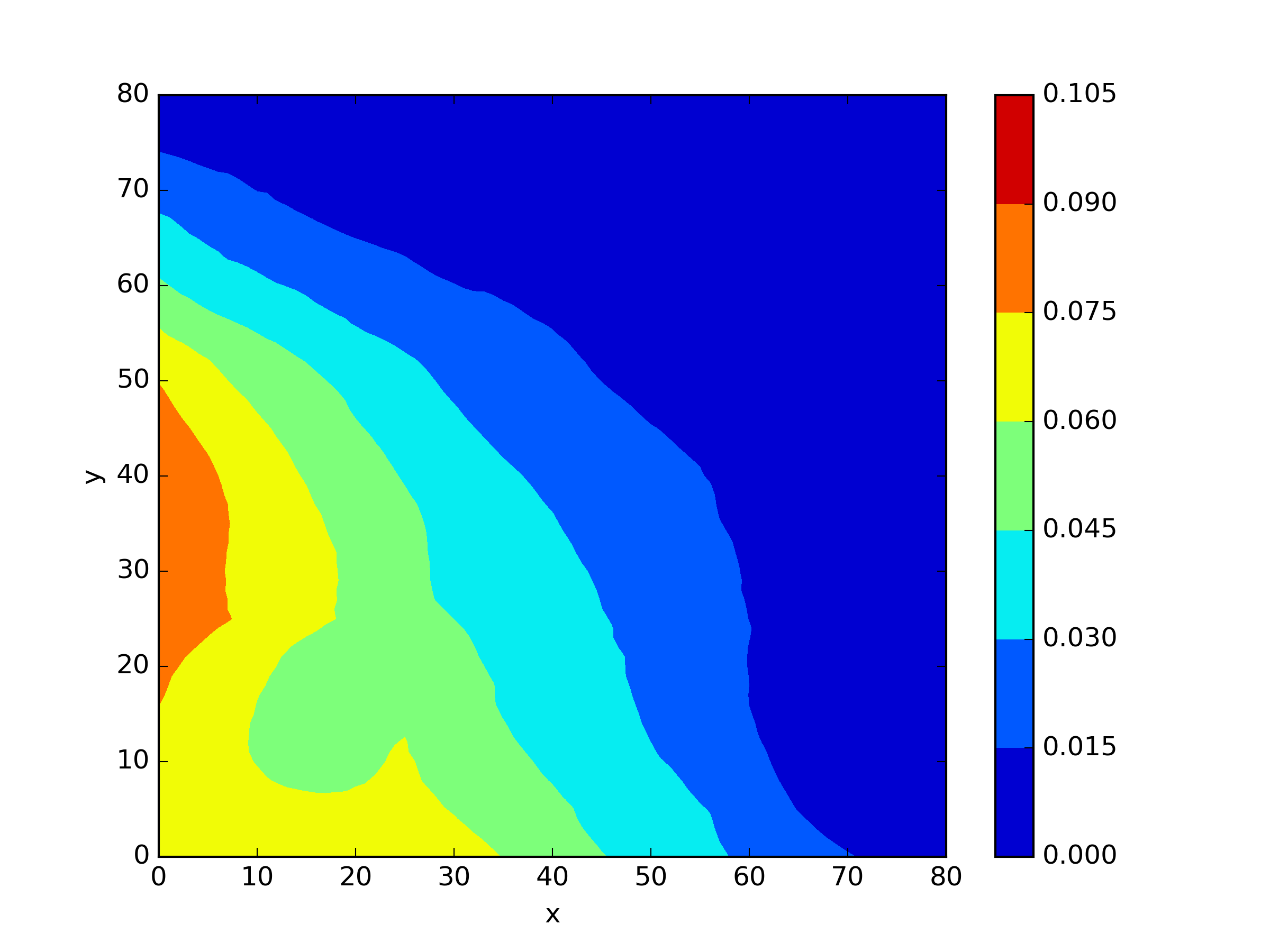}
		\end{minipage}
		
		\begin{minipage}{0.3\textwidth}
			\centering
			\includegraphics[width=4.7cm]{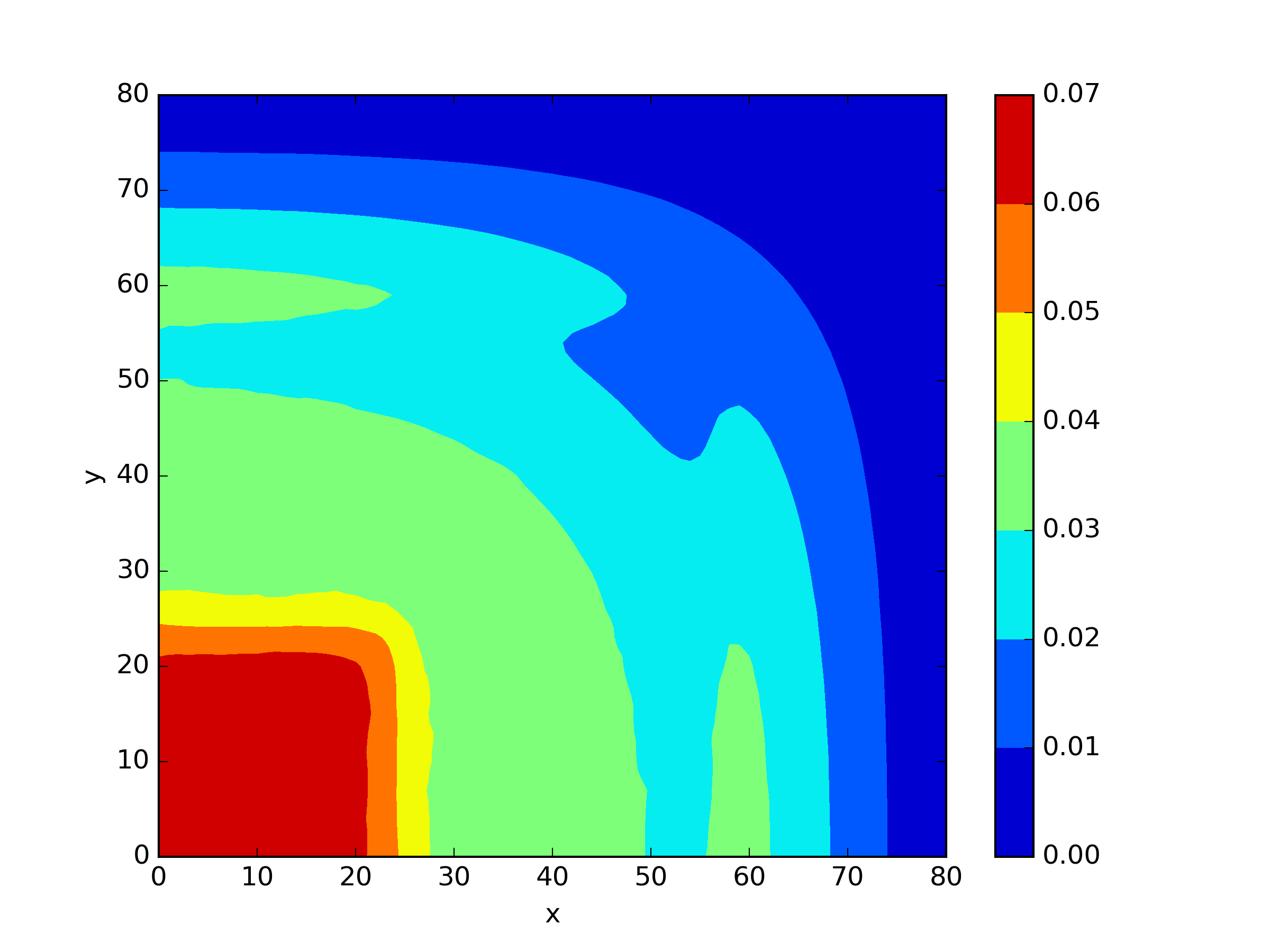}
		\end{minipage}
		\begin{minipage}{0.3\textwidth}
			\centering
			\includegraphics[width=4.7cm]{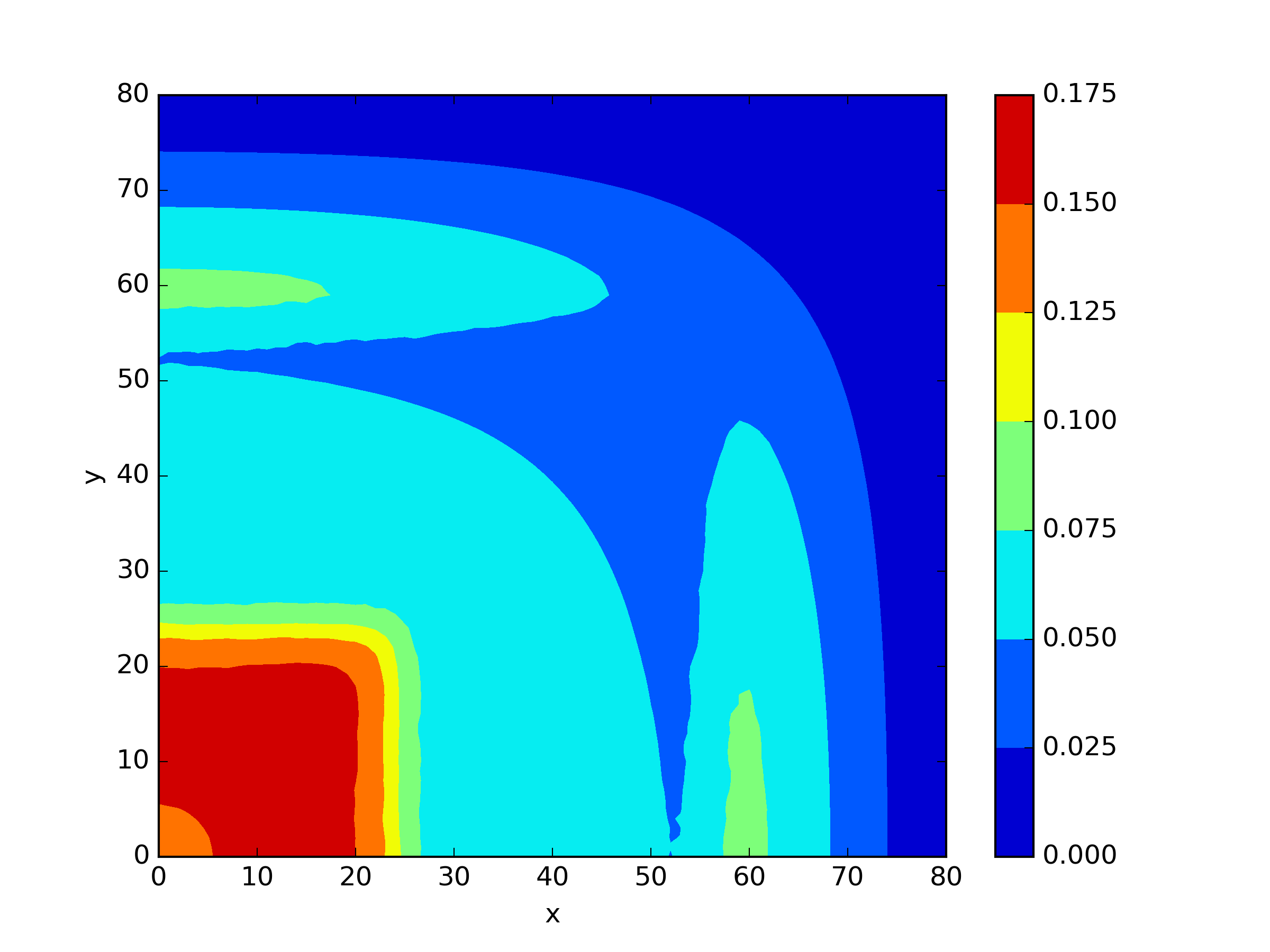}
		\end{minipage}
		\begin{minipage}{0.3\textwidth}
			\centering
			\includegraphics[width=4.7cm]{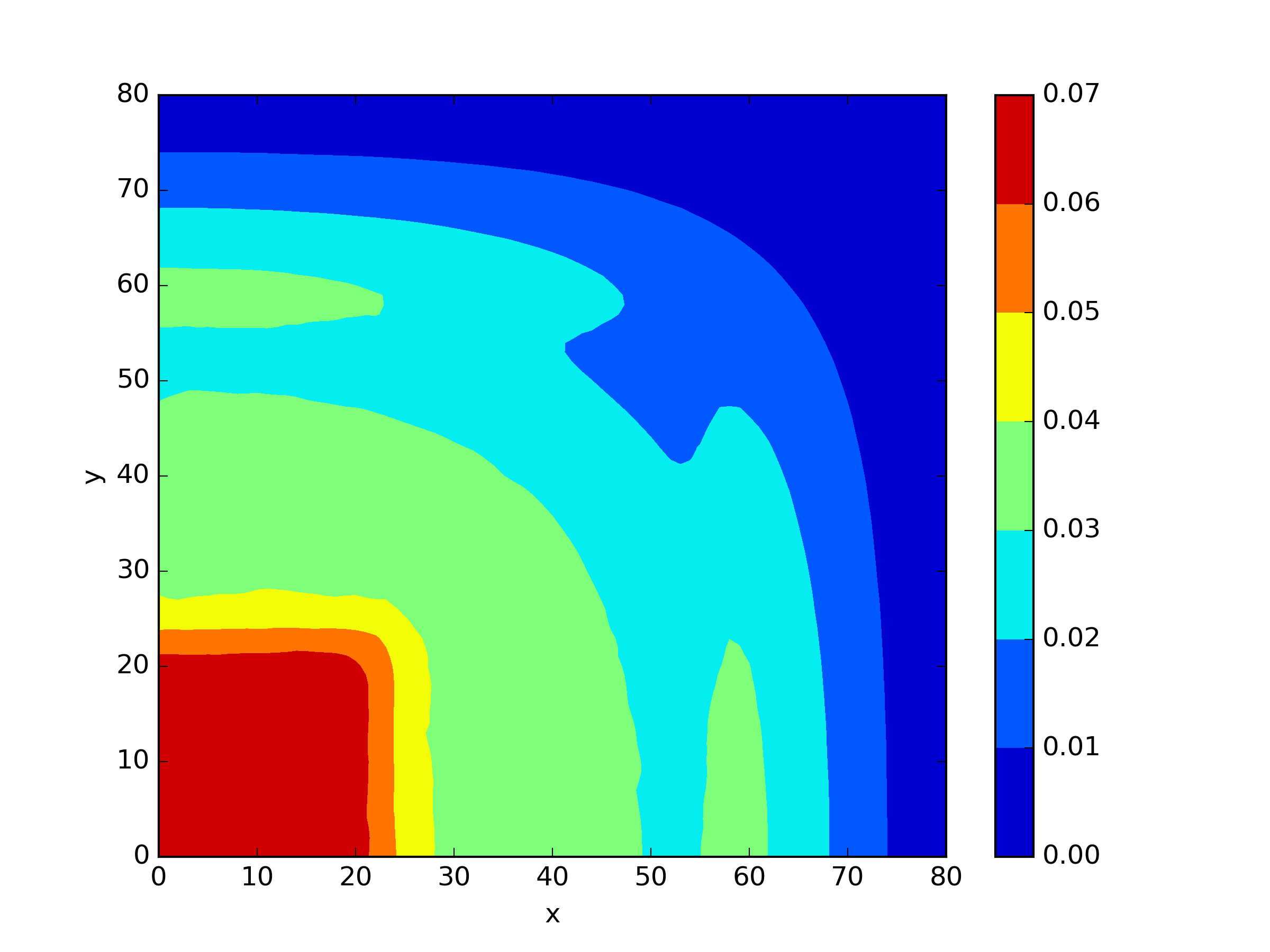}
		\end{minipage}
		
		\begin{minipage}{0.3\textwidth}
			\centering
			\includegraphics[width=4.7cm]{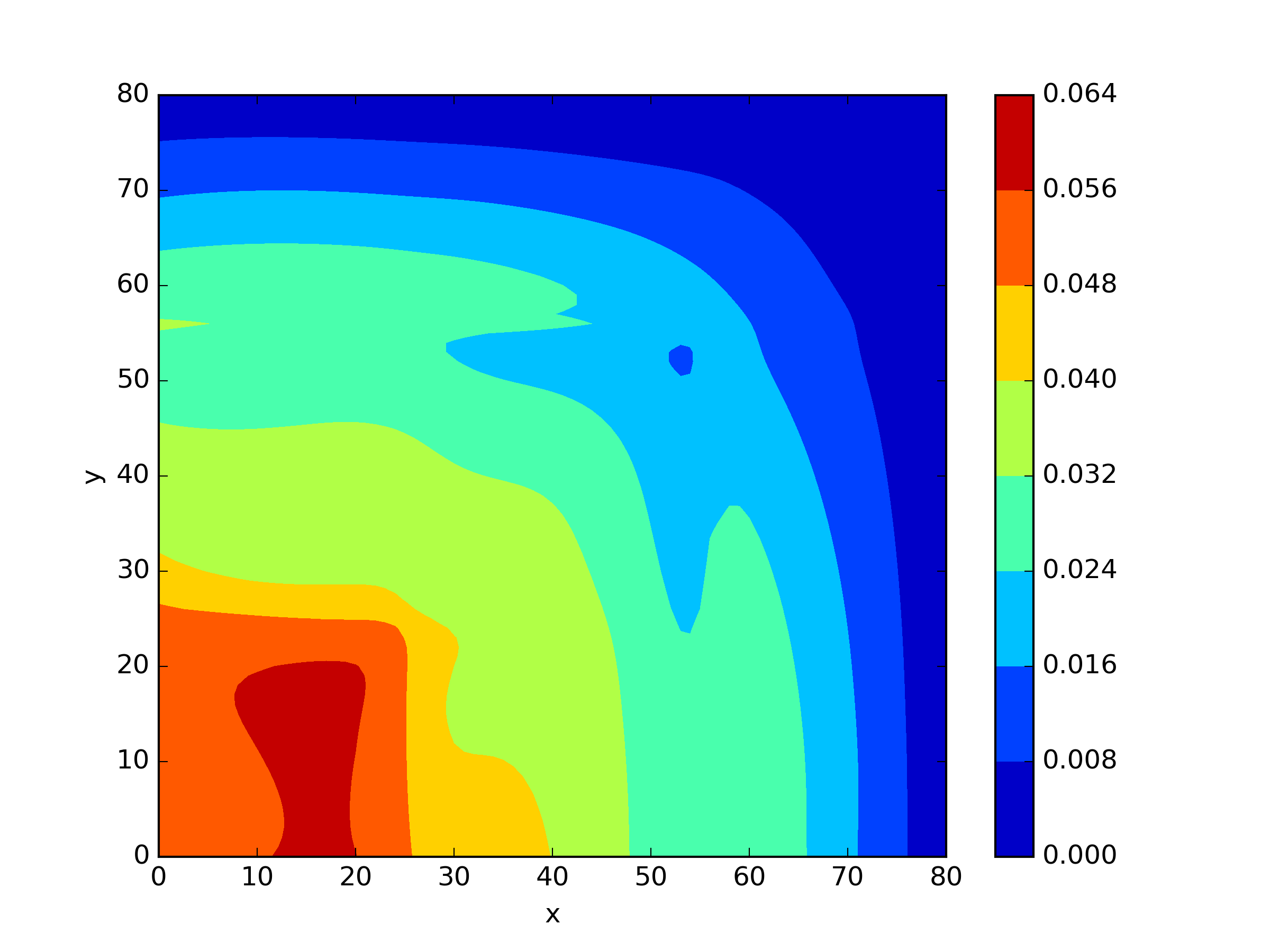}
		\end{minipage}
		\begin{minipage}{0.3\textwidth}
			\centering
			\includegraphics[width=4.7cm]{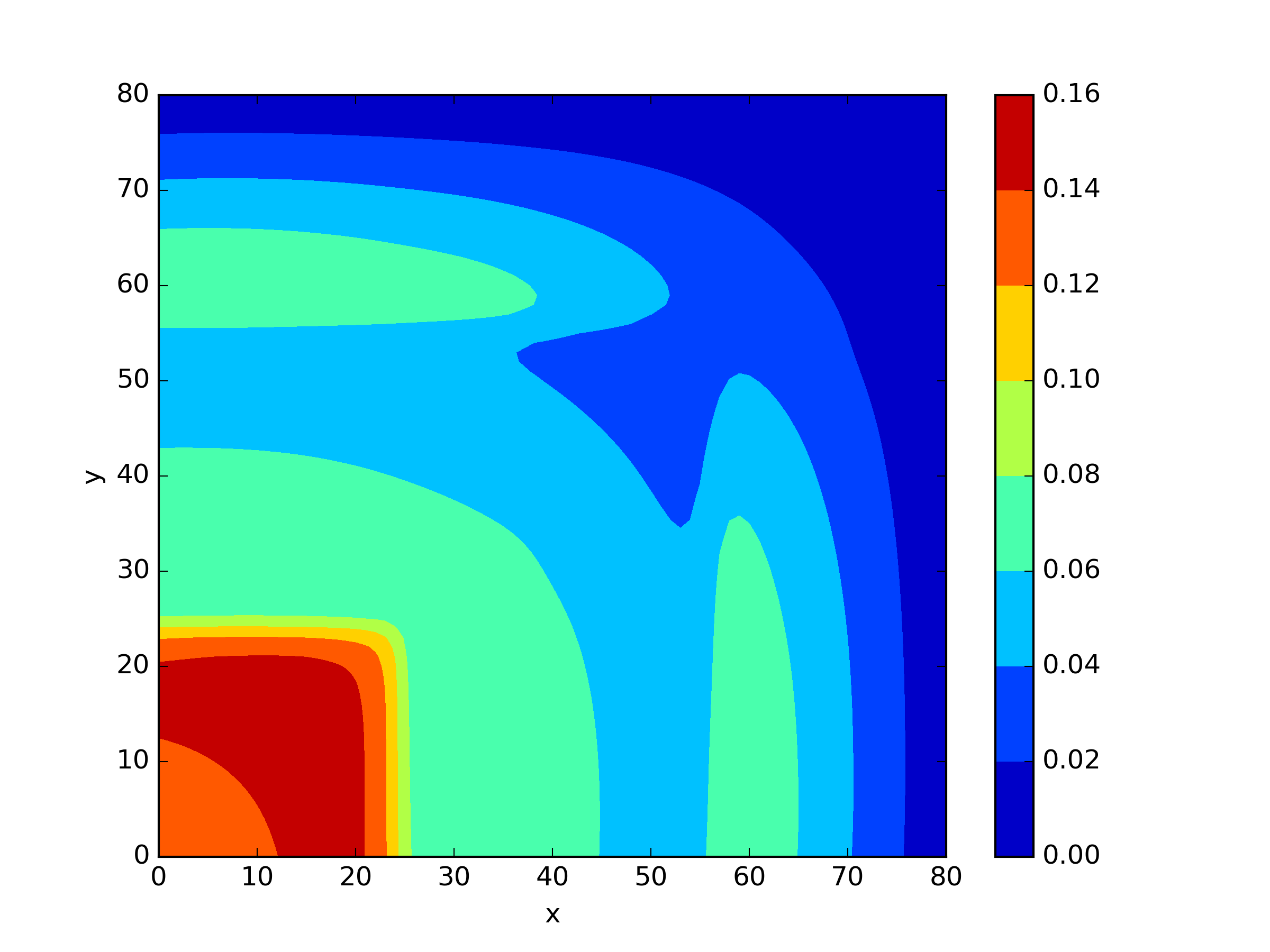}
		\end{minipage}
		\begin{minipage}{0.3\textwidth}
			\centering
			\includegraphics[width=4.7cm]{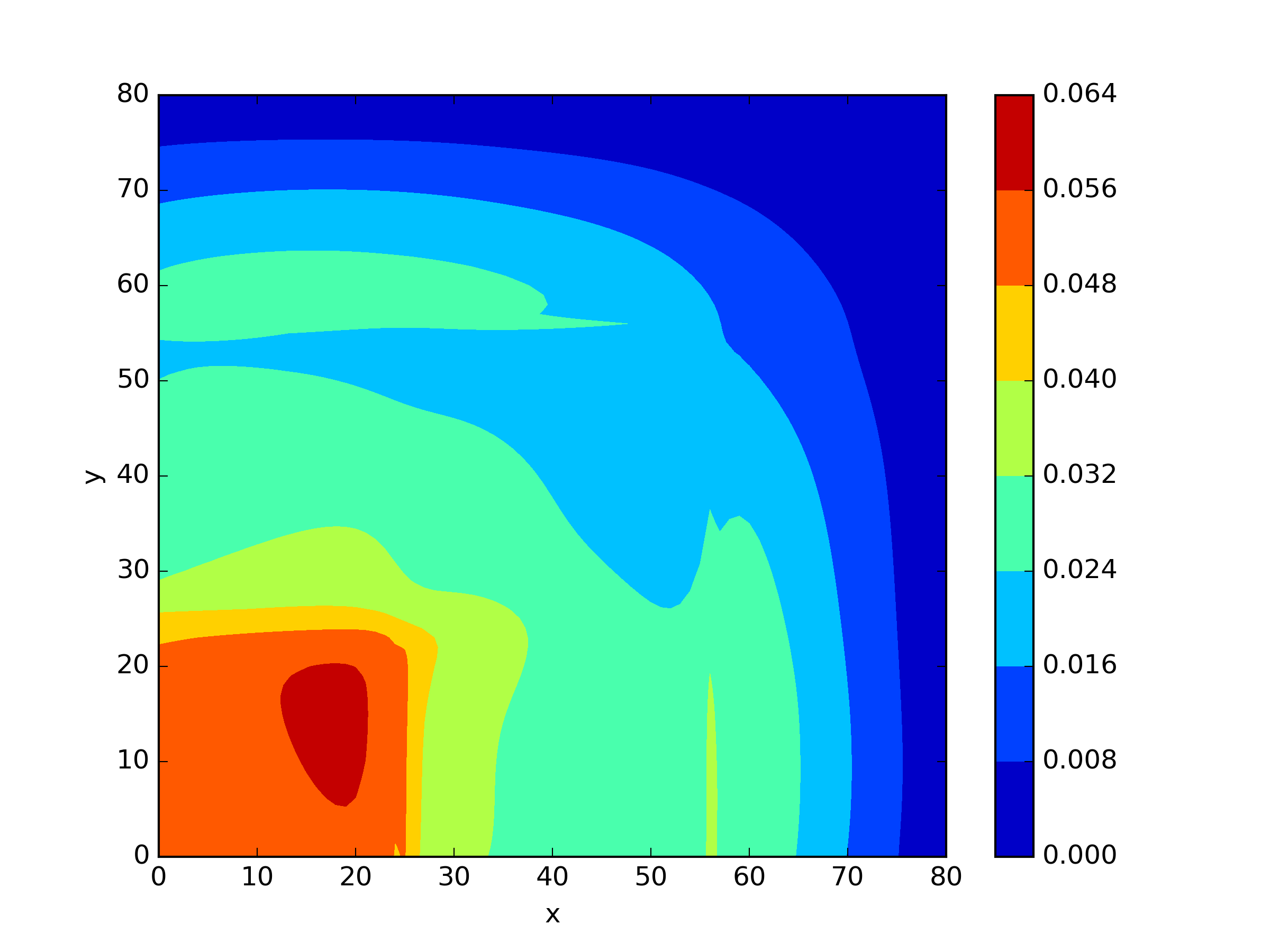}
		\end{minipage}
		
		\begin{minipage}{0.3\textwidth}
			\centering
			\includegraphics[width=4.7cm]{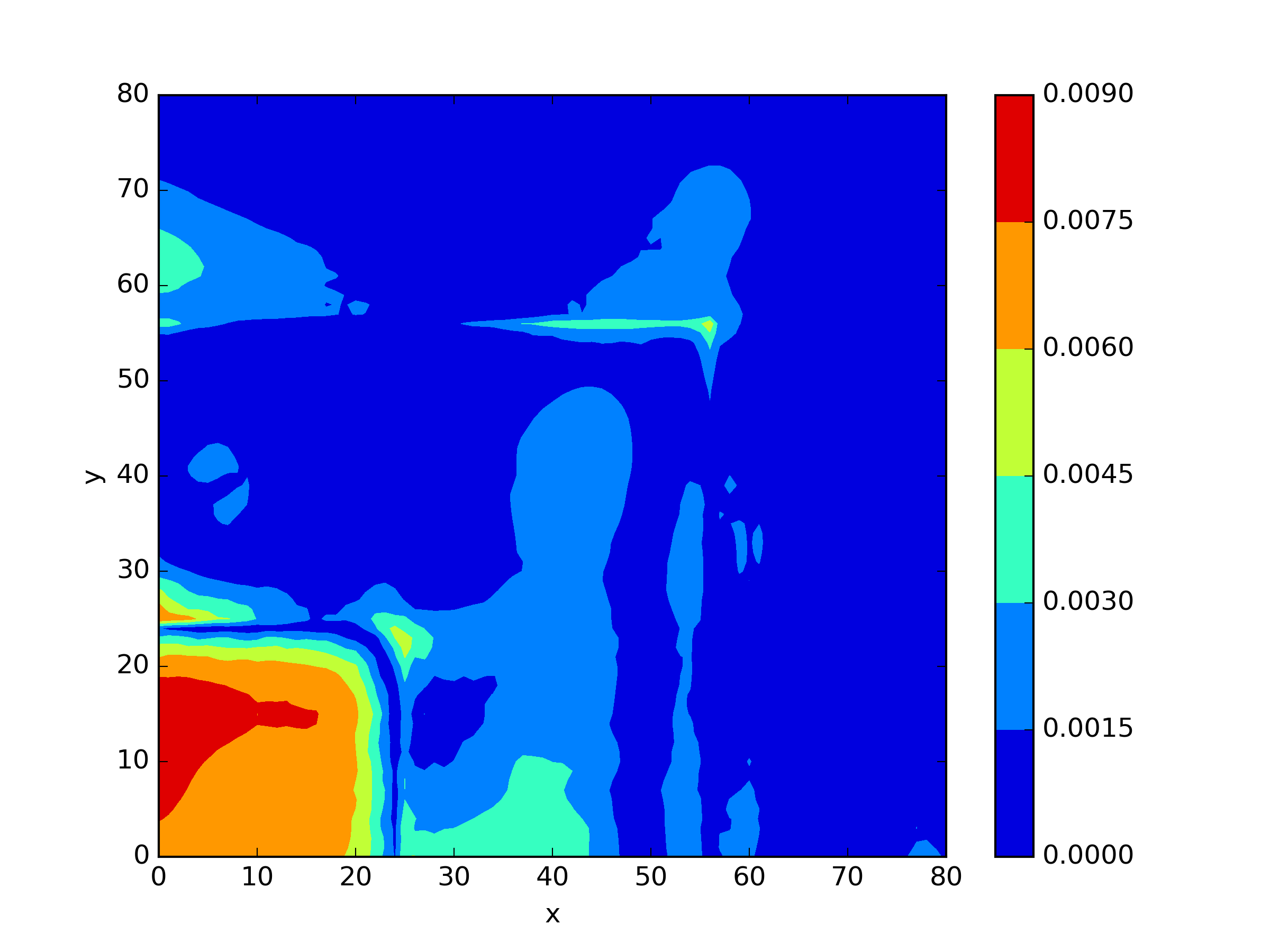}
		\end{minipage}
		\begin{minipage}{0.3\textwidth}
			\centering
			\includegraphics[width=4.7cm]{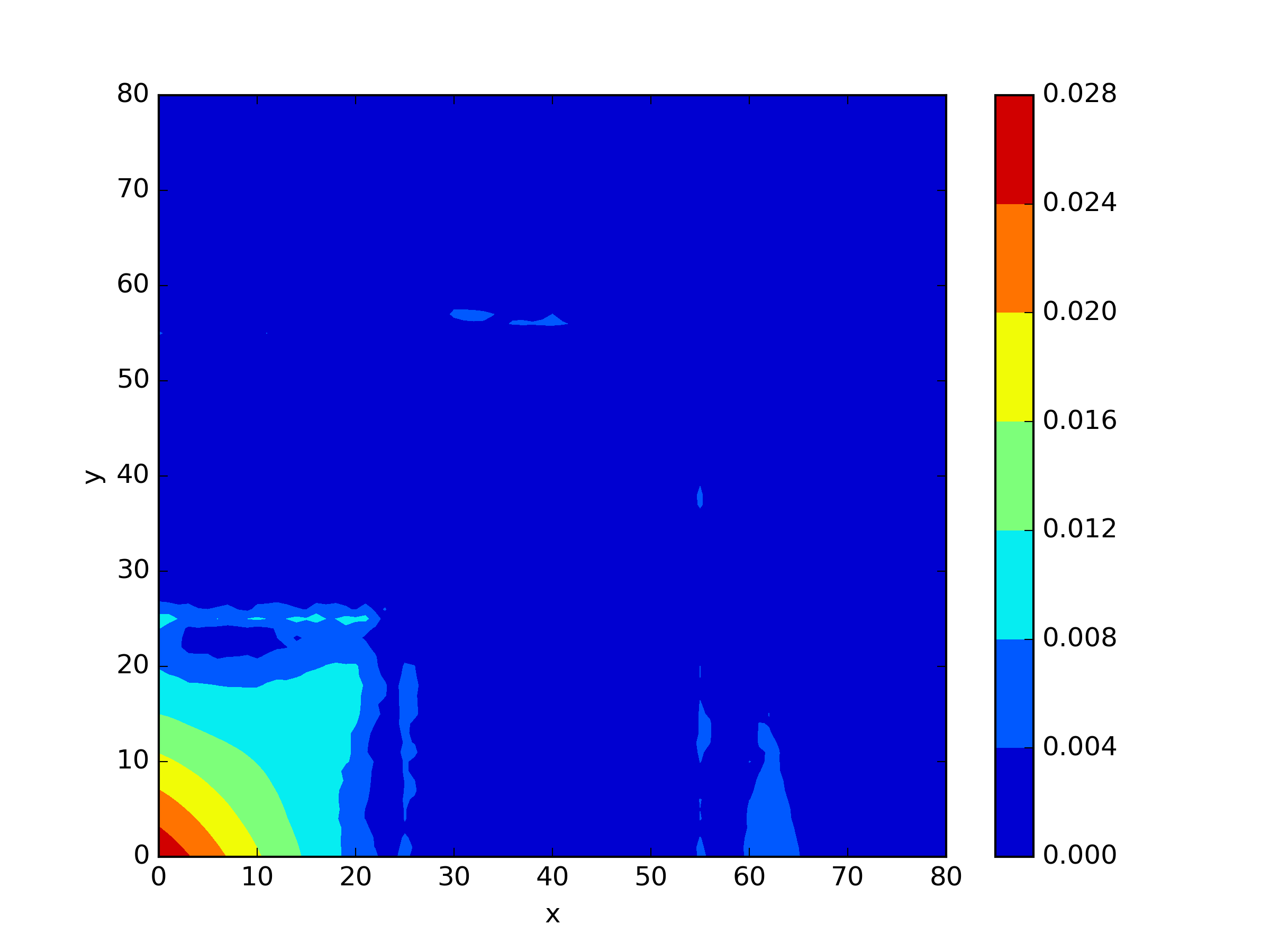}
		\end{minipage}
		\begin{minipage}{0.3\textwidth}
			\centering
			\includegraphics[width=4.7cm]{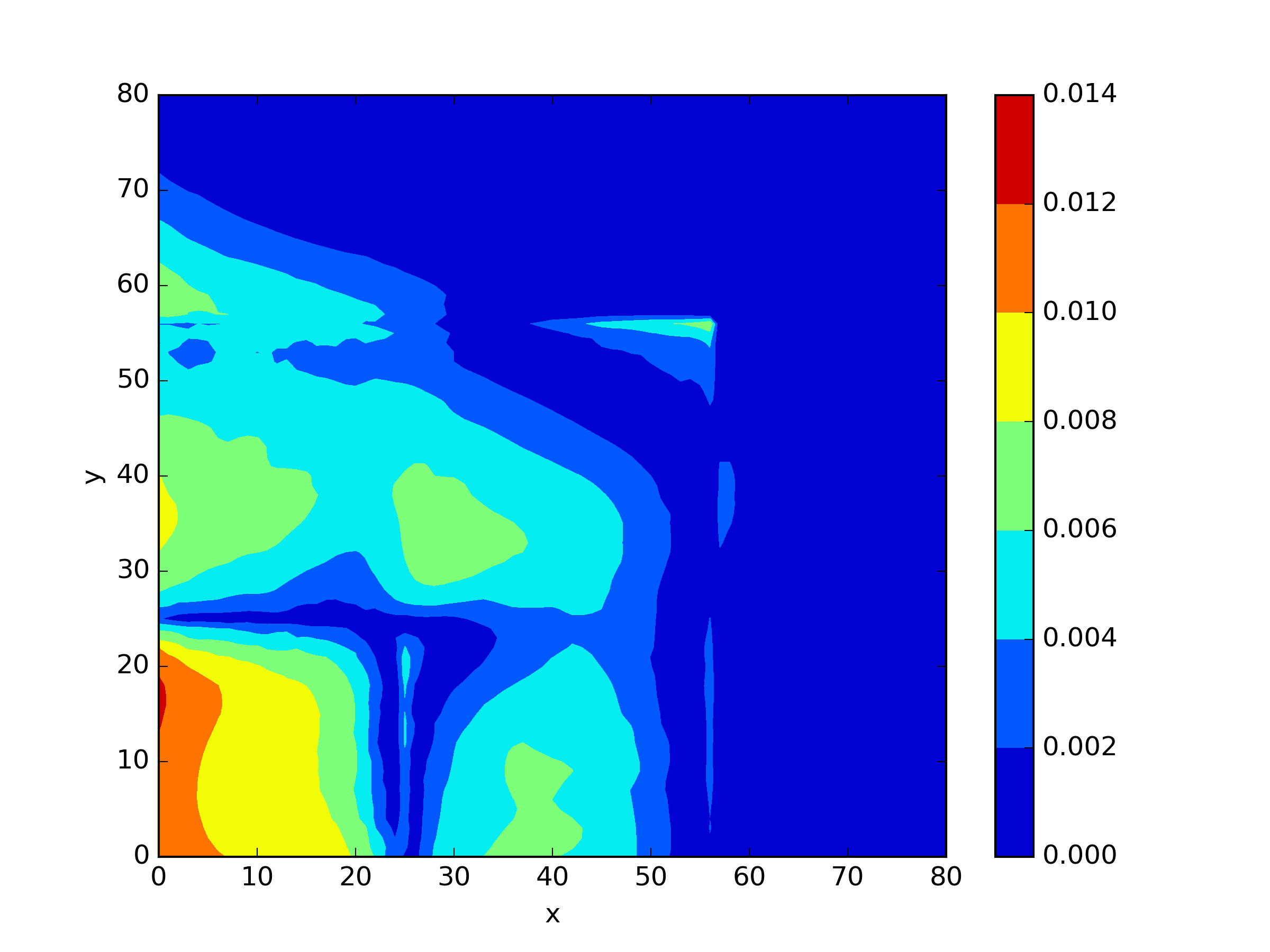}
		\end{minipage}
		\caption{The first column, second column, and third column in the figure depict the results of the solution at z = 20, 70, and 120, respectively. The six rows of images, from top to bottom, represent the reference solution of $\phi_1$, the neural network solution of $\phi_1$, the absolute error between the two, the reference solution of $\phi_2$, the neural network solution of $\phi_2$, and the absolute error between the two.}
		\label{pic-TWIGL-3D}
	\end{figure*}
	
	\begin{figure*}[h]
		\begin{minipage}{0.24\textwidth}
			\centering
			\includegraphics[width=4.25cm]{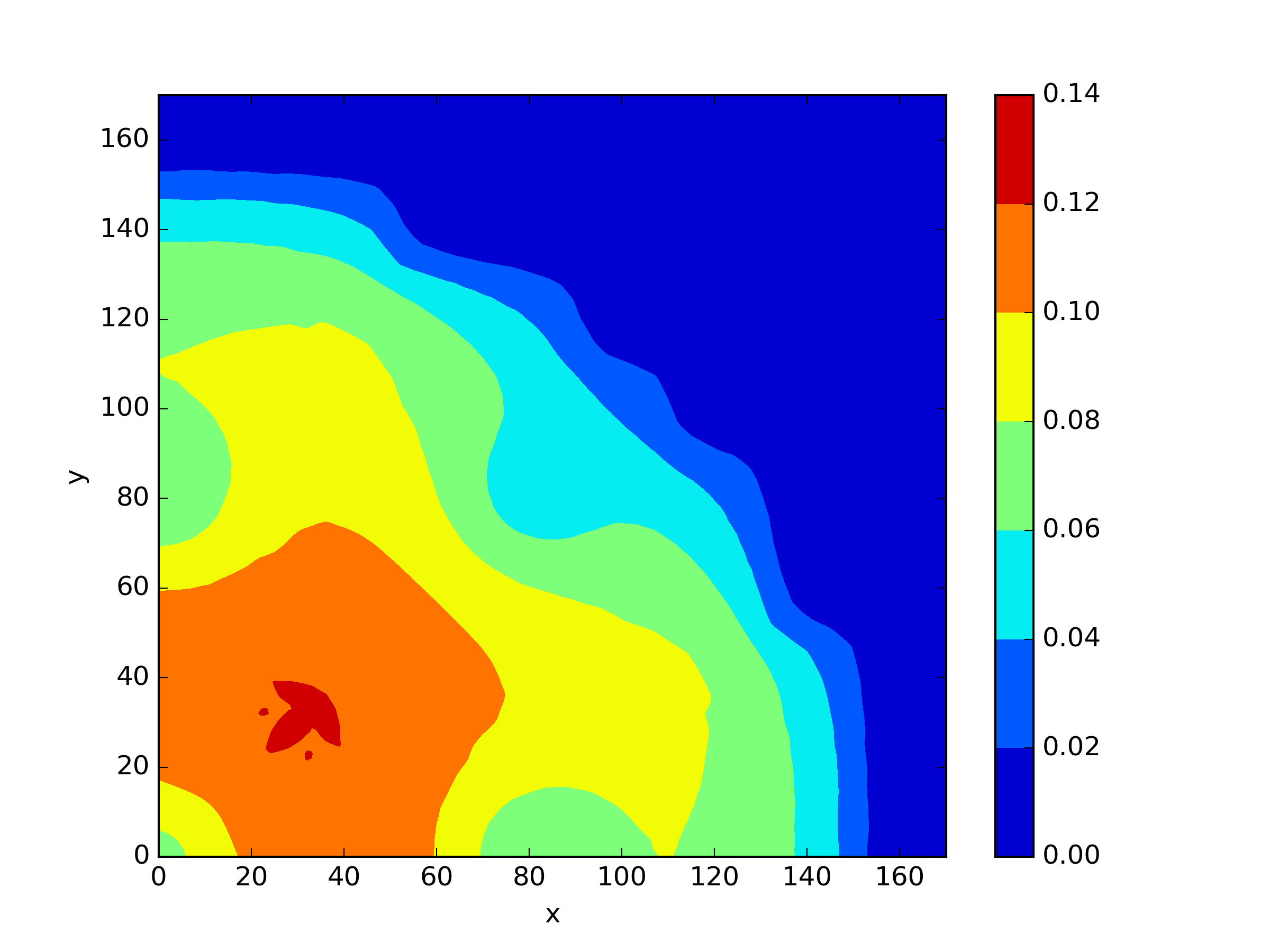}
		\end{minipage}
		\begin{minipage}{0.24\textwidth}
			\centering
			\includegraphics[width=4.25cm]{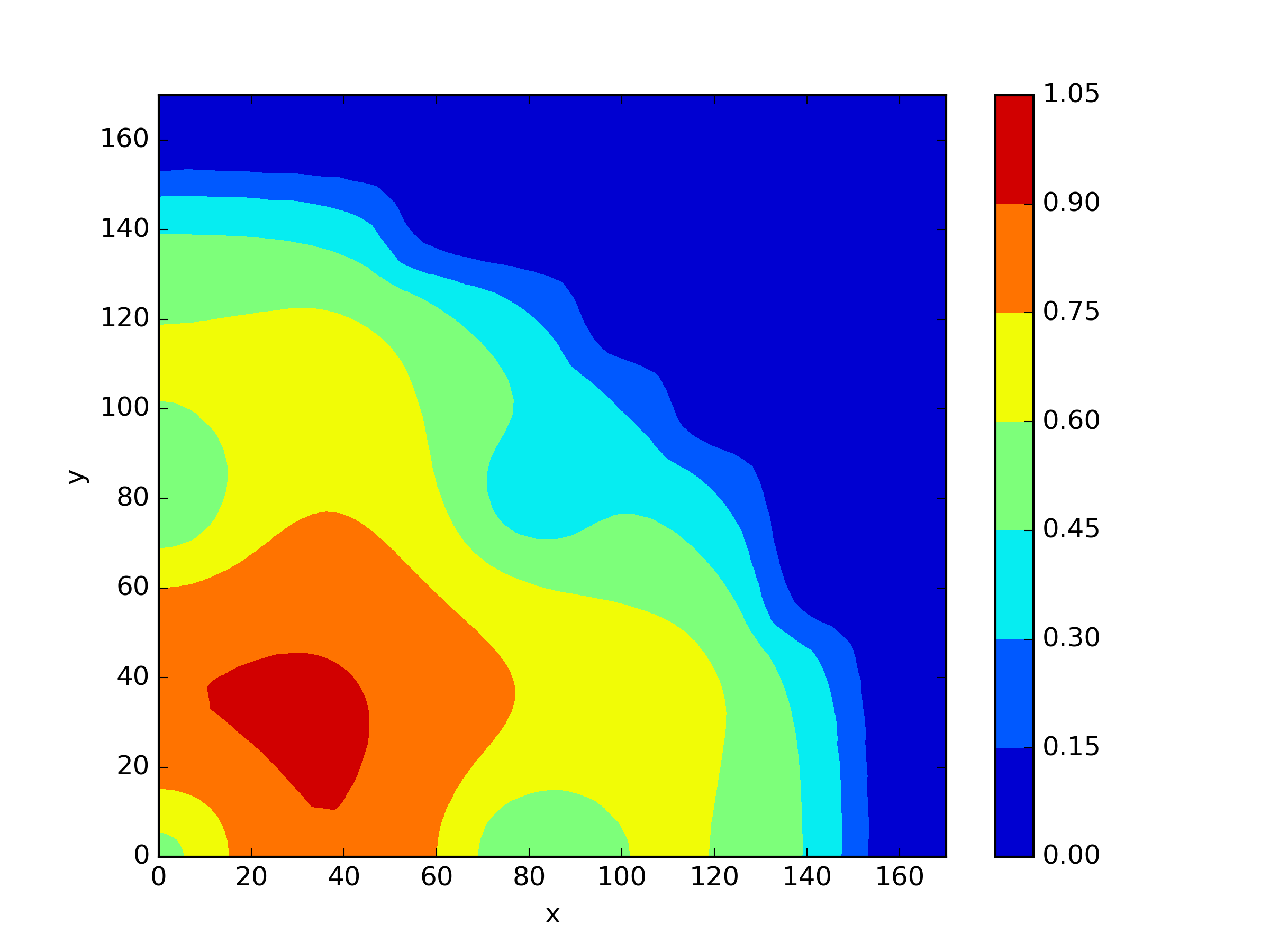}
		\end{minipage}
		\begin{minipage}{0.24\textwidth}
			\centering
			\includegraphics[width=4.25cm]{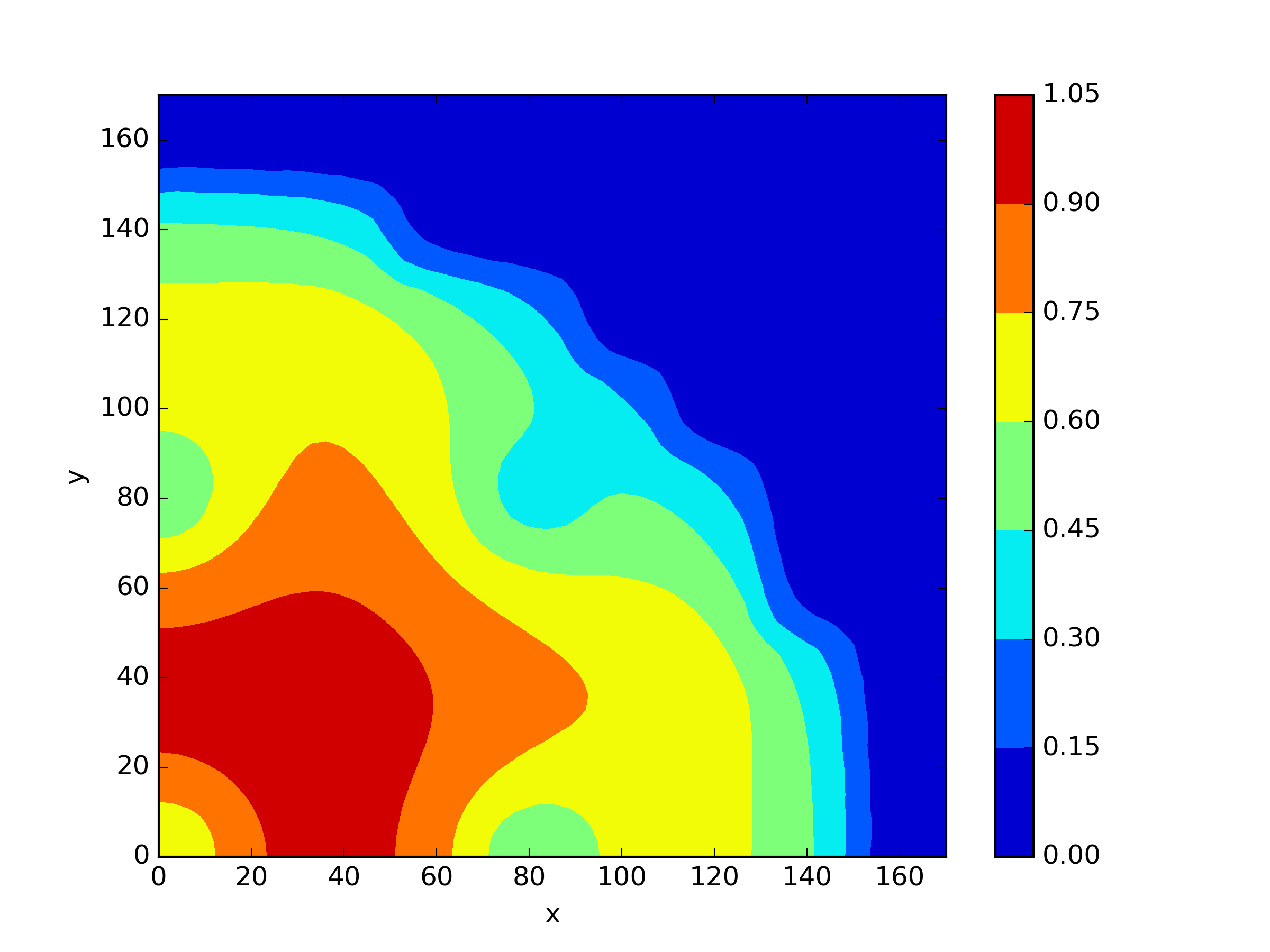}
		\end{minipage}
		\begin{minipage}{0.24\textwidth}
			\centering
			\includegraphics[width=4.25cm]{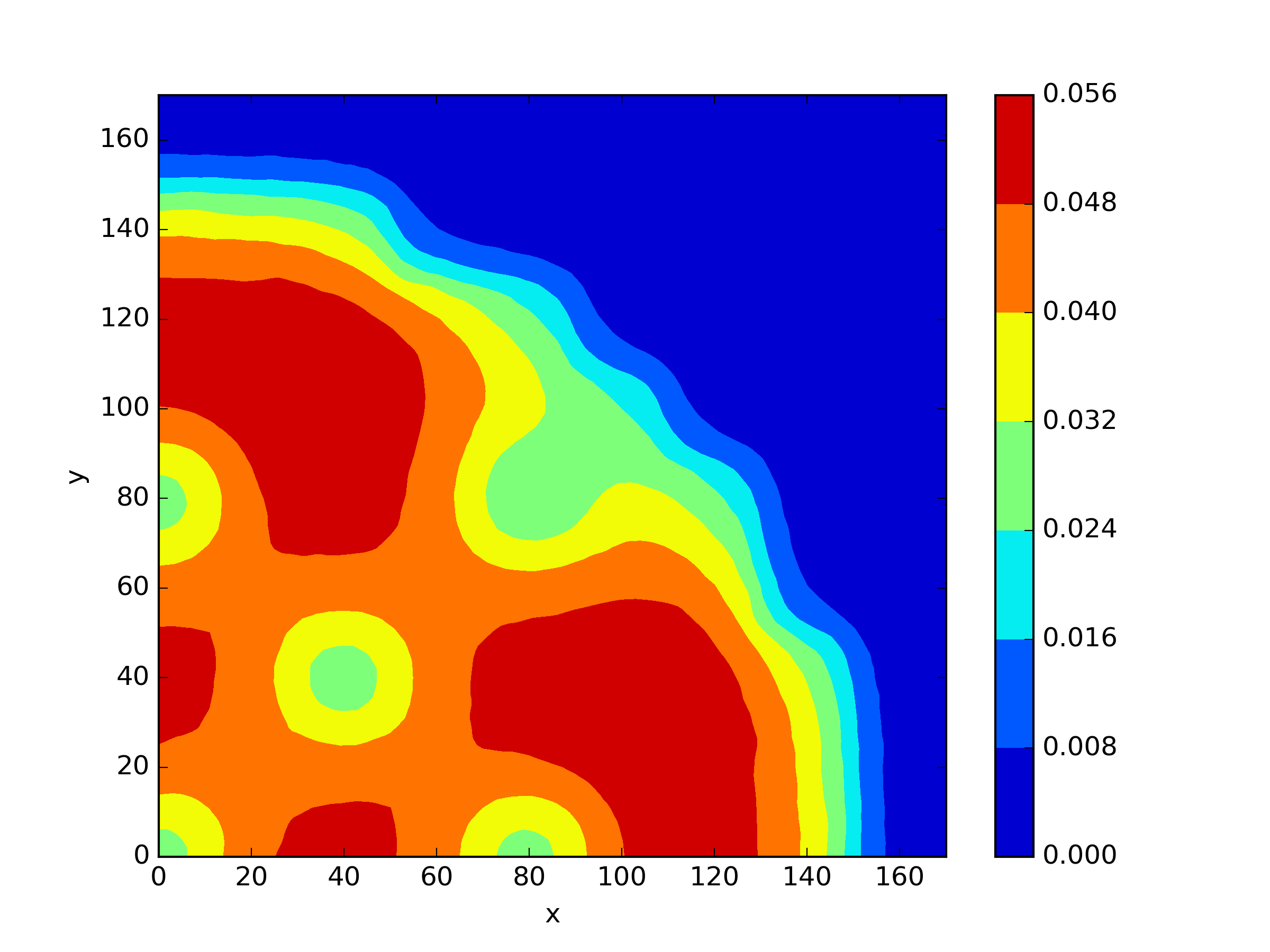}
		\end{minipage}
		
		\begin{minipage}{0.24\textwidth}
			\centering
			\includegraphics[width=4.25cm]{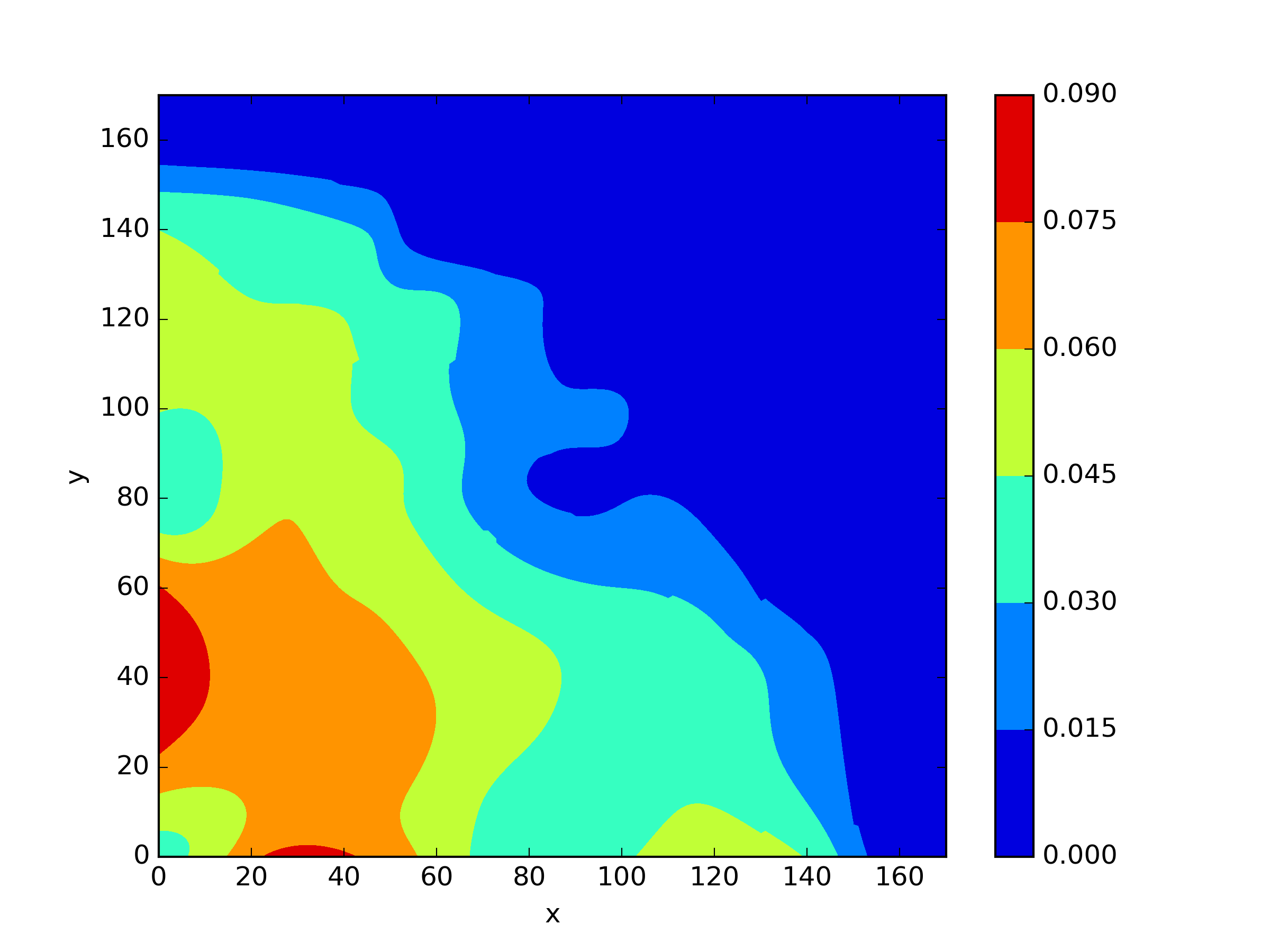}
		\end{minipage}
		\begin{minipage}{0.24\textwidth}
			\centering
			\includegraphics[width=4.25cm]{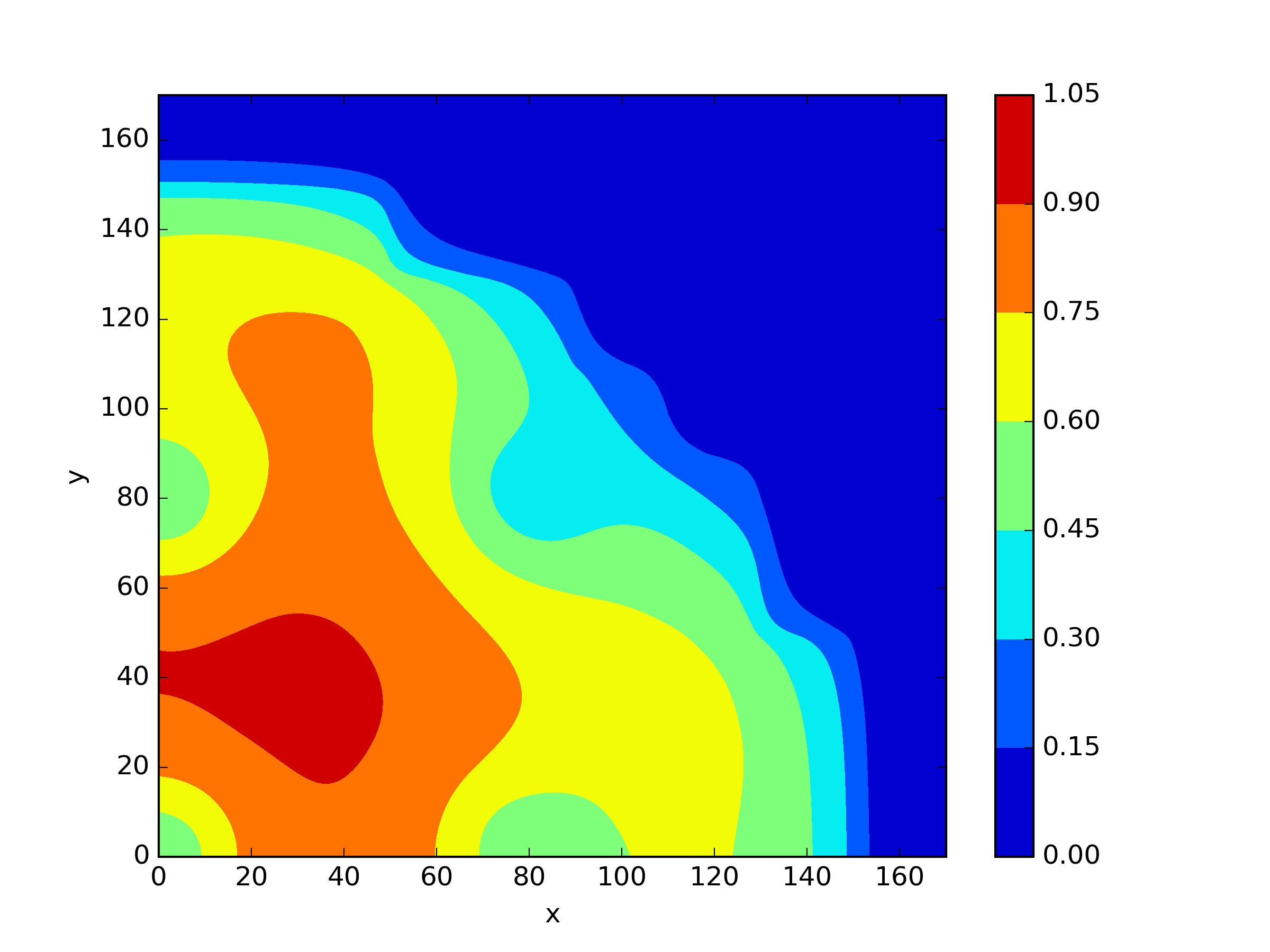}
		\end{minipage}
		\begin{minipage}{0.24\textwidth}
			\centering
			\includegraphics[width=4.25cm]{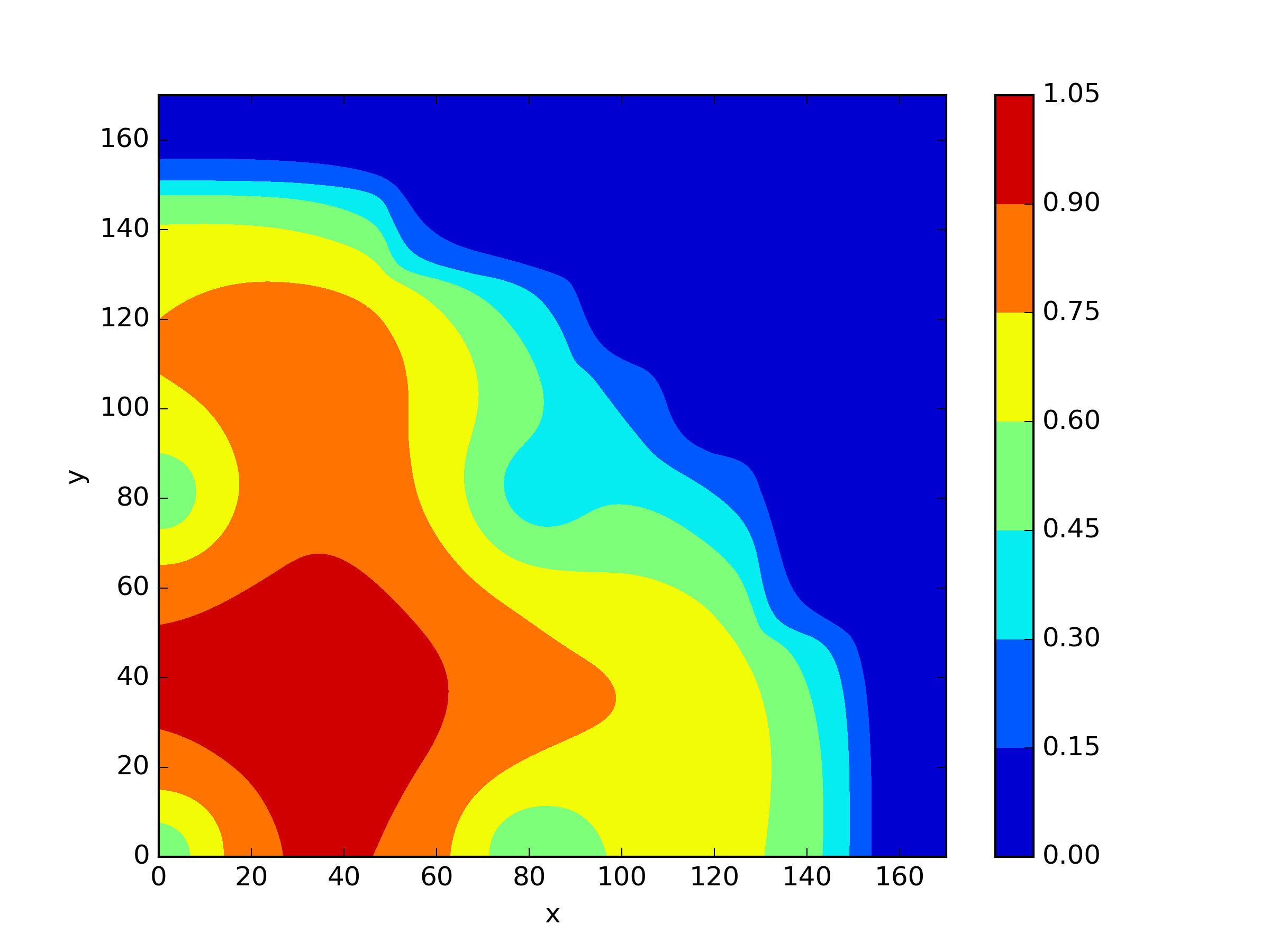}
		\end{minipage}
		\begin{minipage}{0.24\textwidth}
			\centering
			\includegraphics[width=4.25cm]{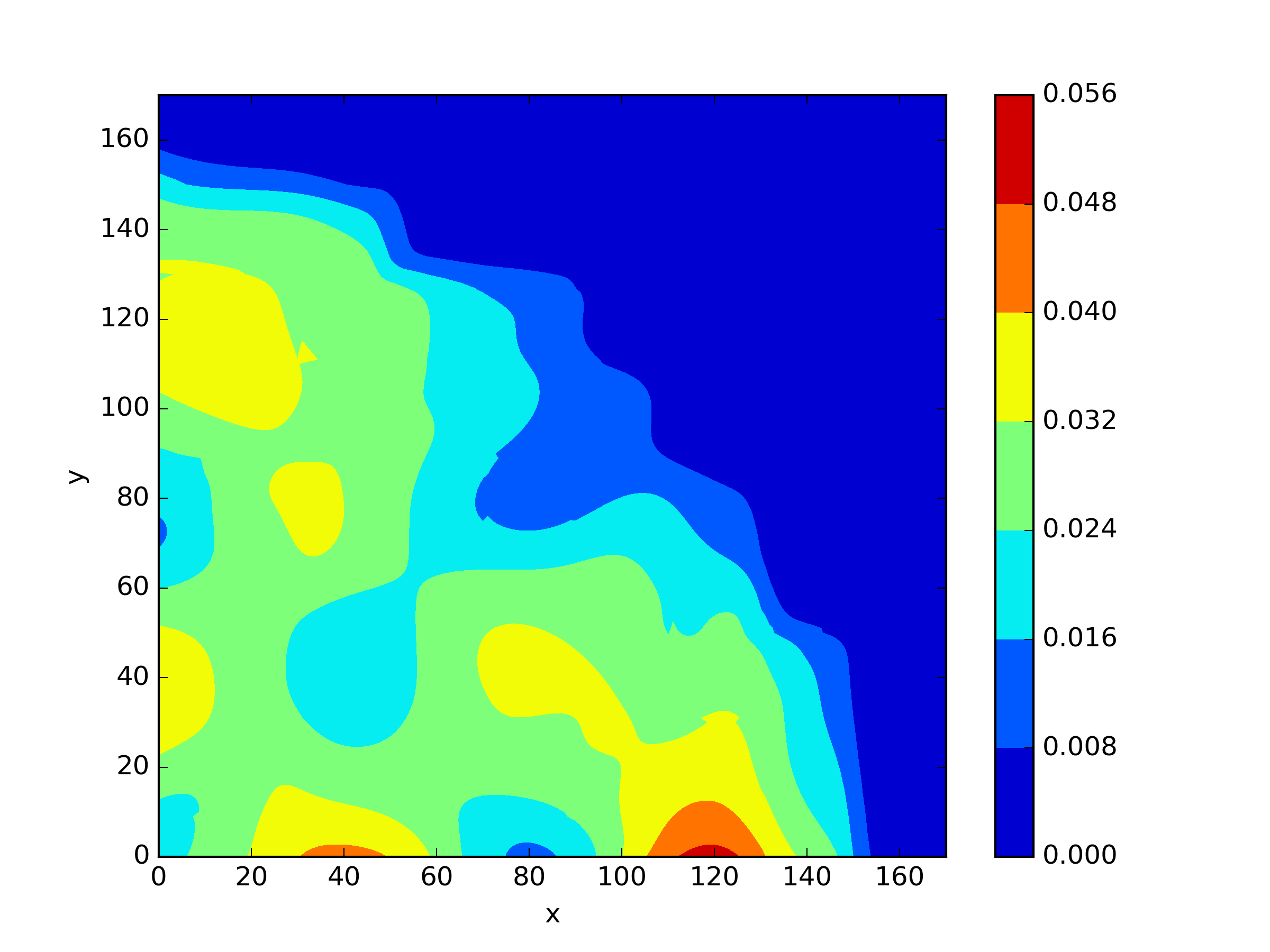}
		\end{minipage}
		
		\begin{minipage}{0.24\textwidth}
			\centering
			\includegraphics[width=4.25cm]{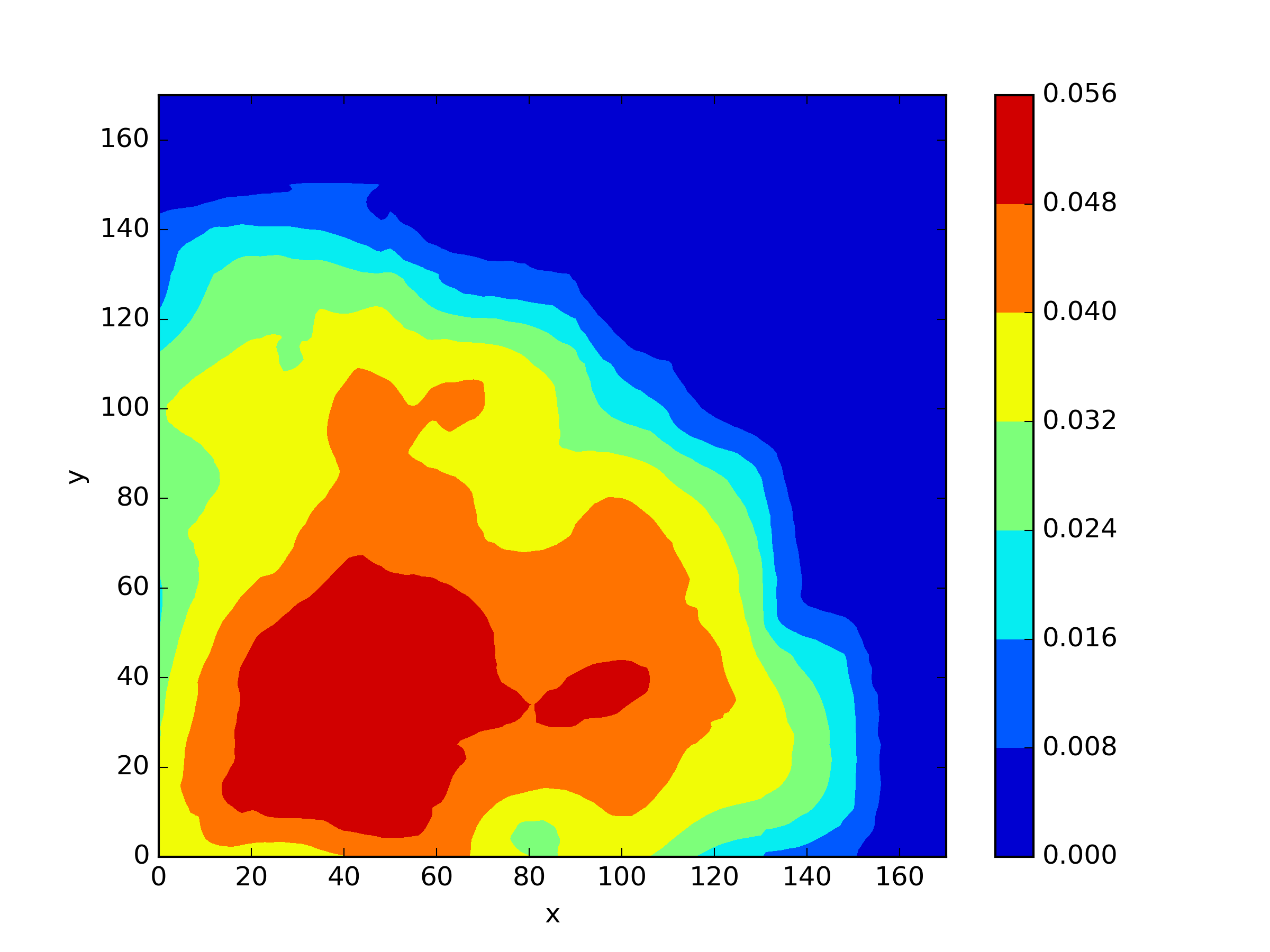}
		\end{minipage}
		\begin{minipage}{0.24\textwidth}
			\centering
			\includegraphics[width=4.25cm]{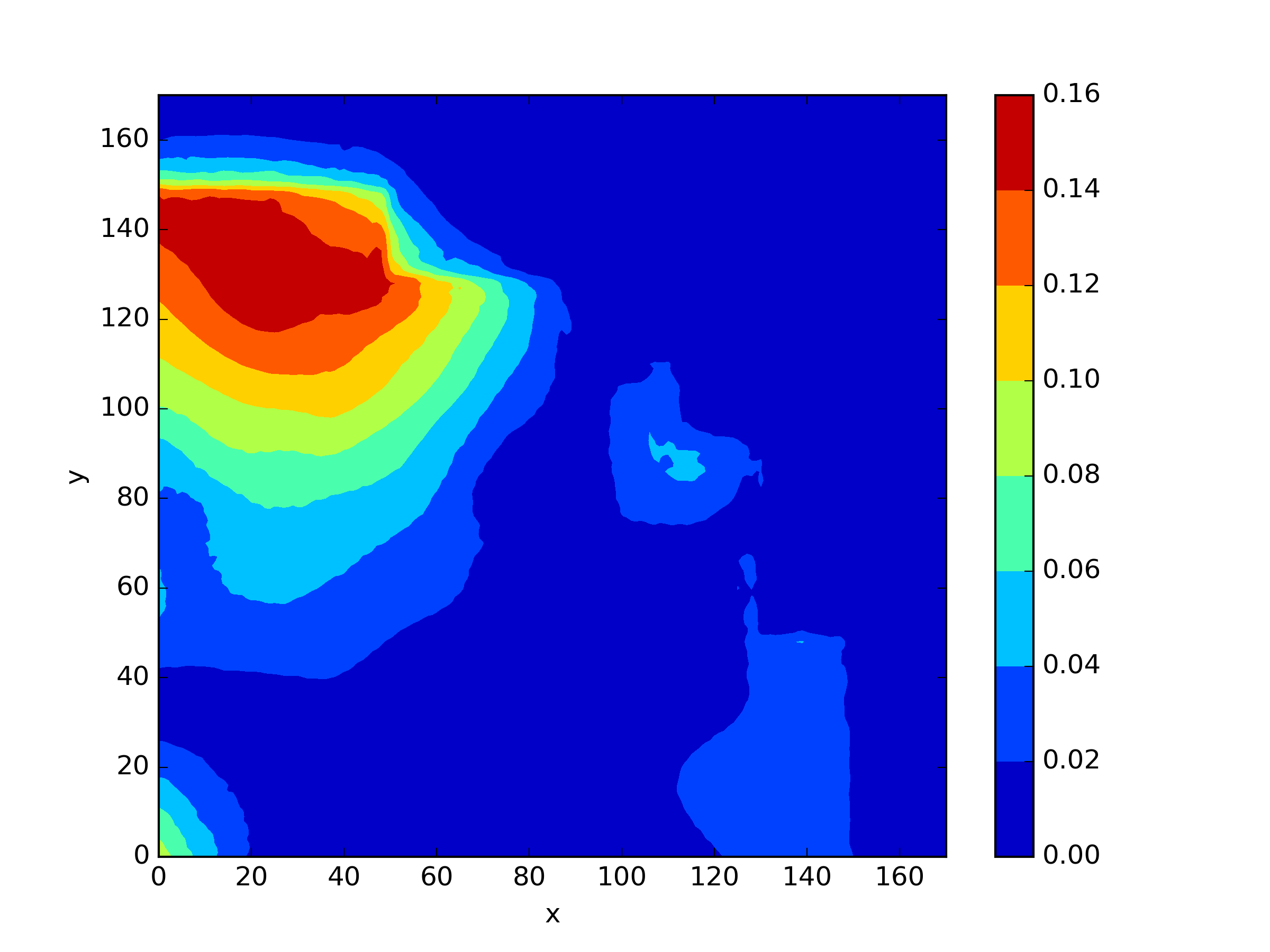}
		\end{minipage}
		\begin{minipage}{0.24\textwidth}
			\centering
			\includegraphics[width=4.25cm]{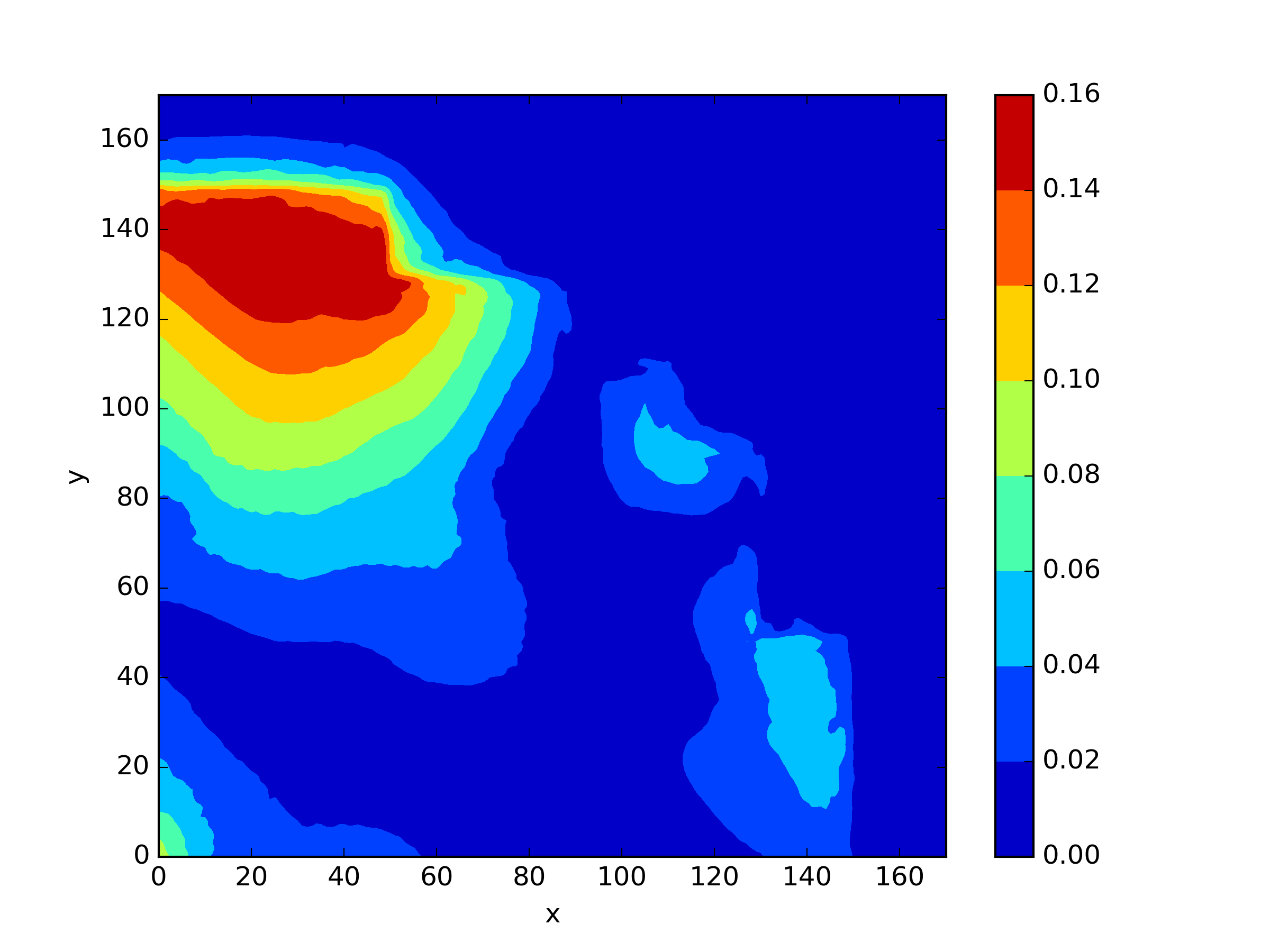}
		\end{minipage}
		\begin{minipage}{0.24\textwidth}
			\centering
			\includegraphics[width=4.25cm]{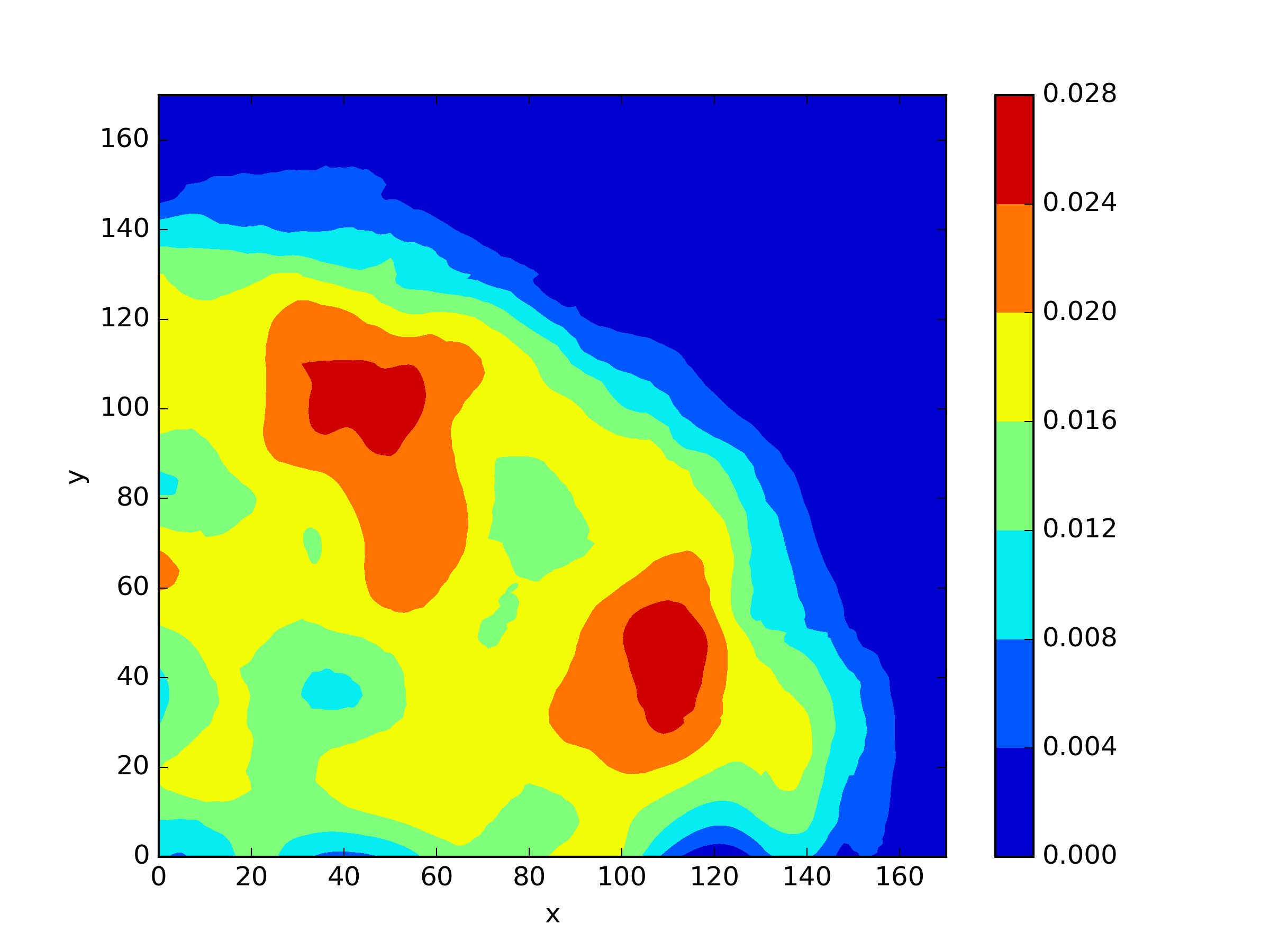}
		\end{minipage}
		
		\begin{minipage}{0.24\textwidth}
			\centering
			\includegraphics[width=4.25cm]{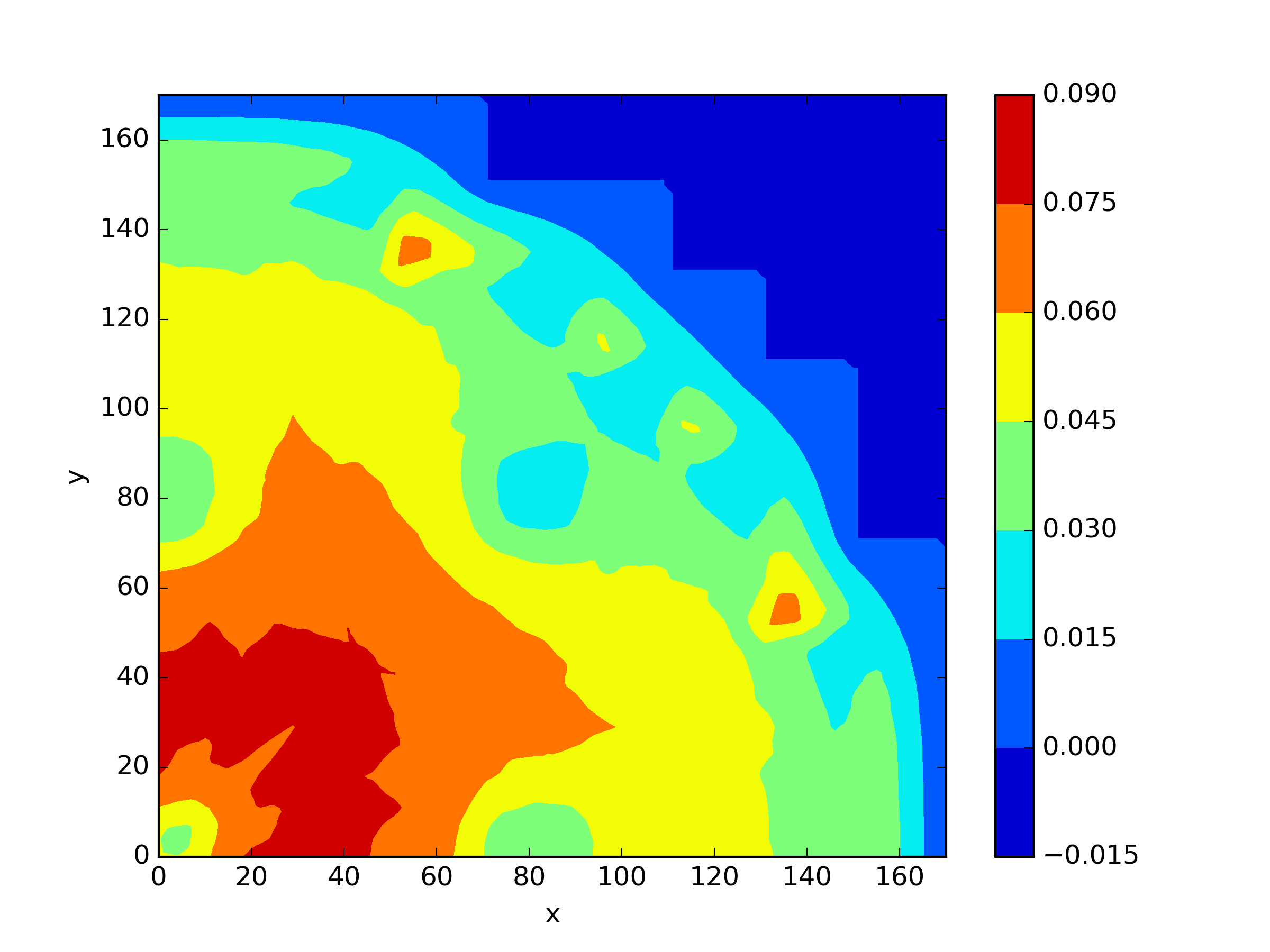}
		\end{minipage}
		\begin{minipage}{0.24\textwidth}
			\centering
			\includegraphics[width=4.25cm]{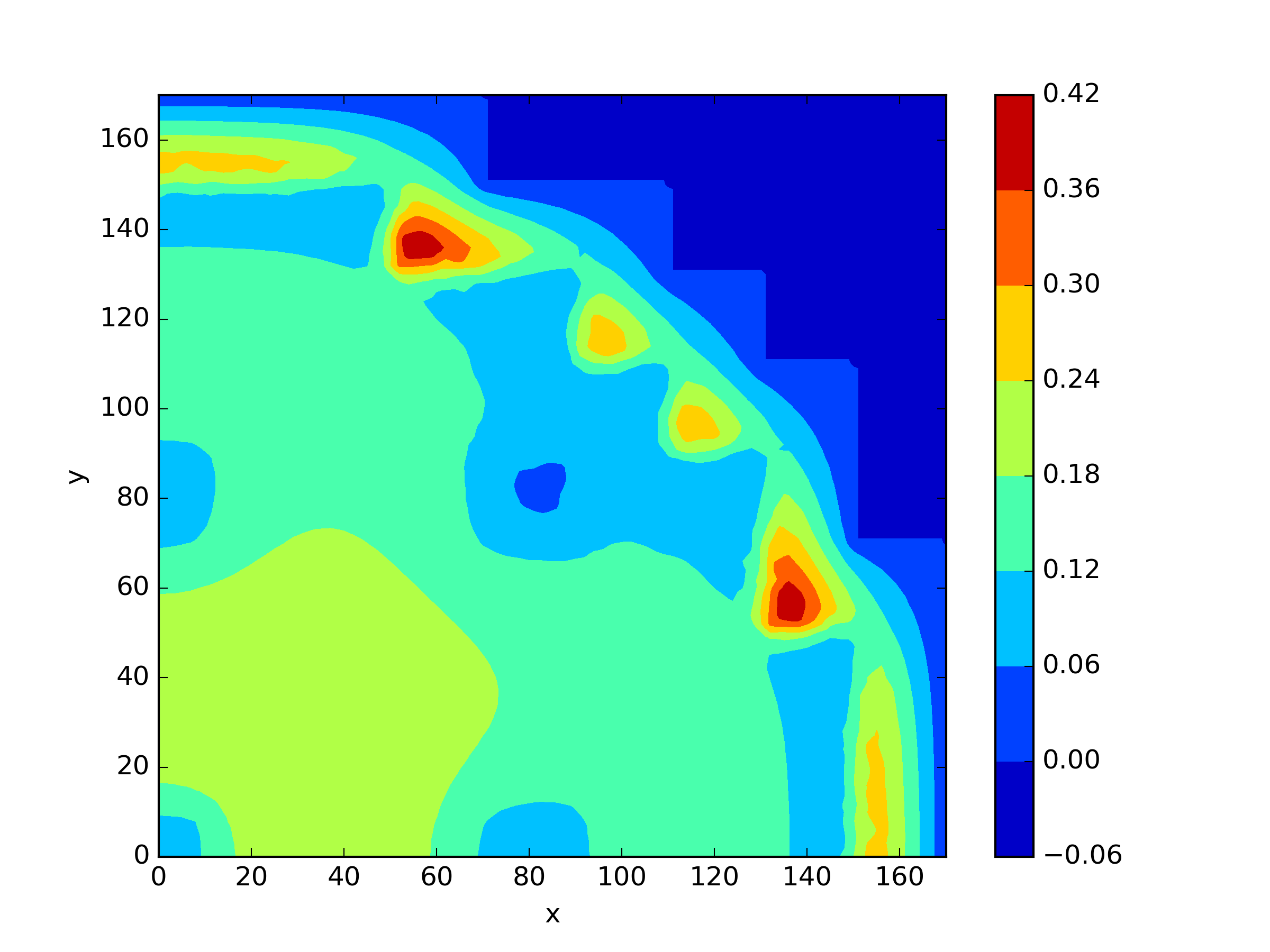}
		\end{minipage}
		\begin{minipage}{0.24\textwidth}
			\centering
			\includegraphics[width=4.25cm]{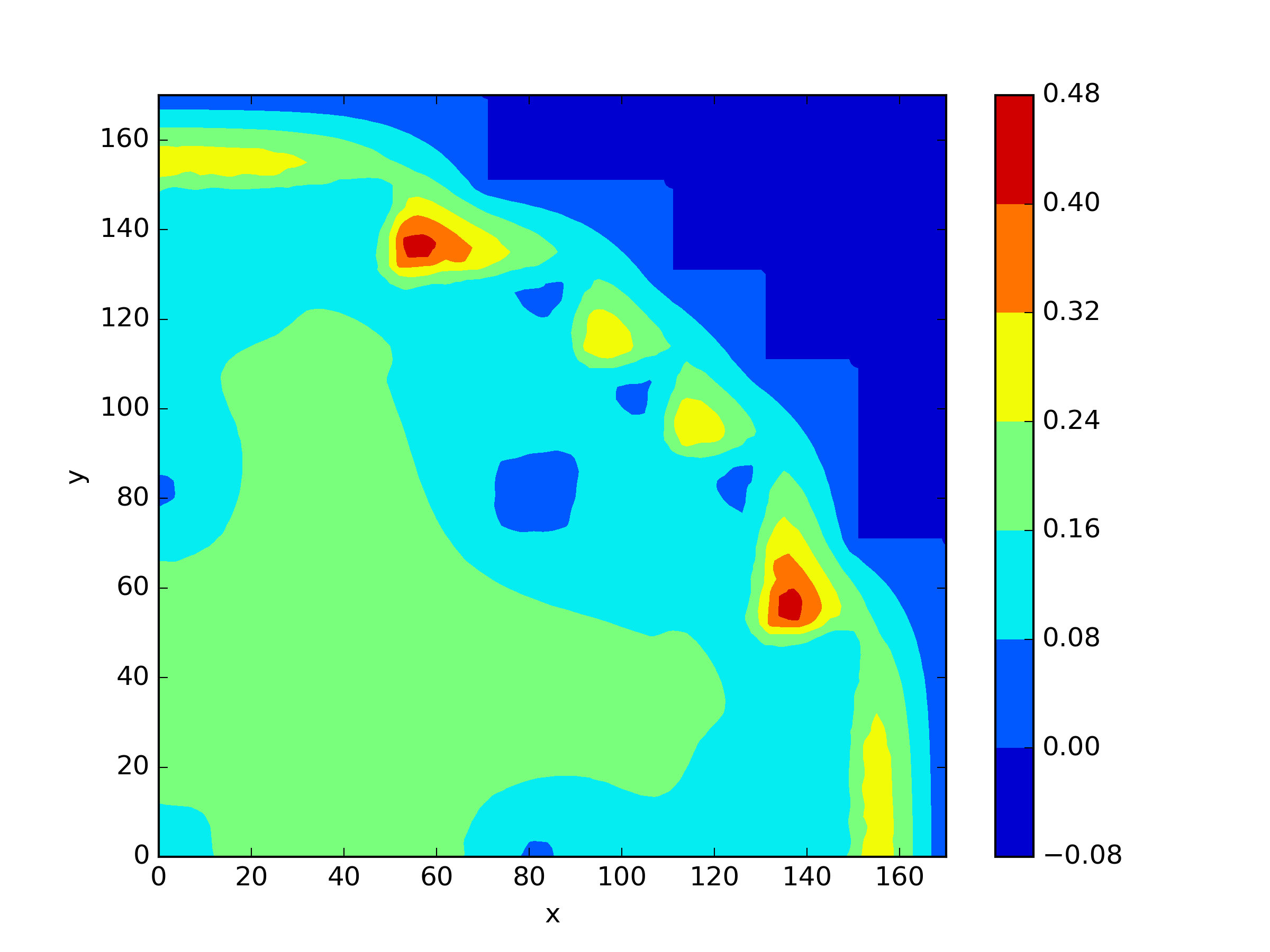}
		\end{minipage}
		\begin{minipage}{0.24\textwidth}
			\centering
			\includegraphics[width=4.25cm]{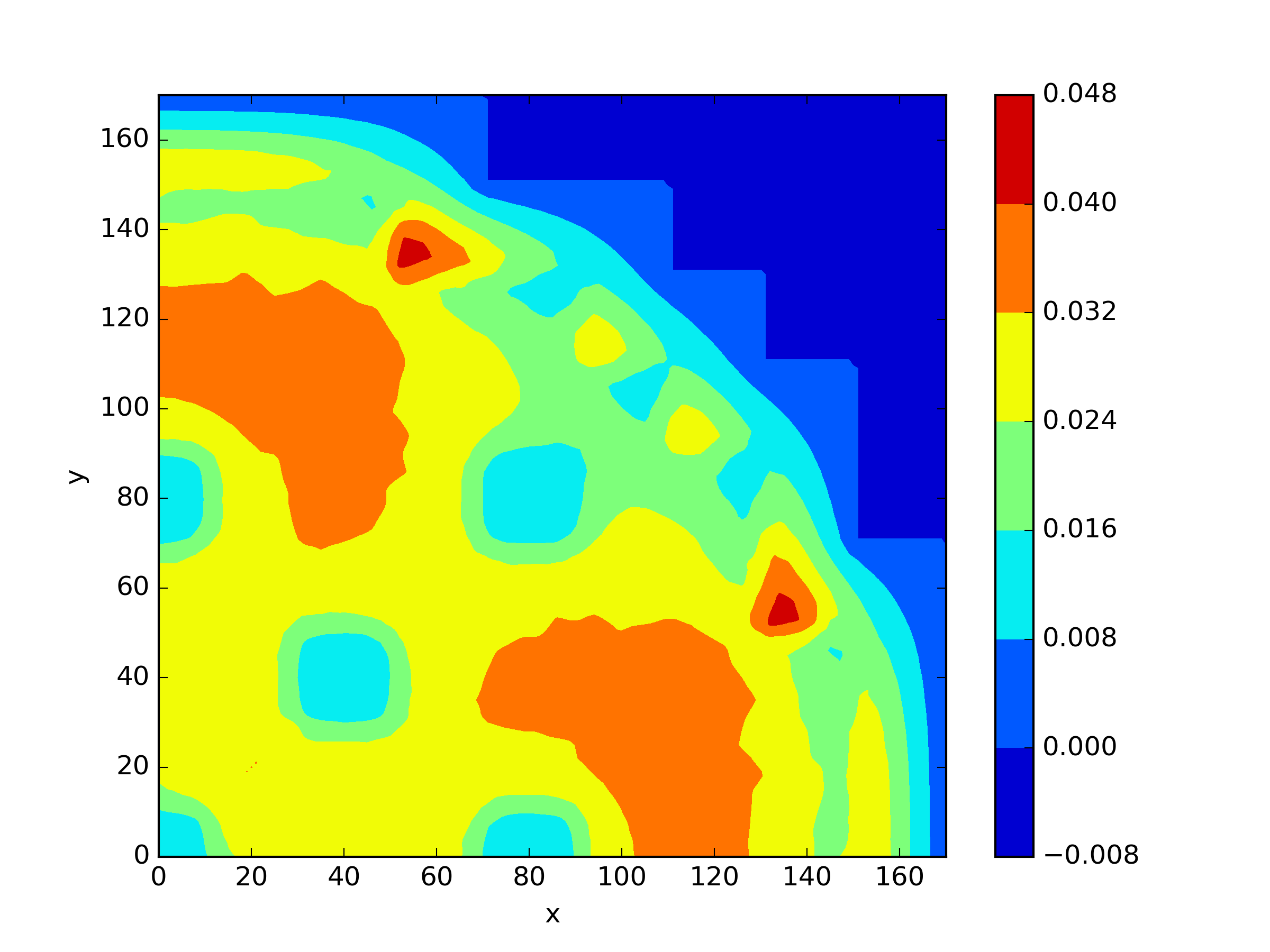}
		\end{minipage}
		
		\begin{minipage}{0.24\textwidth}
			\centering
			\includegraphics[width=4.25cm]{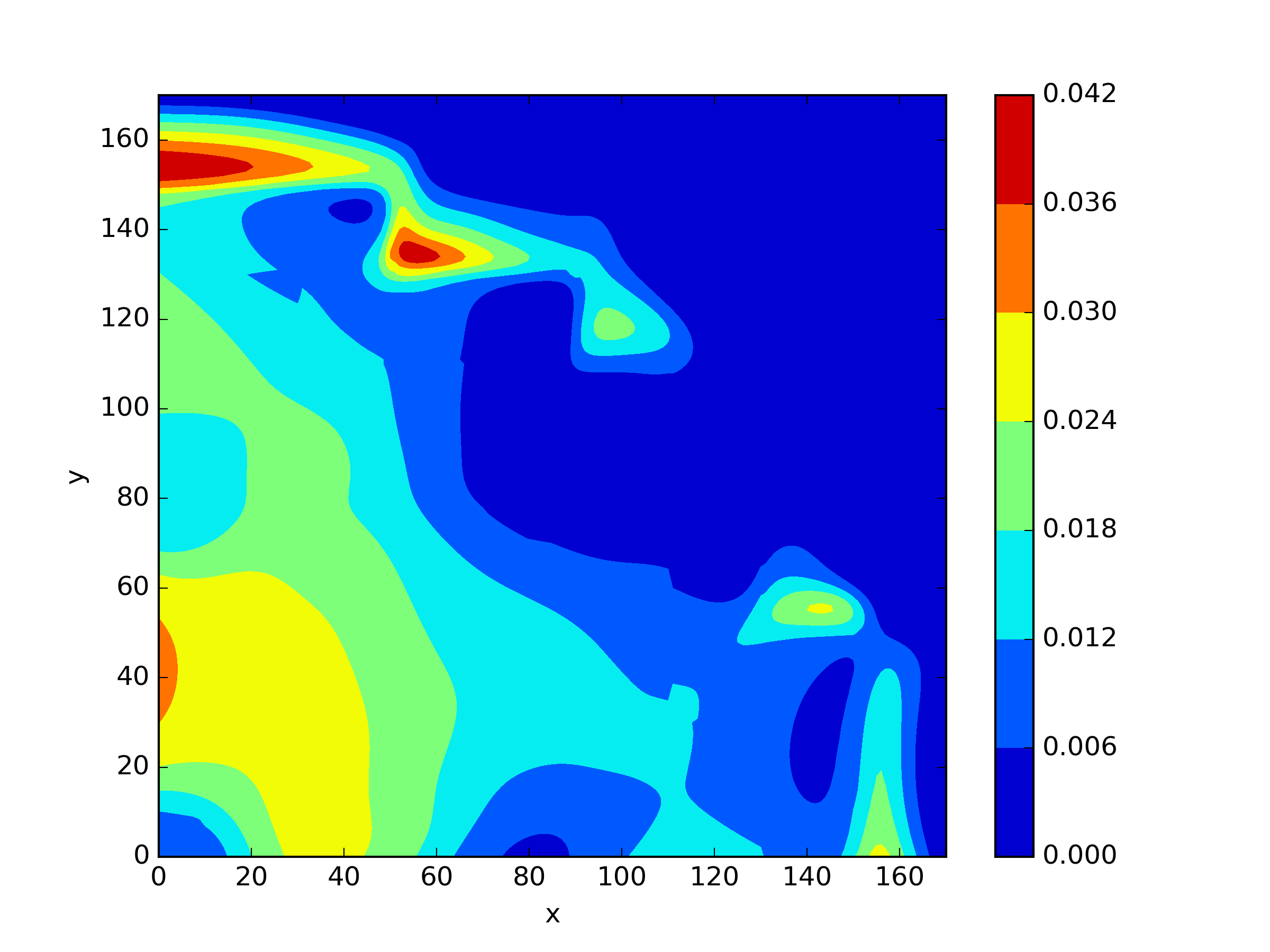}
		\end{minipage}
		\begin{minipage}{0.24\textwidth}
			\centering
			\includegraphics[width=4.25cm]{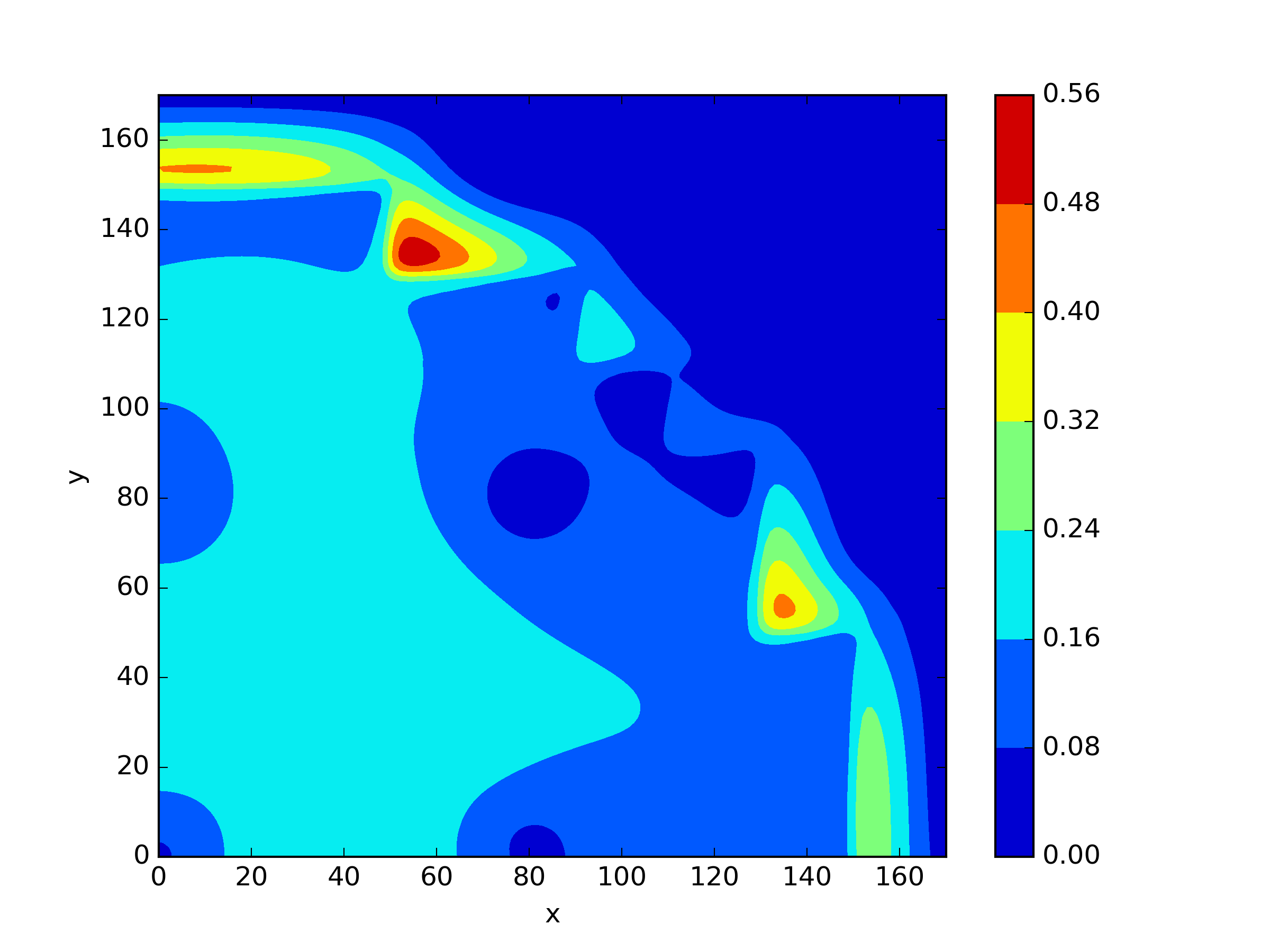}
		\end{minipage}
		\begin{minipage}{0.24\textwidth}
			\centering
			\includegraphics[width=4.25cm]{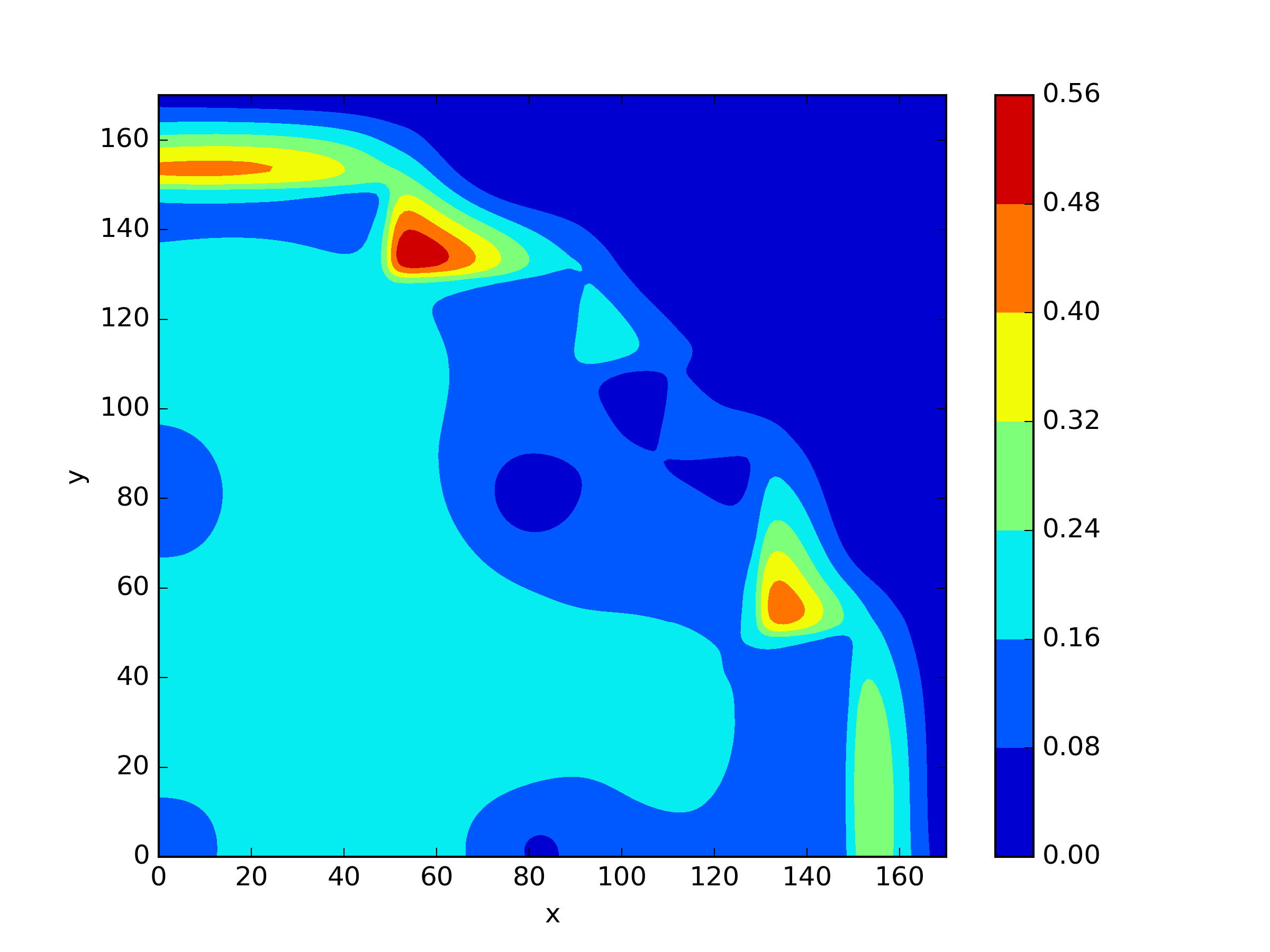}
		\end{minipage}
		\begin{minipage}{0.24\textwidth}
			\centering
			\includegraphics[width=4.25cm]{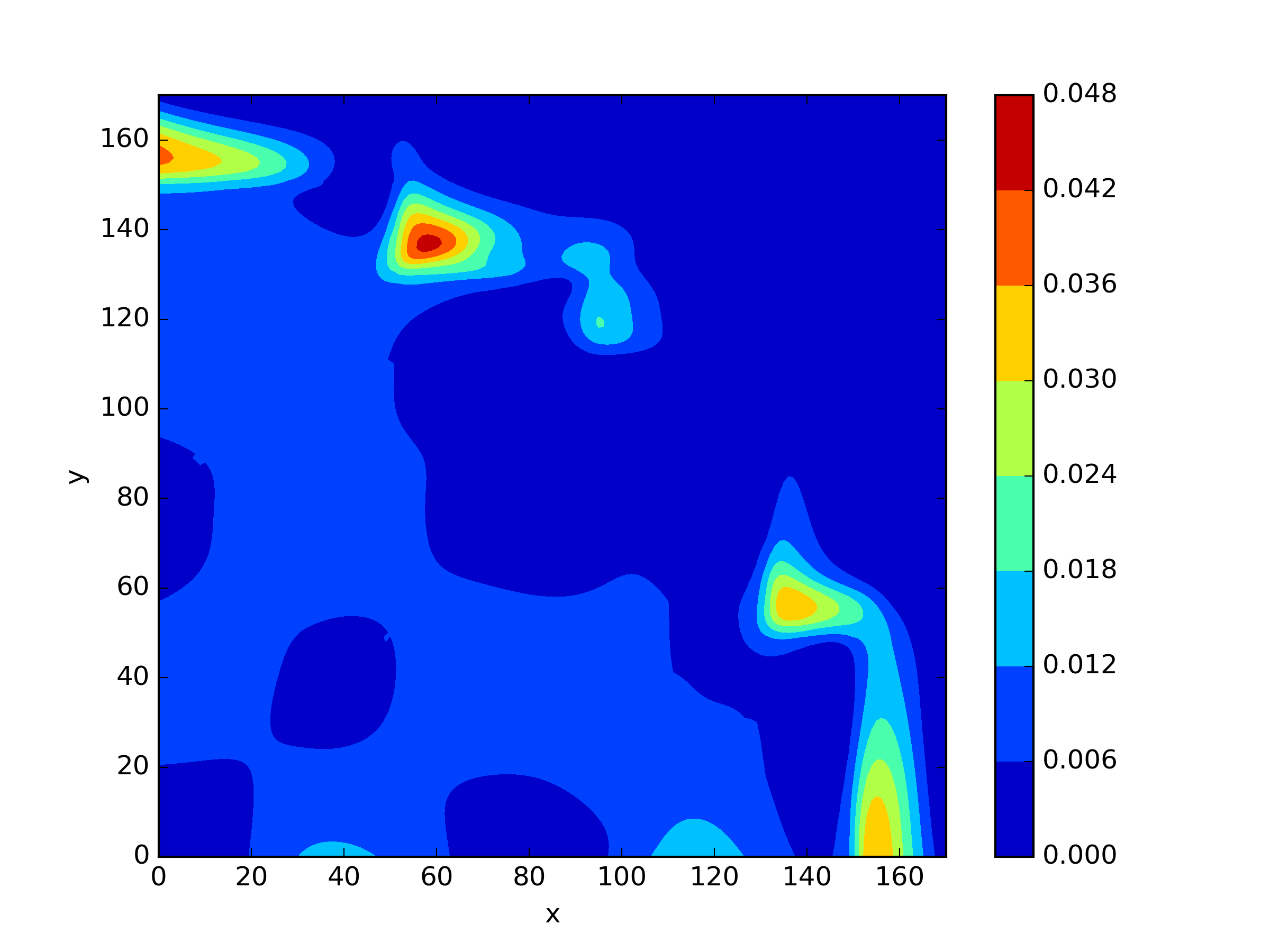}
		\end{minipage}
		
		\begin{minipage}{0.24\textwidth}
			\centering
			\includegraphics[width=4.25cm]{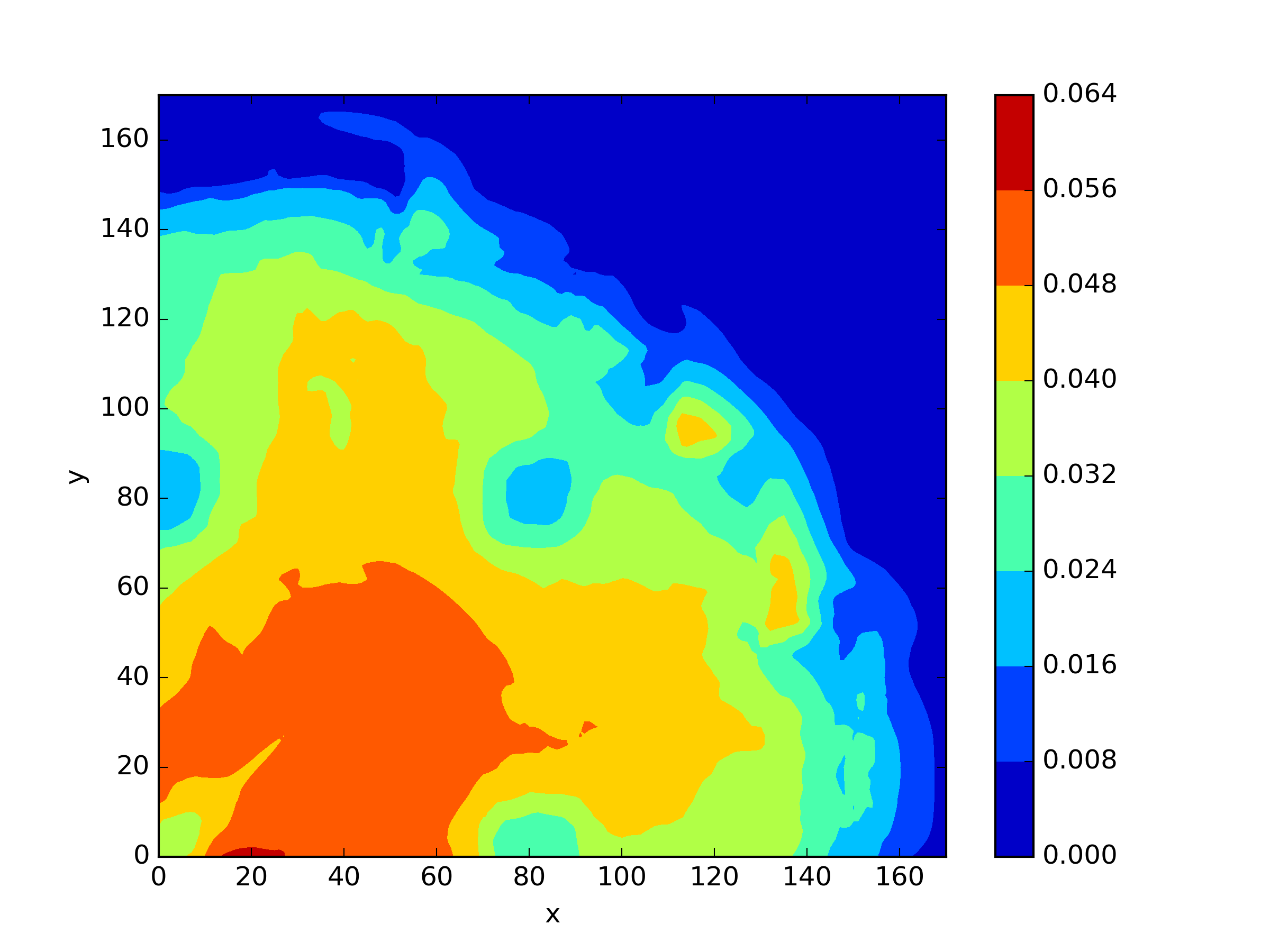}
		\end{minipage}
		\begin{minipage}{0.24\textwidth}
			\centering
			\includegraphics[width=4.25cm]{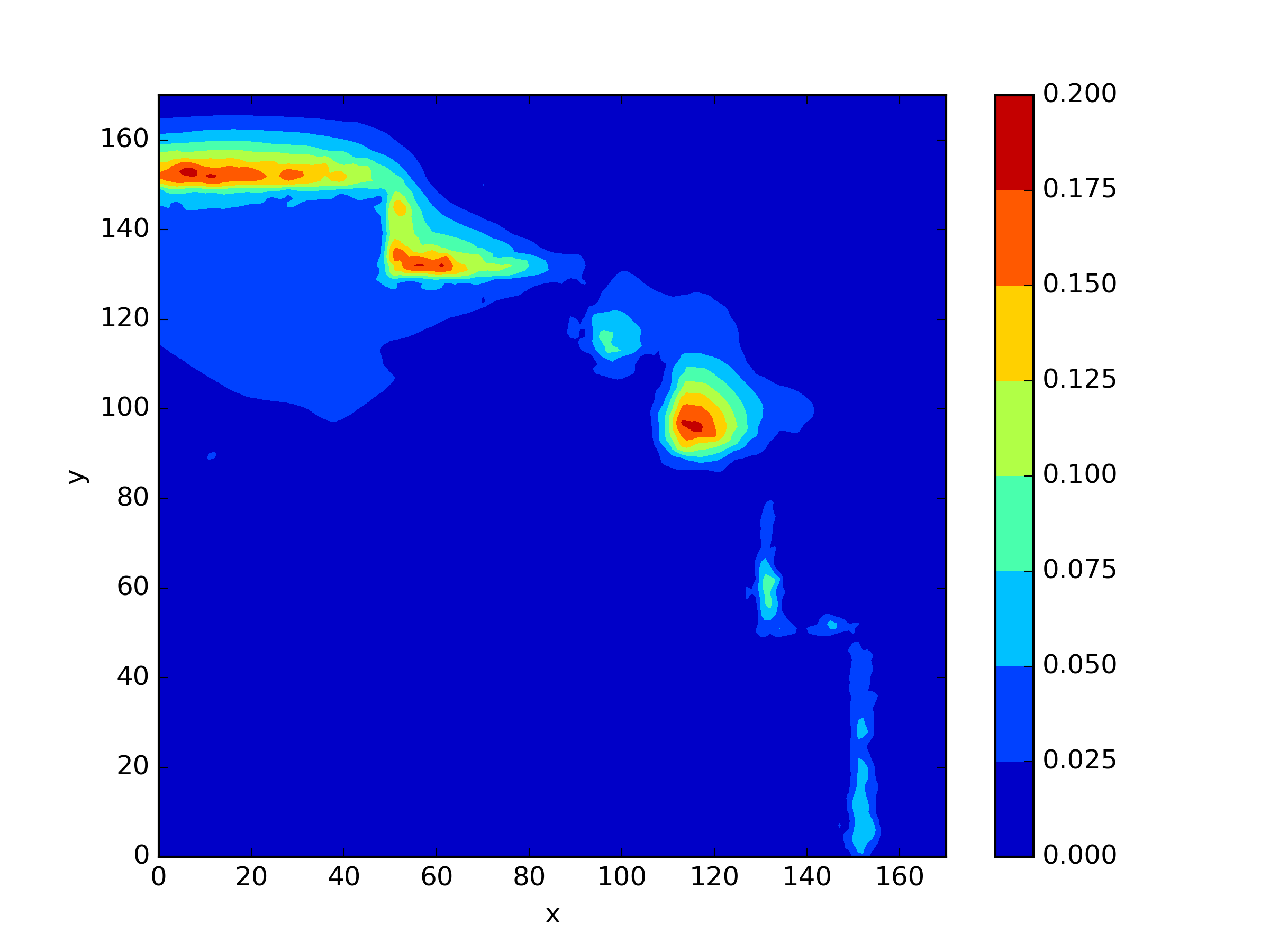}
		\end{minipage}
		\begin{minipage}{0.24\textwidth}
			\centering
			\includegraphics[width=4.25cm]{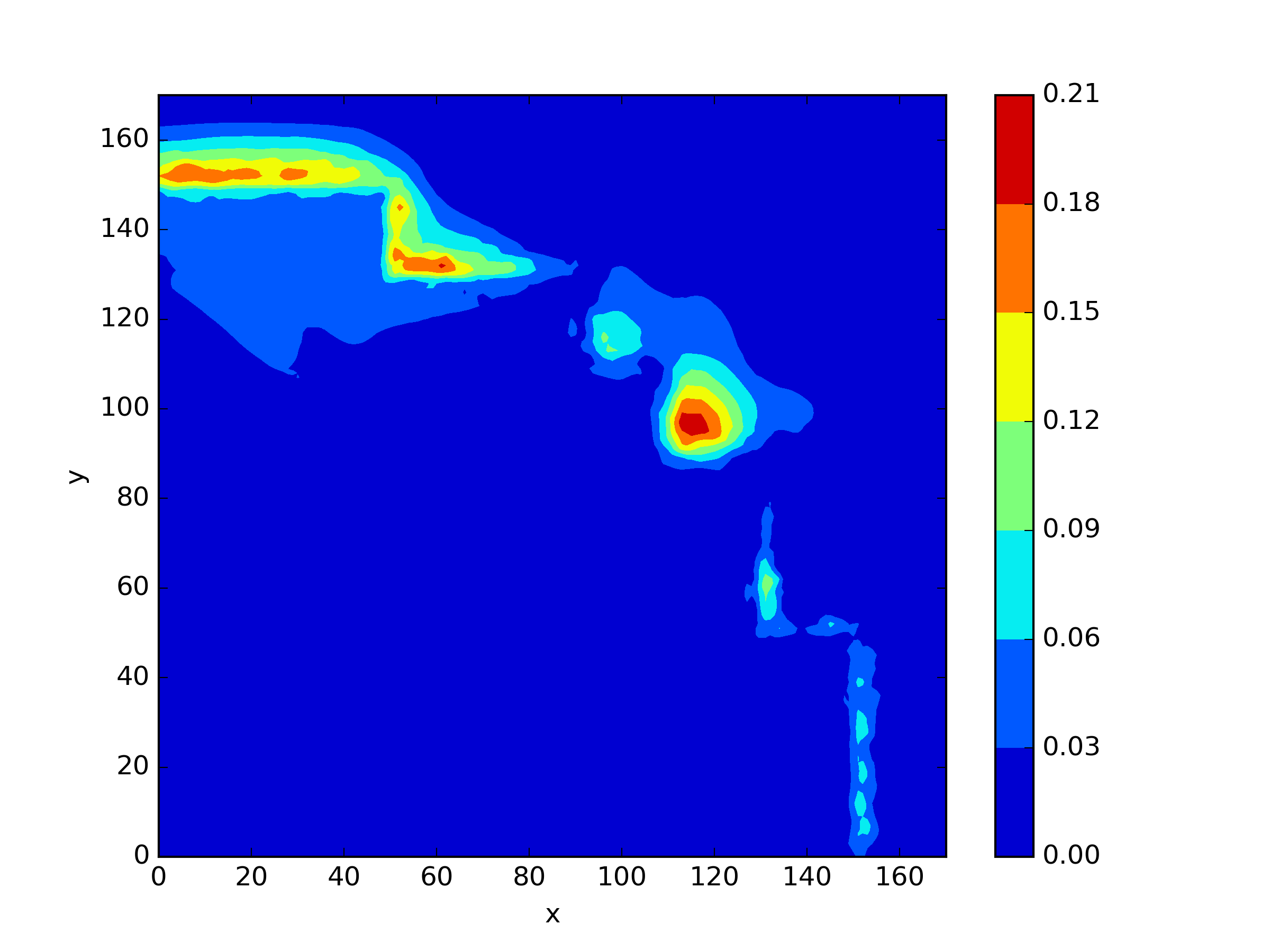}
		\end{minipage}
		\begin{minipage}{0.24\textwidth}
			\centering
			\includegraphics[width=4.25cm]{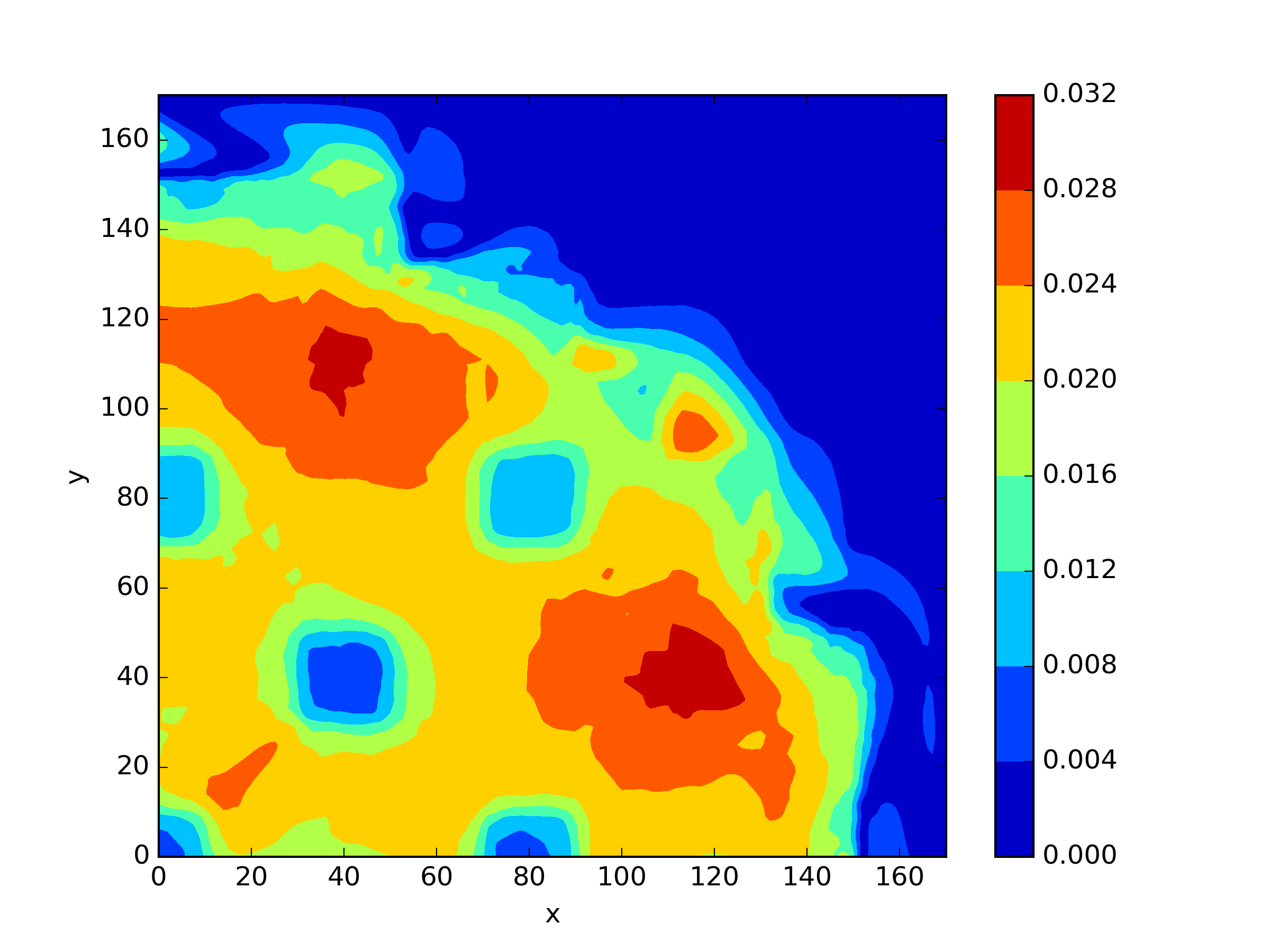}
		\end{minipage}
		\caption{The first column, second column, third column and fourth column in the figure depict the results of the solution at z = 20, 130, 190 and 360, respectively. The six rows of images, from top to bottom, represent the reference solution of $\phi_1$, the neural network solution of $\phi_1$, the absolute error between the two, the reference solution of $\phi_2$, the neural network solution of $\phi_2$, and the absolute error between the two.}
		\label{pic-IAEA-3D}
	\end{figure*}
	
	The numerical results Table \ref{tab-3D} Figure \ref{pic-TWIGL-3D} and Figure \ref{pic-IAEA-3D} indicate that our neural network algorithm has captured most of features of the 3D case in both the TWIGL and IAEA problems. The significant infinite error in $\phi_2$ in the IAEA problem is primarily due to the sharp variations in the solution within the reflector region, where $\nu\Sigma_{f,2}$ vanishes resulting in disappearance of right hand side of \eqref{eq-diff-eigen}. In contrast, the solution within the core region is relatively easier to obtain, as evidenced by the comparison of the solutions between IAEA and TWIGL. There is still room for improvement, as we are in the face of realistic problems.
	
	\section{Conclusion} \label{Conclusion}
	The main contribution of this paper are the application of neural networks to solve multi-group neutron diffusion eigenvalue problem in nuclear reactor physics and proposal of two new residual loss function. We have successfully solved multiple physical problems in one-dimensional and two-dimensional within an acceptable range of error, even with random initial guess of eigenvalue and eigenfunctions. Two new loss functions, decoupling residual loss function \eqref{loss-res-3} and direct iterative residual loss function \eqref{loss-res-2}, are proposed and compared with the previous PC-GIPMNN method for two-group neutron diffusion eigenvalue problem. 
	
	Numerous numerical experiments based on benchmark problem in nuclear reactor physics reveal the following points:
	\begin{itemize}
		\item For the same problem, increasing the number of sampling points leads to higher accuracy. However, for a fixed neural network structure, the accuracy will not infinitely improve. Nevertheless, increasing the number of sampling points can accelerate the convergence of the solution.
		\item For the same interface eigenvalue problem, the introduction of interface conditions is crucial for solving the problem itself. The accurate modeling of interface conditions significantly impacts the quality of the solution.
		\item As the geometric complexity and scale of the solved problem increase, the difficulty of the solution also rises. However, neural networks still possess the capability to capture high-precision solutions even in challenging scenarios.
	\end{itemize}
	
	These findings highlight the importance of sampling points, interface conditions, and the ability of neural networks when combined with decoupling residual loss function to handle complex and large-scale problems while achieving acceptable accurate solutions. 
	
	\begin{table*}[h]
		\caption{Mean value of time of training every 10 epochs for different problems.}
		\centering
		\begin{tabular}{cc}
			\hline
			Case & Time \\
			\hline
			One-dimensional problem     & 0.4644 s\\
			Two-dimensional TWIGL problem    & 1.6541 s\\
			Two-dimensional TWIGL-3R problem & 2.0659 s\\
			Two-dimensional IAEA problem     & 2.4642 s\\
			Three-dimensional TWIGL problem    & 24.1786 s\\
			Three-dimensional IAEA problem     & 293.1382 s\\
			\hline
		\end{tabular}
		\label{tab-time}
	\end{table*}
	
	On the other hands, the ability of neural network to solve complex three-dimensional problems with high accuracy can be challenging due to the increased computational complexity. More sophisticated approaches and large training dataset may be required to improve the accuracy of the neural network solution for three-dimensional problems. We have observed that there is still a contradiction between a large number of sampling points and minimizing the solution time (Table \ref{tab-time}) within the framework of this paper. In the case of three-dimensional problems, the accuracy achieved by the neural network is not yet satisfactory. 
	
	Many researchers have discovered challenges and confusions also existing in using neural networks to solve PDEs during their research process \cite{grossmann2023can} \cite{krishnapriyan2021characterizing}. {At the same time, maintaining solution accuracy as the problem becomes larger and more complex is a goal that many researchers are actively pursuing, as referenced in \cite{Tangkejun2023jcp} and \cite{Yangyu2024movingsampling}.} 
	
	{It is interesting to further improve the accuracy of decoupling residual loss function and direct iterative residual loss function applied to this kind of problem. The experiments in this paper show that it is crucial to combine appropriate physical information as a loss function.}
	
	The Boltzmann transport equation in nuclear reactor physics involves the transport and interactions of neutrons in reactor materials, considering position, direction, energy, and time. We conduct a series of research on issues related to neutron diffusion models {which is derived from Boltzmann transport equation under simplified assumptions} with the aim of exploring the efficient integration of deep learning with precise nuclear reactor physics models. This, in turn, enables the development of deep learning frameworks for neutron transport coupled with thermal-hydraulics, structural mechanics, and other physical fields.
	
	Hence, for our future work:
	\begin{itemize}
		\item{Hoping to have a deeper understanding of the physical properties of the transport equation and deriving a new type of loss function which could significantly overcome the obstacle caused by larger and more complex geometry domain.}
		
		\item {Exploring techniques to strike a balance between the number of training points, time of solving problems, and solution accuracy.}
		
		\item {Paying attention to the problem of multi-physics coupling solved by neural network, hoping to overcome the problems encountered in this process by hybrid model.}
	\end{itemize} 
	{More than these, we also remain attentive to the research of sampling strategy and other neural network theory.}
	
	\section*{Acknowledgments}
	This research was partially supported by the Natural Science Foundation of Shanghai (No.23ZR1429300), Innovation Funds of CNNC (Lingchuang Fund, Contract No.CNNC-LCKY-202234), the Project of the Innovation Center of Science Technology and Industry(Contract No.HDLCXZX-2023-HD-039-02), Natural Science Foundation of Sichuan Province, China (2023NSFSC0075) , the National key R \& D Program of China (No.2022YFE03040002) and the National Natural Science Foundation of China (No.11971020, No.12371434).

\end{document}